\documentclass[11pt,twoside, a4paper]{article}

\usepackage{latexsym}
\usepackage[english]{babel}
\usepackage{t1enc}
\usepackage{mathrsfs}
\usepackage{ifthen}
\usepackage{url}
\usepackage{epsfig}
\usepackage[nottoc,notlof,notlot]{tocbibind}

\usepackage{amsbsy}
\usepackage{extarrows}
\usepackage{amsfonts}
\usepackage{amsmath}
\usepackage{bbm}
\usepackage{amssymb}
\usepackage{amsthm}
\usepackage{amsxtra}
\usepackage{bbm}
\usepackage{stmaryrd}
\usepackage{enumitem}

\usepackage[margin=2cm]{geometry}

\usepackage{verbatim}

\usepackage[utf8x]{inputenc}
\usepackage[T1]{fontenc}
\usepackage{enumitem} 

\usepackage{hyperref} 

\numberwithin{equation}{section}
\numberwithin{figure}{section}

\newtheoremstyle{thm-style-oskari}
{7pt}      
{7pt}      
{\itshape} 
{}         
{\scshape} 
{.}        
{.5em}     
{}         

\theoremstyle{thm-style-oskari}
    \newtheorem{theorem}{Theorem}[section]
    \newtheorem{proposition}[theorem]{Proposition}
    \newtheorem{corollary}[theorem]{Corollary}
    \newtheorem{lemma}[theorem]{Lemma}
    \newtheorem{definition}[theorem]{Definition}
    \newtheorem{example}[theorem]{Example}
    \newtheorem{convention}[theorem]{Convention}
    \newtheorem{remark}[theorem]{Remark}

\newenvironment{Proof}[1][Proof]{\begin{proof}[\sc{#1}]}{\end{proof}}

\newcommand{\bels}[2] {
        \begin{equation} \label{#1} \begin{split} 
                #2 
        \end{split} \end{equation}
        }
\newcommand{\bes}[1]{
        \begin{equation*}  \begin{split} 
                #1 
        \end{split} \end{equation*}
        }

\newcommand{\bs}[1]{\boldsymbol{\mathrm{#1}}} 
\newcommand{\bb}{\mathbb} 
\newcommand{\bbm}{\mathbbm} 
\renewcommand{\rm}{\mathrm} 
\renewcommand{\cal}{\mathcal} 
\newcommand{\scr}{\mathscr} 
\renewcommand{\frak}{\mathfrak} 
\newcommand{\ul}[1]{\underline{#1} \!\,} 
\newcommand{\ol}[1]{\overline{#1} \!\,} 
\newcommand{\wh}{\widehat}
\newcommand{\wt}{\widetilde}

\newcommand{\eps}{\varepsilon}
\newcommand{\ord} {\mathcal{O}}

\newcommand{\pt}{\partial} 

\renewcommand{\P}{\mathbb{P}}
\newcommand{\E}{\mathbb{E}}
\newcommand{\R}{\mathbb{R}}
\newcommand{\C}{\mathbb{C}}
\newcommand{\N}{\mathbb{N}}

\newcommand{\D}{\mathbb{D}}

\newcommand{\ee}{\mathrm{e}} 
\newcommand{\ii}{\mathrm{i}} 
\newcommand{\dd}{\mathrm{d}}

\newcommand{\p}[1]{({#1})}
\newcommand{\pb}[1]{\bigl({#1}\bigr)}
\newcommand{\pB}[1]{\Bigl({#1}\Bigr)}
\newcommand{\pbb}[1]{\biggl({#1}\biggr)}

\renewcommand{\sb}[1]{\bigl[{#1}\bigr]}
\newcommand{\sB}[1]{\Bigl[{#1}\Bigr]}
\newcommand{\sbb}[1]{\biggl[{#1}\biggr]}

\newcommand{\cb}[1]{\bigl\{{#1}\bigr\}}
\newcommand{\cB}[1]{\Bigl\{{#1}\Bigr\}}
\newcommand{\cbb}[1]{\biggl\{{#1}\biggr\}}

\newcommand{\abs}[1]{\lvert #1 \rvert}
\newcommand{\absb}[1]{\big\lvert #1 \big\rvert}
\newcommand{\absB}[1]{\Big\lvert #1 \Big\rvert}
\newcommand{\absbb}[1]{\bigg\lvert #1 \bigg\rvert}

\newcommand{\norm}[1]{\lVert #1 \rVert}
\newcommand{\normb}[1]{\big\lVert #1 \big\rVert}
\newcommand{\normB}[1]{\Big\lVert #1 \Big\rVert}
\newcommand{\normbb}[1]{\bigg\lVert #1 \bigg\rVert}

\newcommand{\tnorms}[2][0]{%
  \ifcase#1\relax
    \lvert\kern-0.25ex\lvert\kern-0.25ex\lvert #2 \rvert\kern-0.25ex\rvert\kern-0.25ex\rvert\or 
    \big\lvert\kern-0.25ex\big\lvert\kern-0.25ex\big\lvert #2 \big\rvert\kern-0.25ex\big\rvert\kern-0.25ex\big\rvert\or   
    \Big\lvert\kern-0.25ex\Big\lvert\kern-0.25ex\Big\lvert #2 \Big\rvert\kern-0.25ex\Big\rvert\kern-0.25ex\Big\rvert\or   
    \bigg\lvert\kern-0.25ex\bigg\lvert\kern-0.25ex\bigg\lvert #2 \bigg\rvert\kern-0.25ex\bigg\rvert\kern-0.25ex\bigg\rvert\or   
    \Bigg\lvert\kern-0.25ex\Bigg\lvert\kern-0.25ex\Bigg\lvert #2 \Bigg\rvert\kern-0.25ex\Bigg\rvert\kern-0.25ex\Bigg\rvert\or   
    \left\lvert\kern-0.25ex\left\lvert\kern-0.25ex\left\lvert #2 \right\rvert\kern-0.25ex\right\rvert\kern-0.25ex\right\rvert
  \fi
} 

\newcommand{\nnorm}[1]{{\vert\kern-0.25ex\vert\kern-0.25ex\vert #1 
    \vert\kern-0.25ex\vert\kern-0.25ex\vert}}

\newcommand{\avg}[1]{\langle #1 \rangle}
\newcommand{\avgb}[1]{\big\langle #1 \big\rangle}
\newcommand{\avgB}[1]{\Big\langle #1 \Big\rangle}
\newcommand{\avgbb}[1]{\bigg\langle #1 \bigg\rangle}

\newcommand{\scalar}[2]{\langle{#1} \mspace{2mu}, {#2}\rangle}

\newcommand{\scalarbb}[2]{\bigg\langle{#1} \,\mspace{2mu},\, {#2}\bigg\rangle}

\DeclareMathOperator{\tr}{Tr}
\DeclareMathOperator{\var}{Var}
\DeclareMathOperator{\supp}{supp}

\DeclareMathOperator{\re}{Re}
\DeclareMathOperator{\im}{Im}

\DeclareMathOperator{\dist} {dist}                
\DeclareMathOperator*{\spec}{Spec}						
\renewcommand{\Re}{\re}
\renewcommand{\Im}{\im} 

\newcommand{\1} {\mspace{1 mu}}
\newcommand{\2} {\mspace{2 mu}}
\newcommand{\msp}[1] {\mspace{#1 mu}}

\newcommand{\rhoDOS}{\rho}
\newcommand{\spradius}{\varrho}

\newcommand{\mtwo}[2]
{
\left(
\begin{array}{cc}
#1 
\\
#2
\end{array}
\right)
}

\newcommand{\vtwo}[2]
{
\left(
\begin{array}{c}
#1 
\\
#2
\end{array}
\right)
}

\newcommand{\qq}[1]{\llbracket #1 \rrbracket}
\newcommand{\smin}{\mathrm{s}_{\min}} 
\DeclareMathOperator{\Span}{span}

\makeatletter
\def\blfootnote{\xdef\@thefnmark{}\@footnotetext}
\makeatother

\begin{document}
\blfootnote{Date: \today} 
\blfootnote{MSC2010 Subject Classifications: 60B20, 15B52, 46Txx.} 

\title{Inhomogeneous Circular Law for Correlated Matrices} 
\author{ 
{\small \begin{tabular}{c} 
\large Johannes Alt\footnote{
{Partial funding from the European Research Council (ERC) under the European Union's Horizon 2020 research and innovation programme (grant agreement No. 715539 RandMat) and from the Swiss National Science Foundation through the NCCR SwissMAP grant is gratefully acknowledged.}
Email: \href{mailto:johannes.alt@unige.ch}{johannes.alt@unige.ch}} \\ 
\small  University of Geneva  \end{tabular}} 
\hspace*{0.5cm} \and \hspace*{0.5cm} 
{\small \begin{tabular}{c} 
\large Torben Krüger\footnote{Partially supported by the Hausdorff Center for Mathematics. Email: \href{mailto:torben.krueger@uni-bonn.de}{torben.krueger@uni-bonn.de}}\\
\small University of Bonn \end{tabular}}
}

\date{} 

\maketitle
\thispagestyle{empty}

\vspace*{-0.8cm} 

\begin{abstract}
We consider non-Hermitian random matrices $X \in \C^{n \times n}$ with general decaying correlations between their entries. 
For large $n$, the empirical spectral distribution is  
 well approximated by a deterministic density, expressed in terms of the solution to a system of two coupled non-linear $n \times n$ matrix equations. 
 This density is interpreted 
as the Brown measure of a linear combination of free circular elements with matrix coefficients on a non-commutative probability space. 
It is radially symmetric, real analytic in the radial variable 
and strictly positive on a disk around the origin in the complex plane with a  discontinuous  drop to zero at the edge. 
 The radius of the disk is given explicitly in terms of the covariances of the entries of $X$. 
We show convergence down to local spectral scales just slightly above the typical eigenvalue spacing with an optimal rate of convergence. 
\end{abstract}

\tableofcontents

\section{Introduction}

Many random matrix models exhibit a strong concentration of measure phenomenon; their empirical eigenvalue distributions  are well approximated by deterministic measures as their sizes tend to infinity.  For Hermitian matrices, the simplest and most prominent example is the celebrated semicircle law for Wigner ensembles with independent and identically distributed (i.i.d.) entries above the diagonal \cite{Wigner1955}. 
Girko's circular law \cite{Girko1985,bai1997}  is its non-Hermitian analogue\footnote{We refer to the survey \cite{BordenaveChafai2012} for a complete account of the history of the circular law until the minimal moment assumptions in 
\cite{tao2010}.}. For  matrices $X=(x_{ij})_{i,j =1}^n$ with centred i.i.d.\ entries, unrestricted by symmetry and with  normalisation $\E\2\abs{x_{ij}}^2=\frac{1}{n}$, it asserts convergence of the eigenvalue distribution to the uniform probability measure on the unit disk in the complex plane.

Establishing similar concentration results and identifying the limiting spectral density  while simultaneously relaxing the two basic assumptions of identical distributions and independence of the entries has since been the focus of many works in random matrix theory. When the entries are independently drawn from different distributions, their variance profile $s_{ij}=\E\2\abs{x_{ij}}^2$ becomes an additional parameter of the model that determines the density through the nonlinear Dyson equation for $n$ unknowns. Since in general no explicit formula for its solution is available, analysing the characteristic properties of the spectral density has attracted considerable attention.

In the Hermitian case, convergence of the empirical spectral measure is well established \cite{AZind,Guionnet-GaussBand,ShlyakhtenkoGBM}  and a classification of the degree of regularity of the asymptotic density as well as  of its possible singularities has been given \cite{AjankiSingularities}.   Even when the independence of matrix entries is dropped and local correlations with sufficient decay are considered this classification persists \cite{Altshape} and concentration of the spectral measure has been proven in broad generality \cite{MR2417889,MR3332852,MR1431189,Erdos2017Correlated,MR1887675,MR2191967,PasturShcerbinaAMSbook,MR2444540,MR2155229}.

There are far fewer results on the  existence   and characteristics of limiting spectral densities  
for non-Hermitian matrices since their  spectral instability makes such questions more challenging compared to the Hermitian situation. 
For random matrices $X$ with centred, independent entries and a general variance profile, 
the convergence of the spectral measure of $X$ to a rotationally symmetric, continuous limiting 
density $\sigma$ was shown in \cite{Cook2018}, and independently in \cite{Altcirc} on all mesoscopic scales in the bulk spectrum under stronger assumptions on the variance profile and regularity of the entry distribution.   
The extension of convergence on mesoscopic scales  to the spectral edges  and optimal control of the spectral radius  was achieved in \cite{AltSpecRad}. 
These three papers avoided the requirement of identical variances imposed earlier.

In the present paper we also depart from the independence assumption on the entries. 
We consider a large class of  centred  non-Hermitian random matrices $X\in \C^{n\times n}$  with general decaying correlations among their entries. 
Throughout this class, the limiting spectral density $\sigma$ is determined solely by  the covariances between the matrix entries and has the following properties:
(i) the density is rotationally symmetric around zero, (ii) its support is a disk centred at 
the origin, 
(iii) the density is real analytic  as a function of the radial variable  inside the disk and has a jump at its boundary.  

The analyticity is a new result even when the entries of $X$ are independent  (apart from the explicitly known circular law case).   
In this case, the other properties are known \cite{Altcirc}. 
 Remarkably, the support of $\sigma$ is always connected in the 
non-Hermitian case, in the independent as well as the correlated setup. 
This is in sharp contrast to the Hermitian case, where the support can be disconnected even for matrices 
with centred, independent entries and a variance profile \cite{AjankiQVE}. 

 The class of random matrices $X$ we consider here contains finite-dimensional approximations  of linear combinations of free circular elements with matrix coefficients on a non-commutative probability space. These are non-normal analogs of   operator-valued semicircular elements introduced in \cite{Voiculescu1995} (see also \cite{Speicher1998}). 
For such linear combinations, one is interested in their Brown measure, a generalisation of the spectral measure of normal operators to general operators in a finite von Neumann algebra. 
It was introduced in \cite{Brown1986} and revived in \cite{HaagerupLarsen2000}. 
Since then  significant attention  has been given to  determining the Brown measure and  understanding its  properties  for specific classes of non-normal operators,  see e.g.\ 
\cite{BelinschiSniadySpeicher2018,BianeLehner2001,DHKBrownMeasure,GWZBrownMeasure,HaagerupLarsen2000,Larsenthesis}. 
In the present work, we prove that the Brown measure of  these  matrix-valued circular elements 
has the properties (i), (ii), (iii) listed above. 
In previous works addition of or multiplication with an $R$-diagonal element (see \cite{NicaSpeicher1997} for the definition)  and its invariance under unitary transformations was crucial in order to introduce  generic directionality into the model.  
In contrast our model and its ensuing analysis are generically non-isotropic due to  the matrix coefficients.

Convergence of the eigenvalue density to a limiting measure is commonly expressed by showing that 
for each ball with   fixed diameter on the scale of the entire spectrum  the fraction of eigenvalues in it agrees 
asymptotically with the mass assigned to this ball by the limiting measure. Such \emph{global law} 
is refined to a \emph{local law}, showing convergence on mesoscopic scales, by allowing the diameter to decrease with $n$ as long as it stays slightly above 
the typical eigenvalue spacing.  
We now review some previous results on local laws for non-Hermitian random matrices with independent entries.  
A bulk local law for random matrices with centred, independent entries of identical variances was 
shown in \cite{Bourgade2014}. Additionally requiring the first three moments of the entry distribution to match a standard Gaussian, the local law including the edge was established in 
\cite{tao2015} and  in \cite{BYY_circular2}. 
The third moment matching condition for the edge local law was then removed in \cite{Y_circularlaw}. 

For the bulk local law, the assumption of identical variances was dropped in \cite{Altcirc}.
In this situation, the limiting density differs substantially from the circular law. 
Under weaker moment assumptions and asymptotically identical variances, a bulk local law with the 
circular law as limiting density was shown in \cite{Gotze2017}. 
In the setup of \cite{Altcirc}, the edge local law was proven in \cite{AltSpecRad}. 

 The availability of a local law has wide ranging implications for the spectral analysis of any random matrix model. In the present paper, we apply it to exclude eigenvalues away from the support of the limiting spectral density $\sigma$, i.e.\ with high probability all eigenvalues concentrate on a disk around the origin whose radius is determined by the covariances of the matrix entries. 
We also obtain the complete isotropic delocalisation of all eigenvectors  associated  to the bulk eigenvalues.
 Furthermore, local laws have been a key ingredient in the  study of  more refined eigenvalue statistics. In the non-Hermitian i.i.d.\ setup they have been crucially used in the proofs of universality of bulk and edge eigenvalues with a four moment matching condition 
in \cite{tao2015}, edge universality with two matching moments in \cite{CipolloniEdge} 
and the central limit theorem for linear statistics in \cite{CipolloniCLT,CipolloniCLTreal}.

 Non-Hermitian random matrices without any symmetry constraint also play an important role in various applications.  In particular, they are used to model  connectivities in food webs and neural networks \cite{Allesina:2015ux,may1972will,PhysRevLett.97.188104,Sompolinsky1988}. Since understanding the stability properties of such systems requires precise knowledge of the eigenvalue locations of the associated random matrix model, our work contributes to this line of research by allowing the correlation among the connectivities to depend on underlying geometric structures.

The analysis of the eigenvalue density of a non-Hermitian random matrix $X$ 
is commonly reduced via Girko's Hermitization trick \cite{Girko1985} to the study of the family of Hermitian matrices 
\bels{Hzeta}{
\bs{H}_\zeta := \begin{pmatrix} 0 & X - \zeta \\ (X-\zeta)^* & 0 \end{pmatrix} 
}
with spectral parameter $\zeta \in \C$.
 Consequently the main task is to control the resolvent $\bs{G}(\zeta,\eta) :=(\bs{H}_\zeta - \ii \eta)^{-1}$ of $\bs{H}_\zeta$ on the imaginary axis via its deterministic approximation $\bs{M} = \bs{M}(\zeta,\eta)$ that solves 
 the associated \emph{matrix Dyson equation} (MDE)
\begin{equation} \label{eq:mde_intro} 
 -\bs{M}^{-1} = \ii \eta + \begin{pmatrix} 0 & \zeta \\ \bar \zeta & 0 \end{pmatrix} 
 + \mathscr{S}[\bs{M}].
\end{equation}
This equation  has a unique solution for every $\zeta \in \C$ and $\eta>0$ if $\Im \bs{M} = \frac{1}{2\ii} ( 
\bs{M} - \bs{M}^* )$ is required to be positive definite \cite{Helton01012007}.
Here, $\mathscr{S}$ is a linear map on $\C^{2n \times 2n}$ defined through 
\begin{equation} \label{eq:def_S_intro} 
\mathscr{S}[\bs{R}] := \begin{pmatrix} \E[X R_{22} X^*] & 0 \\ 0 & \E[X^* R_{11} X] \end{pmatrix}, \qquad \qquad \bs{R} = 
\begin{pmatrix} R_{11} & R_{12} \\ R_{21} & R_{22} \end{pmatrix} 
 \end{equation} 
for any deterministic matrix $\bs{R} \in \C^{2n\times 2n}$ with $n\times n$--blocks $R_{11}$, $R_{12}$, $R_{21}$, $R_{22}$. 
 The operator $\mathscr{S}$ captures the covariances between the entries of $X$. 

The main  tool developed in the present paper is a precise stability result for  the non-linear high dimensional matrix equation  \eqref{eq:mde_intro}. 
From \cite{AjankiCorrelated,Erdos2017Correlated}, we know that $\bs{G}=\bs{G}(\zeta,\eta)$ satisfies a perturbed version of \eqref{eq:mde_intro} with a small error term when $n$ becomes large.
Thus, $\bs{G}$ is close to $\bs{M}$ if stability of \eqref{eq:mde_intro} against small perturbations is controlled. 
Moreover, the limiting spectral density $\sigma$ for $X$ is obtained as a derivative of $\Im \bs{M}$ with respect to $\abs{\zeta}$, where $\zeta$ is the spectral parameter of $X$. 
Thus, any  analysis of $\sigma$ also requires  stability of \eqref{eq:mde_intro}. 

In previous works, the matrix structure of \eqref{eq:mde_intro} was crucially simplified 
due to more restrictive assumptions on $X$. 
If the entries of $X$ are independent, then \eqref{eq:mde_intro} reduces 
to a  vector-valued equation for the diagonal of $\bs{M}$  and, thus, the Dyson equation is formulated on the commutative algebra of diagonal matrices.  For identical variances, 
all diagonal entries of $\bs{M}$ coincide, yielding a single scalar equation. 

In the matrix setup, a general version of \eqref{eq:mde_intro} and its stability have been studied 
in \cite{AjankiCorrelated} under a strong 
irreducibility condition on $\mathscr{S}$, which is called \emph{flatness}. However, $\mathscr{S}$ as defined in \eqref{eq:def_S_intro} does not 
fulfil this flatness condition due to its special block structure, making the equation inherently unstable. 
This issue was overcome in \cite{Altcirc} for vector Dyson equation, i.e.\ when the entries of $X$ are independent. 
Owing to the commutative structure of this vector case, an additional symmetry of $\bs{M}$ 
could be exploited to obtain the stability against perturbations respecting this symmetry. 

The analysis in the present situation necessitates tackling, at the same time, both main challenges from 
\cite{AjankiCorrelated} and \cite{Altcirc}, 
the non-commutativity of the MDE and the instability due to the specific block structure, respectively. 
 The genuinely non-commutative structure of the MDE is a major obstacle throughout the entire argument 
requiring the introduction of appropriately symmetrised objects, which are much more complicated 
than their counterparts in the commutative setup. 
To resolve the instability we perform a non-linear transformation of the MDE that allows to restrict the analysis to  
the manifold of perturbations that respect the additional symmetry of $\bs{M}$.
 This transformation is also applicable in the context of other non-normal models, e.g. non-Hermitian polynomials in several non-commutative variables. Furthermore, it is crucial to show that $\sigma$ is 
a real analytic function of $\abs{\zeta}^2$.

\vspace*{0.3cm} 

\noindent \textbf{Acknowledgement:} 
The authors are grateful to L\'aszl\'o Erd{\H o}s for  inspiring and insightful discussions.

\section{Main results}
\subsection{Correlated random matrices} 

Let $X\in\C^{n \times n}$ be a random matrix with centred, $\E x_{ij}=0$, entries. For the index set we write
\[
\llbracket n\rrbracket\,:=\, \{1, \dots, n\}\,.
\]
Within our main results we will refer to the following assumptions on the entries of $X$. 
Some of them are stated in terms of the covariances between the entries of $X$. These covariances are encoded in the two operators 
 $\cal{S},\cal{S}^*: \C^{n \times n} \to \C^{n \times n} $ on the space of $n \times n$-matrices, defined through 
 \bels{definition cal S}{
\cal{S}A\,:=\, \E\2X A X^*, 
 \qquad \cal{S}^* A \, := \, \E \2X^* A X \,.
} 
\begin{enumerate}[label=\bf A\arabic*]
\item \label{assum:A1} \emph{Finite moments:} All moments  of the entries of $\sqrt{n}X$ are finite, i.e.\  there is a sequence of positive constants $C_\nu$ such that 
\bels{bounded moments}{
\E\,\abs{x_{i j}}^\nu\,\le\, C_\nu \2n^{-\nu/2}\,,
}
for all $i,j\in \llbracket n\rrbracket$ and $\nu \in \N$.
\item \label{assum:A2} \emph{Decay of correlation:}  The index set $\llbracket n\rrbracket$ is equipped with a pseudo-metric $d$ that satisfies for a fixed $p \in \N$ the sub-$p$-dimensional volume growth condition
\begin{equation} \label{eq:growth_condition} 
\abs{\{j \in \llbracket n\rrbracket: d(i,j) \le \tau\}} \le C\2\tau^p , \qquad \tau \ge 1\,, \; i \in \llbracket n\rrbracket,
\end{equation}
with a constant $C>0$.
Furthermore, the correlations among  the entries of $\sqrt{n}\2X$ decay in the product metric $d \times d$ on $\llbracket n\rrbracket^2$  faster than any power law, i.e.\ there is a sequence of positive constants $C_\nu$ such that 
\bels{Decay of correlation}{
\rm{Cov}(f_1(\sqrt{n}\2X),f_2(\sqrt{n}\2X))\,\le\, \frac{C_\nu\norm{f_1}_2 \norm{f_2}_2}{1+d\times d(\supp f_1, \supp f_2)^\nu}, \qquad \nu \in \N\,,
}
for any two measurable functions $f_i: \C^{A_i} \to \C$ with  $\norm{f_i}_2^2:=\E \abs{f_i(\sqrt{n}\2X)}^2<\infty$, where $A_i = \supp f_i \subset \llbracket n\rrbracket^2$.
\item \label{assum:A3} \emph{Flatness:}   There is a  constant $c>0$ such that for any two deterministic vectors $x, y \in \C^n$ we have 
\bels{RM Flatness}{
\E\,\abs{ \scalar{ x}{X y}}^2\,\geq\, \frac{c}{n}\, \norm{x}^2\norm{y}^2,
}
where $\scalar{\1\cdot\1}{\1\cdot\1}$ and $\norm{\2\cdot\2}$ denote the standard Euclidean scalar product and norm on $\C^n$, respectively.
\item \label{assum:A4} \emph{Smallest singular value:} 
 For each $\eps >0$ and $\nu \in \N$, there is $C_{\eps,\nu}>0$ such that 
\begin{equation} \label{eq:log_min_lambda_i_prec_1} 
 \P \big( \smin (X-\zeta) \leq \ee^{-n^\eps} \big) \leq C_{\eps,\nu} n^{-\nu} 
\end{equation}
for all $n \in \N$ and all $\zeta \in \C$.  
Here, $\smin(X-\zeta)$ denotes the smallest singular value of $X-\zeta$.  
\end{enumerate}
\begin{enumerate}[label=\bf A\arabic*']
\setcounter{enumi}{3} 
\item \label{assum:bounded_density} 
\emph{Bounded conditional density:} 
There are $q \in (1, \infty]$ and $\kappa >0$ such that, for each pair $(i,j) \in \llbracket n\rrbracket^2$, there is a probability density $\psi_{ij} \in L^{q}(\C)$  (or $\psi_{ij} \in L^{q}(\R)$ if $X \in \R^{n \times n}$) 
which satisfies 
$\E \norm{\psi_{ij}}_{q} \leq n^\kappa$ and 
\begin{equation} \label{eq:rep_density} 
\P \big( \sqrt{n} x_{ij} \in B \big| (x_{kl})_{(k,l) \in \llbracket n\rrbracket^2 \setminus \{ (i,j)\}} \big) = \int_B \psi_{ij}(z) \dd z 
\end{equation}
almost surely for all measurable $B \subset \C$ (or $B \subset \R$). 
\end{enumerate} 
\begin{enumerate}[label=\bf A\arabic*]
\setcounter{enumi}{4} 
\item \label{assum:S_inverse_bound} 
There is $c>0$ such that the spectral radius $\varrho(\cal{S})$ of $\cal{S}$ satisfies $\varrho(\cal{S})\geq c$. 
Moreover, there is an $n$-independent monotonically decreasing function $f \colon (0,\infty) \to (0,\infty)$ 
such that 
\begin{equation} \label{eq:assum_S_inverse_bound} 
\norm{(\tau - \cal{S})^{-1}}  \leq f(\tau)/\tau 
\end{equation} 
for all $\tau> \varrho(S)$ and for all $n \in \N$. 
\end{enumerate}

We remark that Assumption~\ref{assum:bounded_density} implies Assumption~\ref{assum:A4} as shown in Proposition~\ref{pro:smallest_singular_value} below. 
Moreover, in Section~\ref{sec:relaxed_assumptions}, we explain how some assumptions can 
be relaxed (see in particular Remark~\ref{rem:weaker_version_A2} for weaker versions of \ref{assum:A2}) and examples satisfying the assumptions listed above. 

The $n$-independent constants appearing in Assumptions~\ref{assum:A1}--\ref{assum:S_inverse_bound}  
 will be called \emph{model parameters} and while many constants in the following depend on these parameters, we consider them as fixed and often do not explicitly mention this dependence.

\begin{remark}
The monotonicity of $f(\tau)$ in Assumption~\ref{assum:S_inverse_bound} is not a restriction 
since multiplying the right-hand side of \eqref{eq:assum_S_inverse_bound} by $\tau$ yields a monotonically 
decreasing function.  
 Furthermore, Assumption~\ref{assum:A3} implies Assumption~\ref{assum:S_inverse_bound} (cf.\ Lemma~\ref{lmm:Resolvent control for cal S}). The weaker Assumption~\ref{assum:S_inverse_bound} is imposed to exclude eigenvalues away from the support of the asymptotic spectral density  of $X$, while Assumption~\ref{assum:A3} is imposed to guarantee convergence of the empirical eigenvalue distribution to this density in the spectral bulk. 
\end{remark}

The first main result is that, with very high probability, $X$ does not have any eigenvalues away from the disk of radius $\sqrt{\varrho(\cal{S})}$ centred at the origin. 
This will be proven in Section~\ref{sec:proof_no_eigenvalues} below.

\begin{theorem}[No eigenvalue outliers] \label{thm:no_eigenvalues} 
Let $X$ satisfy \ref{assum:A1}, \ref{assum:A2} and \ref{assum:S_inverse_bound}. 
Then, for every $\nu \in \N$ and $\tau_\ast>0$, there exists a constant $C_\nu>0$ such that
\[
\P\sb{
\abs{\zeta}^2 \le\2 \spradius(\cal{S})+\tau_\ast \text{ for all } \zeta \in \spec X 
}
\,\ge\, 1 - C_\nu \, n^{-\nu}\,,
\]
uniformly for all $n \in \N$. 
\end{theorem}

The next theorem states that, for large $n$, the empirical spectral distribution $\frac{1}{n} \sum_{\zeta \in \spec X} \delta_\zeta$ is well approximated by a deterministic probability density $\sigma$ on the complex plane.

\begin{theorem}[Global inhomogeneous circular law]  \label{thm:global_law_for_X} 
Let $X$ satisfy \ref{assum:A1} -- \ref{assum:A4}. 
Then there is a (possibly $n$-dependent) deterministic probability density $\sigma \colon \C \to [0,\infty)$ such that the empirical spectral distribution of $X$ approaches $\sigma(\zeta) \dd^2 \zeta$  weakly in probability for $n \to \infty$. That is, for every bounded, continuous function $f \colon \C \to \C$ and $\eps>0$, we have 
\[ \lim_{n\to \infty} \P \bigg( \absbb{\frac{1}{n} \sum_{\zeta \in \spec X} f(\zeta) - \int_{\C} f(\zeta) \sigma(\zeta) \dd^2\zeta} 
 > \eps \bigg) = 0. \] 
\end{theorem} 

The proof of Theorem~\ref{thm:global_law_for_X} will be presented in Section~\ref{sec:proof_global_law_X} below.
The density $\sigma$ will be  explicitly  defined in \eqref{definition of sigma} below in terms of the solution to a system of two coupled $n \times n$-matrix equations 
 determined by the operators $\cal{S}$ and $\cal{S}^*$ from \eqref{definition cal S}. 
The existence and uniqueness of this solution is stated in the following proposition, whose proof is deferred to the end of Subsection~\ref{Subsec:Resolvent control on scr L} below. 
\begin{proposition}[Existence and uniqueness]
\label{prp:Existence and uniqueness}
Let $X$ satisfy \ref{assum:A1}--\ref{assum:A3}, $\cal{S}$, $\cal{S}^*$ be defined as in \eqref{definition cal S} and $\tau \in  [0,\spradius)$, where $\spradius=\spradius(\cal{S})$ is the spectral radius of $\cal{S}$.  
Then the coupled system of matrix equations 
\bels{Dyson equation at eta = 0}{
\frac{1}{V_1}\,=\,\cal{S}V_2+\frac{\tau}{\cal{S}^*V_1}\,,
\qquad
\frac{1}{V_2}\,=\, \cal{S}^*V_1+\frac{\tau}{\cal{S}V_2}\,,
}
has a unique solution $V_1(\tau)=V_1,V_2(\tau)=V_2 \in \C^{n \times n}$ such that both $V_i$ are positive definite  and satisfy the constraint
\bels{constraint at eta = 0}{
 \tr V_1 = \tr V_2. 
}
This solution can be extended to real analytic functions $V_1,V_2:(-c, \spradius) \to \C^{n \times n}$ with some $n$-independent constant $c>0$. 
\end{proposition}
We will refer to \eqref{Dyson equation at eta = 0} as the \emph{Dyson equation} since as we will see later in Section~\ref{Sec:Dyson equation and its stability} it is equivalent to a Dyson equation that describes the limit of the resolvent of self-adjoint random matrices. 
Our first theorem expresses the density $\sigma$ in terms of the solution to \eqref{Dyson equation at eta = 0} and shows that its support is a disk centred at the origin of the complex plane. It is proven at the end of Section~\ref{Subsec:Solution close to the edge}.  

\begin{theorem}[Density] 
\label{thr:Density}
Let $X$ satisfy \ref{assum:A1}--\ref{assum:A3}, $V_1$, $V_2$  be  the unique positive definite solution of \eqref{Dyson equation at eta = 0} with \eqref{constraint at eta = 0}  and $\spradius=\spradius(\cal{S})$. Then the radially symmetric function $\sigma: \C \to \R$ given by
\bels{definition of sigma}{
\sigma(\zeta)\,:=\, \frac{1}{\pi \1n }\frac{\dd}{\dd \tau}  \tr\sbb{\frac{\tau}{\tau+ (\cal{S}^*V_1(\tau)) (\cal{S}V_2(\tau))}}_{\tau = \abs{\zeta}^2}\msp{-30}\times\, \bbm{1}(\abs{\zeta}^2 < \spradius)\,,
}
is non-negative and inherits its analyticity as a function of $\tau =\abs{\zeta}^2$  in the disk $\D_{\sqrt{\spradius}}= \{ \zeta \in \C \colon \abs{\zeta}^2 < \spradius\}$ from $V_1$  and  $V_2$.  Furthermore, $\sigma$ is a probability density, $\int_{ \C} \sigma(\zeta)\2 \dd^2 \zeta\,=\, 1$ and  is uniformly bounded and bounded away from zero, i.e.\ there are $n$-independent constants $c,C>0$ such that
\bels{lower and upper bounds on sigma}{
c \,\le\, \sigma(\zeta)\,\le\, C\,,\qquad \zeta \in  \D_{\sqrt{\spradius}}\,.
}
 In particular, $\supp \sigma = \overline{\D_{\sqrt{\spradius}}}$ and
at the boundary $\abs{\zeta}^2=\spradius$, the density $\sigma$ has a jump height 
\bels{jump height of sigma}{
\lim_{\abs{\zeta} \uparrow \sqrt{\spradius}} \sigma(\zeta) \,=\, \frac{1}{\pi\1\spradius\2n }\frac{\tr [S_1S_2]^2}{\tr[(S_1S_2)^2]}\,,
}
expressed in terms of the right and left Perron-Frobenius eigenmatrices of $\cal{S}$, i.e.\ $\cal{S}S_2 =\spradius \2S_2$ and $\cal{S}^*S_1 =\spradius\2S_1$. 
\end{theorem}
\begin{definition}[Self-consistent density of states]
We call the probability density $\sigma$, defined through \eqref{definition of sigma}, the \emph{self-consistent density of states} associated to $\cal{S}$ or to $X$.
\end{definition}

In order to formulate  the local law in the spectral bulk we introduce observables around a fixed spectral parameter $\zeta_0 \in \C$ on mesoscopic scales $n^{-\alpha}$ with $\alpha \in (0,1/2)$. 
For any function $f\colon \C \to \C$  we define 
\[
f_{\zeta_0,\alpha}\colon \C \to \C\,, \qquad f_{\zeta_0,\alpha}(\zeta) := n^{2\alpha} f(n^\alpha(\zeta -\zeta_0) )\,.
\]
For any $r>0$, we denote the disk of radius $r$ centred at the origin by 
 $\D_r := \{ \zeta \in \C \colon \abs{\zeta} < r \}$.  

\begin{theorem}[Local inhomogeneous circular law] 
\label{thr:Local inhomogeneous circular law}
Let $X$ be a centred  non-Hermitian  random matrix satisfying \ref{assum:A1}--\ref{assum:A4}. 
Let $\alpha \in (0,1/2)$, $\eps$, $\tau_* >0$ and $\nu \in \N$.
Then there is a constant $C>0$ such that
\begin{equation} \label{eq:local_law_X}
\P\sbb{
\absbb{\frac{1}{n}\sum_{ \zeta \in \spec X} f_{\zeta_0,\alpha}( \zeta  )-\int_\C f_{\zeta_0,\alpha}(\1\zeta\1)\2 \sigma(\1\zeta\1)\dd^2 \zeta}
\2\le\2 n^{-1+2\alpha+\eps}\norm{\Delta f}_{\rm{L}^1}
}
\,\ge\, 1 - C \, n^{-\nu}\,,
\end{equation}
uniformly for every $n \in \N$, every $\zeta_0 \in \C$ with $\abs{\zeta_0}^2 \le \spradius(\cal{S})-\tau_*$ and for every $f \in \rm{C}^2_0(\C)$ satisfying $\supp f \subseteq \D_\varphi$ 
 and $\norm{\Delta f}_{\rm{L}^{1 + \beta}} \leq n^D \norm{\Delta f}_{\rm{L}^1}$ 
with some fixed $\varphi$, $\beta > 0$ and $D \in \N$.  
\end{theorem}

The proof of Theorem~\ref{thr:Local inhomogeneous circular law} will be given in Section~\ref{sec:proof_local_law_X} below. 
Under the stronger Assumption~\ref{assum:bounded_density}, 
the condition $\norm{\Delta f}_{\rm{L}^{1 + \beta}} \leq n^D\norm{\Delta f}_{\rm{L}^1}$ 
in Theorem~\ref{thr:Local inhomogeneous circular law} is not necessary 
(as explained at the beginning of its proof). 
However, if the eigenvalue distribution has a discrete component, then control on $\norm{\Delta f}_{\rm{L}^1}$ alone may not ensure convergence of the linear statistics of $f$ in the eigenvalues in \eqref{eq:local_law_X} which coincides with the integral of $\Delta f$ against the log-potential of the empirical spectral measure (see \eqref{eq:girko_1} below).

As a corollary we prove complete delocalisation of the eigenvectors of $X$. In the case of independent entries eigenvector delocalisation was first proven in \cite{rudelson2015}.

\begin{corollary}[Isotropic eigenvector delocalisation] \label{cor:delocalization} 
Let $X$ satisfy \ref{assum:A1}--\ref{assum:A3}. 
For any $\tau_*>0$, let $\cal{U}_{\tau_\ast}$ denote the set of eigenvectors of $X$ with corresponding eigenvalue in $\D_{\sqrt{\spradius-\tau_\ast}}$ with $\spradius=\spradius(\cal S)$. Then for any $\eps>0$ and $\nu \in \N$ there exists a constant $C_{\eps,\nu}$ such that
\[
\P\sB{\abs{\scalar{v}{u}} \le n^{-1/2 +\eps}\norm{u}\norm{v} \text{ for all } u \in \cal{U}_{\tau_\ast}} \ge 1-C_{\eps,\nu}\, n^{-\nu}\,,
\]
for all $n \in \N$ and all $v \in \C^n$.  
\end{corollary}

Corollary~\ref{cor:delocalization} will be proven in Section~\ref{subsec:proof_consequences_local_law} below.

\subsection{Brown measure of matrix-valued circular elements}

We now illustrate how the probability density defined in \eqref{definition of sigma} is interpreted as the Lebesgue density of the Brown measure associated to a matrix linear combination of circular operators and thus how Theorem~\ref{thr:Density} provides information about this measure.  
To that end, let $(\cal M, \tau)$ be a  tracial $W^*$-probability space\footnote{For this and other basic notions in free probability theory, we refer to the recent monograph \cite{MingoSpeicher}.}. 
For $\ell \in\N$, free circular elements $c_1$, \ldots, $c_\ell \in \cal M$ 
and deterministic matrices $a_1$, \ldots, $a_\ell \in \C^{n\times n}$, we consider the 
operator 
\begin{equation} \label{eq:Kronecker_circular} 
X = \sum_{j=1}^\ell a_{j} \otimes c_j \in \cal M^{n\times n}. 
\end{equation}
We are interested in the spectral distribution of $X$. Since $X$ is non-normal, we 
consider the Brown measure, a generalisation of the spectral measure for normal operators.
The \emph{Brown measure} $\mu_X$ of $X$ is the unique compactly supported probability measure on $\C$ such that 
\begin{equation} \label{eq:def_Brown_measure} 
 \int_\C \log \abs{\lambda- \zeta} \, \mu_X(\dd \zeta) = \log D(X- \lambda) 
\end{equation}
for all $\lambda \in \C$, where $D$ is the Fuglede-Kadison determinant on $(\cal M^{n\times n}, \avg{\,\cdot\,} \otimes \tau)$ defined by 
\begin{equation} \label{eq:def_fuglede_kadison_det} 
 D (Y) := \lim_{\eps \downarrow 0} \exp( \avg{\,\cdot\,} \otimes \tau ( \log(Y^* Y + \eps)^{1/2} )) \in [0,\infty),
\end{equation}
for any $Y \in \cal M^{n\times n}$. 
The Brown measure was originally introduced in \cite{Brown1986} and revived in \cite{HaagerupLarsen2000}. 
The Fuglede-Kadison determinant was first defined in \cite{FugledeKadison1952}. 
For an introduction to both of these objects, we refer to \cite[Section 11]{MingoSpeicher}.

In the next result, we express the Brown measure $\mu_X$ of $X$ from \eqref{eq:Kronecker_circular} in terms of the operators $\cal S$ and $\cal S^*$ on $\C^{n\times n}$ defined through 
\begin{equation} \label{eq:def_self_energy_operator_Brown} 
 \cal S [R] := \sum_{j=1}^\ell a_{j} R a_{j}^*, \qquad \cal S^* [R]:=\sum_{j=1}^\ell a_{j}^* R a_{j} 
\end{equation}
for any $R \in \C^{n\times n}$. 
In particular, we identify the support of $\mu_X$ and classify its regularity.

\begin{proposition}[Regularity of $\mu_X$] \label{pro:Brown_measure} 
Let $X \in \cal M^{n\times n}$ be defined as in \eqref{eq:Kronecker_circular}. 
We assume that there are constants $C>c>0$ such that 
\begin{equation} \label{eq:condition_flatness_Brown} 
  c \avg{R} \leq \cal S[R]  \leq C \avg{R},   
\end{equation}
for all positive semidefinite $R \in \C^{n\times n}$.
Then the Brown measure $\mu_X$ of $X$ is given by 
\begin{equation}\label{eq:mu_X_equals_sigma} 
 \mu_X(\dd \zeta) = \sigma(\zeta) \dd^2 \zeta, 
\end{equation}
where $\sigma$ is defined via \eqref{definition of sigma} with $\cal S$ and $\cal S^*$ from \eqref{eq:def_self_energy_operator_Brown}. 
In particular, the Brown measure of $X$ has all properties of $\sigma$ stated in Theorem~\ref{thr:Density}. 
\end{proposition}

The proof of Proposition~\ref{pro:Brown_measure} is presented in Section~\ref{sec:proof_Brown_measure} below.

\subsection{Relaxed assumptions and examples} \label{sec:relaxed_assumptions} 
 In this subsection, we explain how the assumptions \ref{assum:A1} -- \ref{assum:S_inverse_bound} 
are related, how some of them can be relaxed 
and provide some concrete examples satisfying these assumptions.

The first result shows that \ref{assum:bounded_density} implies \ref{assum:A4} and follows directly from Proposition~\ref{pro:smallest_singular_value_general} below. 

\begin{proposition}[Smallest singular value of $X-\zeta$] \label{pro:smallest_singular_value} 
If $X$ satisfies Assumption~\ref{assum:bounded_density}, then it also 
satisfies Assumption~\ref{assum:A4}.  
\end{proposition}

\begin{remark}[Relaxing Assumption~\ref{assum:A2}]  \label{rem:weaker_version_A2} 
We chose to assume a decay of correlation within the matrix $X$ in the form \ref{assum:A2} because it is easy to state. However, for our proof it suffices to assume that the decay of correlation \eqref{Decay of correlation} holds with a fixed power $\nu>12 p$ with $p \in \N$ from \eqref{eq:growth_condition}, provided higher order cumulants of the matrix entries of $X$ satisfy a certain compatibility condition. This compatibility condition is  \cite[equation (3b)]{Erdos2017Correlated}, where $d$ is interpreted as the pseudometric from Assumption~\ref{assum:A2} and $W$ is replaced by $\sqrt{n}X$. In this case $\bs{H}_\zeta$ from \eqref{Hzeta} satisfies \cite[Assumptions (C) and (D)]{Erdos2017Correlated}  (see also  \cite[Remark~2.7]{Erdos2017Correlated}).
On the other hand, Assumption~\ref{assum:A2} implies \cite[Assumption (C)]{Erdos2017Correlated} and  a modified version of \cite[Assumption (D)]{Erdos2017Correlated} 
 by a similar argument as was used in \cite[Example 2.10]{Erdos2017Correlated}. This is made explicit in Lemma~\ref{lmm:A2 implies C and D} below.
\end{remark} 

In analogy to \cite[Example 2.12]{Erdos2017Correlated} we also provide a simple description of our assumptions for the case of Gaussian random matrices while relaxing the polynomial decay of correlations from \eqref{Decay of correlation} to be of order $\nu=2$ when $d(i,j)=\abs{i-j}$ is the standard metric on $\qq{n}$. 

\begin{example}[Results for Correlated Gaussian matrices] 
Let $X \in \C^{n\times n}$ be a random matrix with  centred Gaussian entries such that 
\[
n(\abs{\E\2 x_{ij} \1x_{lk}}+\abs{\E\2 \ol{x}_{ij} \1x_{lk}})\le \frac{C}{\abs{i-l}^2 +\abs{j-k}^2}
\]
for all $i,j,l,k \in \qq{n}$, as well as $\E\2 \abs{\tr BX}^2\ge \frac{c}{n} \tr B^*B$ for all $B \in \C^{n \times n}$, where $c,C>0$ are some positive constants.
Then the conclusions of Theorem~\ref{thm:no_eigenvalues}, Theorem~\ref{thm:global_law_for_X}, 
Theorem~\ref{thr:Local inhomogeneous circular law} and Corollary~\ref{cor:delocalization} 
hold for $X$. 
\end{example} 

Next, we will formulate a condition for block matrices that ensures Assumption~\ref{assum:bounded_density}.  
We denote by $E_{ij} \in \C^{N \times N}$ the matrix whose $(i,j)$-entry is 1 and whose other entries are zero, that is, $E_{ij} = (\delta_{ik}\delta_{jl})_{k,l \in \qq{N}}$. 
In the following lemma, we write $\mathbf{z}$ for a matrix-valued variable $\mathbf{z} = (z_{\gamma \delta})_{\gamma,\delta \in \qq{K}}$. 
We denote by $\dd \mathbf{z}$ integration with respect to all entries of $\mathbf{z}$
and $\wh{\dd^2 z_{\alpha \beta}}$ denotes the omission of the integration over $z_{\alpha \beta}$. 

\begin{lemma}[Block matrices]  \label{lem:A4_for_block_matrices} 
Let $K \in \N$ be fixed. 
Let $\{x_{ij} \colon i,j \in \qq{N} \}$ be a family of independent random matrices in $\C^{K\times K}$ 
satisfying $\E x_{ij} = 0$ for all $i,j \in \qq{N}$. 
We assume that, for all $i,j \in \qq{N}$, the matrix  $x_{ij}\sqrt{NK}$ has a density $f_{ij}$ on $\C^{K\times K}$, i.e.\ 
\[ \P \big(  x_{ij} \sqrt{NK} \in B \big) = \int_{B} f_{ij} ( \mathbf{z}) \dd \mathbf{z} \] 
for all measurable subsets $B \subset \C^{K\times K}$. If there are $q>0$ and $C>0$ such that 
\begin{equation} \label{eq:condition_density_block_matrix} 
 \int_{\C^{K\times K  -1}} \bigg( \int_\C f_{ij}(\mathbf{z})^q \dd^2 z_{\alpha\beta} \bigg)^{1/q} \dd^2 z_{11} \dd^2 z_{12} \ldots \wh{\dd^2 z_{\alpha\beta}} \ldots \dd^2 z_{KK} \leq N^C 
\end{equation}
for all $\alpha,\beta \in \qq{K}$ then Assumption \ref{assum:bounded_density} is satisfied for the block matrix 
\begin{equation} \label{eq:def_Kronecker_matrix_X} 
 X = \sum_{i,j \in \qq{N}} x_{ij} \otimes E_{ij}. 
\end{equation}
\end{lemma} 

An analogous statement holds when $x_{ij} \sqrt{N K}$ has a density on $\R^{K\times K}$ instead of $\C^{K\times K}$. 
Lemma~\ref{lem:A4_for_block_matrices} will be proven in Section~\ref{sec:check_examples} below.

\subsection*{Notations} 

Here we introduce some notations that will be used throughout the paper. We start with basic notations for matrices. We equip the space of $d \times d$-matrices  with the normalised scalar product
\[
\scalar{A}{B} := \frac{1}{n} \tr A^* B\,, \qquad A,B \in \C^{d \times d},
\]
 corresponding norm  $\norm{A}_{\rm{hs}}^2:= \scalar{A}{A}$
and use the short hand $\avg{A} = \frac{1}{d}\tr A$ for the normalised trace. By $\norm{A}$ we denote the operator norm induced by the standard Euclidean metric on $\C^d$. More generally, for linear operators $A: \cal{A} \to \cal{B}$ from a normed space $\cal{A}$ to a normed space $\cal{B}$, we indicate the corresponding operator norm by writing $\norm{A}_{\cal{A} \to \cal{B}}$ and simply $\norm{A}_{\cal{A}}$ in case $\cal{A}=\cal{B}$. Since we often work with $2\times 2$-block matrices having block dimension $n$, we will frequently use the block notation from \eqref{eq:def_S_intro}, where $\bs{R} \in \C^{2n \times 2n}$ and $R_{ij} \in \C^{n \times n}$. 

For nonnegative quantities $\phi,\psi$ we use the comparison relation $\phi\lesssim\psi$ whenever $\phi\le C\psi$ with an $n$-independent constant $C>0$. This constant is uniform in all parameters except the model parameters from Assumptions  \ref{assum:A1}--\ref{assum:A4} and possibly other parameters that are either clearly indicated or obvious from the context. In particular, $C$ is uniform in the spectral parameter $\zeta$ within the domain under consideration. If $c\2\psi\le \phi\le C\psi$ we write $\phi \sim \psi$ and $\phi = \psi + \ord(\nu)$ is a short hand for $\abs{\phi-\psi}\lesssim \nu$. We also use the comparison relation for positive definite matrices, where it is interpreted in a quadratic form sense.

\section{Inhomogeneous circular law} \label{sec:inhomogeneous_circular_law}

In this section we prove Theorems~\ref{thm:no_eigenvalues} and \ref{thm:global_law_for_X}. These proofs will illustrate how Girko's Hermitization trick 
translates these questions to Hermitian random matrices which will be analysed via their resolvents and the associated matrix Dyson equation. 
The proof of Theorem~\ref{thm:global_law_for_X} is a prototype of the more complicated proof of Theorem~\ref{thr:Local inhomogeneous circular law} in Section~\ref{sec:proof_local_law_X} below. 

The fundamental observation due to Girko \cite{Girko1985} is that $\zeta \in \C$ is an eigenvalue of $X$ if and only if the kernel of $\bs{H}_\zeta$ is nontrivial, 
where the Hermitian matrix $\bs{H}_\zeta$ is defined through
\begin{equation} \label{eq:def_H_zeta} 
 \bs{H}_\zeta :=\mtwo{0 & X-\zeta }{(X-\zeta)^* & 0}.  
\end{equation}
The family $(\bs{H}_\zeta)_{\zeta \in \C}$ is called the \emph{Hermitization} of $X$. 
All spectral information about the kernel of $\bs{H}_\zeta$ is captured by the resolvent $\bs{G} = \bs{G}(\zeta, \eta)$ 
of 
$\bs{H}_\zeta$ defined by 
\bels{definition resolvent G}{ \bs{G}(\zeta, \eta)\,:=\, (\bs{H}_\zeta -\ii \1\eta)^{-1}, }
where $\zeta \in \C$ and $\eta>0$.

We will see in Proposition~\ref{pro:global_law_H_averaged} below that the resolvent $\bs{G}$ is well approximated by the matrix $\bs{M} = \bs{M}(\zeta, \eta)\in \C^{2n\times 2n}$ 
which is the unique solution of the \emph{matrix Dyson equation} (MDE) 
\begin{equation} \label{eq:MDE} 
 - \bs{M}^{-1} =  \ii \eta \bs{1}  + \bs{Z}(\zeta,\bar \zeta) + \mathscr{S}[\bs{M}], \qquad \qquad \eta > 0, \quad \zeta \in \C \,,
\end{equation}
under the constraint that the imaginary part $\Im \bs{M} = \frac{1}{2\ii} ( \bs{M} - \bs{M}^*)$ is 
positive definite. 
 Here, the matrix-valued function $\bs{Z}\colon \C^2 \to \C^{2n \times 2n}$ and the  
\emph{self-energy operator} $\scr S$, a linear operator on $ \C^{2n \times 2n}$,  are defined through
 \bels{definition scr S and Z}{ 
 \bs{Z}(\zeta, \omega)\,:=\, \mtwo{0& \zeta}{\omega & 0}, \qquad \qquad 
\mathscr{S}\bigg[\begin{pmatrix} A_{11} & A_{12}\\ A_{21} & A_{22} \end{pmatrix}  \bigg] := \begin{pmatrix} \cal{S} [A_{22}] & 0 \\ 0 & \cal{S}^*[A_{11}] \end{pmatrix}, 
} 
where all blocks in these matrix representations are of size $n\times n$ (see \eqref{definition cal S} for the definitions of $\mathcal{S}$ and $\mathcal{S}^*$).
The existence and uniqueness of $\bs{M}$ have been shown in \cite{Helton01012007}.

We represent $\bs{M}$ in terms of the $2 \times 2$-block structure corresponding to  the right-hand side of \eqref{eq:MDE}.  
For this purpose we first introduce the matrices $V_1$, $V_2\in \C^{n\times n}$ which are the unique solution of 
\begin{subequations}
\label{Dyson equation}
\begin{align}
\label{V1 equation}
\frac{1}{V_1(\tau,\eta)}\,&=\, \eta +\cal{S}V_2(\tau,\eta)+\frac{\tau}{\eta+\cal{S}^*V_1(\tau,\eta)}\,,
\\
\label{V2 equation}
\frac{1}{V_2(\tau,\eta)}\,&=\, \eta +\cal{S}^*V_1(\tau,\eta)+\frac{\tau}{\eta+\cal{S}V_2(\tau,\eta)}\,,
\end{align}
\end{subequations}
for any $\eta >0$ and $\tau \geq 0$ under the constraint that $V_1$ and $V_2$ are positive definite. 
We note that \eqref{Dyson equation} is a regularised version of the Dyson equation \eqref{Dyson equation at eta = 0}, used for the definition of $\sigma$ in \eqref{definition of sigma}, 
with some regularisation parameter $\eta>0$. 
Moreover, we introduce the auxiliary matrix   
\bels{definition of U}{
U(\tau,\eta)\,:=\,  \frac{1}{\tau+(\eta + \cal{S}^*V_1(\tau,\eta))(\eta + \cal{S}V_2(\tau,\eta))}\,. 
}
Then we obtain that 
\bels{definition of M}{
 \bs{M}(\zeta,\eta)\,=\, \mtwo{\ii\1 V_1(\abs{\zeta}^2,\eta) &-\zeta\2U(\abs{\zeta}^2,\eta)}{-\ol{\zeta}\2U(\abs{\zeta}^2,\eta)^*&\ii\1 V_2(\abs{\zeta}^2,\eta)}\,
 }
since the right-hand side of \eqref{definition of M} satisfies \eqref{eq:MDE} and has a positive definite imaginary part. 
Thus solving \eqref{eq:MDE} for $\bs{M}$ with positive imaginary part is equivalent to solving \eqref{Dyson equation} for positive definite $V_1,V_2$. 
From \eqref{definition of M} we easily get that 
\begin{equation} \label{eq:block_representation_im_M} 
 \Im \bs{M}(\zeta,\eta) = \begin{pmatrix} V_1 & 0 \\ 0& V_2 \end{pmatrix}. 
\end{equation}

\subsection{Exclusion of eigenvalues away from the disk}  \label{sec:proof_no_eigenvalues} 

We now prepare the proof of Theorem~\ref{thm:no_eigenvalues}. First we  note  that if $X$ satisfies \ref{assum:A1} and \ref{assum:A2} then, 
 for all positive definite $R \in \C^{n\times n}$, we have  
\begin{equation} \label{eq:upper_bound_flatness} 
\cal{S}[R] \lesssim \avg{R}, \qquad \qquad \cal{S}^*[R] \lesssim \avg{R}.
\end{equation}
The next lemma describes the behaviour of $\Im \bs{M}(\zeta,\eta)$ when $\abs{\zeta}^2 > \varrho(\cal{S})$. 

\begin{lemma} \label{lem:Im_M_sim_eta} 
Let $X$ satisfy \ref{assum:A1}, \ref{assum:A2} and \ref{assum:S_inverse_bound}. 
Then, for every $\delta >0$, we have 
\bels{scaling of im M outside}{ \Im \bs{M}(\zeta,\eta) \sim_{\delta} \frac{\eta}{\abs{\zeta}^2}  }
for all $\eta \in (0,1]$ and $\zeta \in \C$ satisfying $\abs{\zeta}^2 \geq \varrho(\cal{S}) +\delta$. 
\end{lemma}

\begin{Proof}[Proof of Lemma~\ref{lem:Im_M_sim_eta}]  
Multiplying \eqref{V1 equation} with $\eta+\cal{S}^*[V_1]$ from the left and \eqref{V2 equation} from the right with $\eta+\cal{S}[V_2]$ as well as realizing that the resulting right hand sides coincide reveal the identity
 \bels{U identity2}{
 (\eta+ \cal{S}^*V_1)\frac{1}{V_1}\,=\, \frac{1}{V_2}(\eta+ \cal{S}V_2)\,.
 }
Taking the inverse on both sides of \eqref{U identity2} and applying this identity to the result of 
multiplying \eqref{V1 equation} with $V_1$ from the right and  with $\eta + \cal{S}^*V_1$ from the left  
yield 
\begin{equation} \label{eq:V_2_symmetric} 
 \eta + \cal{S}V_2 = (\eta + \cal{S} V_2)V_1 (\eta + \cal{S}V_2) + \tau V_2. 
\end{equation}
We reorganize the terms in \eqref{eq:V_2_symmetric}, use that $\tau - \cal{S}$ is invertible as $\tau > \varrho(\cal{S})$ 
and obtain 
\[ V_2  = (\tau- \cal{S})^{-1} \big( \eta - (\eta + \cal{S} V_2)V_1 (\eta + \cal{S}V_2) \big) 
\leq \eta (\tau - \cal{S})^{-1} 1. \] 
Here, we used in the last step that $(\tau - \cal{S})^{-1}$ is positivity preserving due to the Neumann series 
and that $(\eta + \cal{S} V_2)V_1 (\eta + \cal{S}V_2) \geq 0$. 
Therefore, we have shown that 
\[ V_2 \leq \eta f(\tau)/\tau \leq \eta f(\varrho(\cal{S}) + \delta)/\tau \lesssim_\delta \eta/\tau \] 
for all $\eta>0$ and all $\tau \geq \varrho(\cal{S}) + \delta$. 
Similarly, we get $V_1 \lesssim_\delta \eta/\tau$. 

Using $V_1 + V_2 \lesssim_\delta \eta /\tau$ and \eqref{eq:upper_bound_flatness} to estimate the right-hand side of \eqref{eq:V_2_symmetric} from above implies 
\[ \eta \lesssim_{\delta} \eta^3 + \tau V_2. \] 
Hence, $V_2 \gtrsim_{\delta} \eta/\tau$ and a similar argument yields $V_1 \gtrsim_{\delta} \eta/\tau$. 
Owing to \eqref{eq:block_representation_im_M}, 
these estimates and $\tau = \abs{\zeta}^2$ complete the proof of Lemma~\ref{lem:Im_M_sim_eta}. 
\end{Proof}

For the upcomming arguments, it is convenient to use the following notion of events that occur with ``very high probability''.

\begin{definition}[With very high probability] \label{def:very_high_probability} 
We say that the (sequence of) events $(A_n)_{n \in \N}$ occur \emph{with very high probability} if 
for every $\nu>0$ there is $C_\nu>0$ such that 
\begin{equation} \label{eq:def_very_high_probability} 
\P \big( A_n \big) \geq 1 - C_\nu n^{-\nu} 
\end{equation}
for all $n \in \N$. 
\end{definition} 
The constants $C_{\nu}$ in \eqref{eq:def_very_high_probability} will typically depend on the model parameters. 
Note that an intersection of $n^C$-many events holding with very high probability also holds with very high probability.

\begin{Proof}[Proof of Theorem~\ref{thm:no_eigenvalues}] 
The theorem will follow from the next lemma and an interpolation argument. 
As we will see in its proof in Appendix~\ref{app:tools_section_3} below, this lemma is a direct 
consequence of \cite[Corollary~2.3]{Erdos2017Correlated} and Lemma~\ref{lem:Im_M_sim_eta}. 

\begin{lemma}[No eigenvalues of $\bs{H}_\zeta$ around zero] \label{lem:no_eigenvalues_H}  
Let $X$ satisfy \ref{assum:A1}, \ref{assum:A2} and \ref{assum:S_inverse_bound}.  
If $\zeta \in \C$ satisfies $\sqrt{\varrho(\cal S)} + \delta \leq \abs{\zeta}^2 \leq \delta^{-1}$ for some $\delta \sim 1$ then 
there is $\eps \sim 1$ such that 
\[ \spec(\bs{H}_{\zeta}) \cap (-\eps, \eps) = \emptyset \] 
with very high probability. 
\end{lemma} 
Since $\rm{Spec}(X)=\{\zeta\in \C :0  \in \rm{Spec}(\bs{H}_\zeta) \}$ we conclude 
 from Lemma~\ref{lem:no_eigenvalues_H}  that with  
very  high probability $X$  has no eigenvalues in the annulus  $A:=\{\zeta: \varrho(\cal{S})+\delta < \abs{\zeta}^2<\delta^{-1}\}$. 
We will now show that there are no eigenvalues of  $X$ outside $\D_{\delta^{-1/2}}$ either. 
For this purpose we apply Lemma~\ref{lem:no_eigenvalues_H} to the Hermitization $t\bs{H}_{\zeta/t}$ of $tX$ for any $t \in [0,1]$.  
We choose a finite subset $Z \subset A$ such that $Z + \D_{n^{-1}}$ covers the entire annulus $A$ and $\abs{Z} \le n^{C}$ for some $C>0$. By a union bound and Lemma~\ref{lem:no_eigenvalues_H}, we find that for any $(t, \zeta) \in n^{-1}\llbracket n\rrbracket \times Z$ with very high probability $\rm{Spec}(t\bs{H}_{\zeta/t}) \cap (-\eps,\eps) = \emptyset$. Thus, by Lipschitz-continuity of $t\bs{H}_{\zeta/t}$ in $t$ and $\zeta$, we have $0 \not \in \cup_{t \in [0,1]}\cup_{\zeta \in A}\rm{Spec}(t\bs{H}_{\zeta/t})$ with very high probability. 
In particular, the eigenvalues of each matrix  along the interpolation $t \mapsto tX$ between the zero matrix and $X$, that continuously depend on $t$, do not cross the annulus. Therefore,  $X$ has the same number of eigenvalues inside the disk with radius $\sqrt{\varrho(\cal{S})+\delta}$ as the zero matrix, namely $n$, i.e.\ it has no eigenvalues outside this disk.
\end{Proof}

\subsection{Global inhomogeneous circular law} \label{sec:proof_global_law_X}

In this section, we prove Theorem~\ref{thm:global_law_for_X}.  
We first derive the basic formula relating the eigenvalue density of $X$ to the Hermitian matrices $\bs{H}_\zeta$ defined in \eqref{eq:def_H_zeta}.
This approach goes back to Girko \cite{Girko1985}.  
Then we motivate and collect all other ingredients required for the proof of Theorem~\ref{thm:global_law_for_X}.

The starting point is a relation for the averaged linear statistics with a test function $f \in C_0^2(\C)$ given by 
\begin{equation} \label{eq:girko_1} 
 \frac{1}{n} \sum_{\xi \in \spec X} f (\xi) = \frac{1}{2\pi n} \sum_{\xi \in \spec X} \int_{\C} \Delta f(\zeta) \log \abs{ \zeta - \xi} \dd^2 \zeta,
\end{equation}  
where we used in the first step that $\log$ is the fundamental solution of the Laplace equation in $\R^2$. 

The right-hand side of \eqref{eq:girko_1} can be expressed purely in terms of the Hermitian matrices $\bs{H}_\zeta$ since 
\begin{equation} \label{eq:linear_stat_log_det_H_zeta} 
 \sum_{\xi \in \spec X} \log\abs{\xi-\zeta} = \log \abs{\det(X-\zeta)} = \frac{1}{2} \log \abs{\det \bs{H}_\zeta}. 
\end{equation}
The resolvent $\bs{G}$ contains all spectral information about $\bs{H}_\zeta$. 
In particular, $\log \abs{\det \bs{H}_\zeta}$ is expressed in terms of $\bs{G}$ via the well-known identity 
\begin{equation} \label{eq:log_det_int_im_G} 
 \log \abs{\det \bs{H}_\zeta} = - 2n \int_0^T \avg{\Im \bs{G}(\zeta, \eta)} \dd \eta + \log \abs{\det (\bs{H}_\zeta - \ii T)} 
\end{equation} 
for any $T>0$ (see \cite{tao2015} for the use of \eqref{eq:log_det_int_im_G} in a similar context). 
Hence, owing to \eqref{eq:girko_1}, \eqref{eq:linear_stat_log_det_H_zeta} and \eqref{eq:log_det_int_im_G}, it suffices to control $\bs{G}$ in order to understand the averaged linear statistics. 
As indicated in Section~\ref{sec:proof_no_eigenvalues}, the resolvent $\bs{G}$ will be well approximated by the solution $\bs{M}$ of the MDE \eqref{eq:MDE} for large $n$. 

We now collect some auxiliary results about $\bs{M}$ and $\sigma$. 
We will need the following bounds on $\bs{M}$ proven at the end of Section~\ref{Subsec:Solution}. 

\begin{lemma}[Bounds on $\bs{M}$] \label{lem:bounds_M} 
Let $X$ satisfy \ref{assum:A1} -- \ref{assum:A3}. 
Then, uniformly for $\zeta \in \C$ and $\eta >0$, we have 
\begin{equation} \label{eq:norm_M} 
 \norm{\bs{M}(\zeta,\eta)} \lesssim \frac{1}{1 + \eta + \abs{\zeta}}. 
\end{equation}
Moreover, for any $T > 0 $ and $\zeta \in \C$, we have 
\begin{equation} \label{eq:integral_im_M} 
  \int_0^T \absbb{\avg{\Im \bs{M}(\zeta,\eta)} - \frac{1}{1+ \eta}} \dd \eta \lesssim \min\bigg\{T,1 +  \frac{\abs{\zeta}}{T}  \bigg\}, \qquad \int_T^\infty \absbb{\avg{\Im \bs{M}(\zeta,\eta)} - \frac{1}{1+ \eta}} \dd \eta \lesssim 
 \frac{1 + \abs{\zeta}}{T} . 
\end{equation}
\end{lemma} 

The self-consistent density of states $\sigma$ introduced in \eqref{definition of sigma} relates to $\Im \bs{M}$ in the way expected from \eqref{eq:girko_1}, \eqref{eq:linear_stat_log_det_H_zeta} and \eqref{eq:log_det_int_im_G} 
as well as $\bs{G} \approx \bs{M}$. 
This is the content of the next lemma. 

\begin{lemma}[$\sigma$ as distributional derivative] \label{lem:sigma_as_dist_derivative} 
Let $X$ satisfy \ref{assum:A1} -- \ref{assum:A3}. Then we have 
\begin{equation} \label{eq:sigma_integration_part_def_L} 
\int_\C f(\zeta) \sigma(\zeta) \dd^2 \zeta = -\frac{1}{2\pi} \int_{\C} \Delta f(\zeta) L(\zeta) 
\dd^2 \zeta, \qquad L(\zeta) := \int_0^\infty \bigg( \Im \avg{\bs{M}(\zeta,\eta)} - \frac{1}{1+ \eta} \bigg) \dd \eta
\end{equation}
 for every $f \in C_0^2(\C)$. 
The integral in the definition of $L$ exists in the Lebesgue sense due to \eqref{eq:integral_im_M}. 
\end{lemma} 

Lemma~\ref{lem:sigma_as_dist_derivative} in particular shows that $\Delta L = -2\pi \sigma$ in the sense of distributions, i.e.\ $L$ is the logarithmic potential of the probability measure $\sigma(\zeta) \dd^2 \zeta$. 
The proof of Lemma~\ref{lem:sigma_as_dist_derivative} requires a very detailed analysis of the MDE, \eqref{eq:MDE}, and its stability properties 
and will be presented in Section~\ref{subsec:proof_sigma_dist_derivative} below. 

To illustrate the basic formula used in the proof of Theorem~\ref{thm:global_law_for_X} below, we combine the identities \eqref{eq:girko_1}, \eqref{eq:linear_stat_log_det_H_zeta}, \eqref{eq:log_det_int_im_G} and 
\eqref{eq:sigma_integration_part_def_L} and, thus, obtain for any $T>0$ that 
\begin{equation} \label{eq:decomposition_preparation} 
\begin{aligned} 
 \frac{1}{n} \sum_{\xi \in \spec X} f (\xi) - \int_\C f (\zeta) \sigma(\zeta) \dd^2 \zeta 
  & = \frac{1}{2\pi} \int_{\C} \Delta f(\zeta) \bigg( \int_{0}^T \Im \avg{\bs{M}(\zeta,\eta) - \bs{G}(\zeta,\eta)} 
\dd \eta \\ 
 & \phantom{=} + \int_T^{\infty} \bigg( \Im \avg{\bs{M}(\zeta,\eta)}  - \frac{1}{1 + \eta} \bigg) \dd \eta  
+ \frac{1}{2n} \log \abs{\det( \bs{H}_\zeta - \ii T)} \bigg) \dd^2 \zeta, 
\end{aligned} 
\end{equation}
where we used that $\int_{\C} \Delta f(\zeta) \dd^2 \zeta = 0$. 
The terms on the right-hand side of \eqref{eq:decomposition_preparation} will be bounded as follows. 
The second term is controlled due to the second bound in \eqref{eq:integral_im_M} and the third 
by a simple argument using \ref{assum:A1}. For the first term, we shall use Proposition~\ref{pro:global_law_H_averaged} below and Assumption~\ref{assum:A4}.  

For technical reasons, we discretise the integral over $\zeta$ in \eqref{eq:girko_1} 
via Lemma~\ref{lem:sampling} below. Afterwards, we apply \eqref{eq:linear_stat_log_det_H_zeta} and \eqref{eq:log_det_int_im_G} to the discretised expression. 
Thus, the final proof of Theorem~\ref{thm:global_law_for_X} 
does not start from \eqref{eq:decomposition_preparation} directly. 
For the discretisation of the $\zeta$-integral, we apply the sampling method formulated in the following lemma. 
 For $a = 2$, it is a special case of \cite[Lemma~36]{tao2015}, which was used in a similar context in \cite{tao2015}.   

\begin{lemma}[Monte Carlo sampling] \label{lem:sampling} 
Let $\Omega \subset \C$ be a bounded subset of positive Lebesgue measure and $\mu$ the normalized Lebesgue measure on $\Omega$. 
Let $F \colon \Omega \to \C$ be a function in $\rm{L}^a(\mu)$ from some $a > 1$.   
For $m \in \N$, let $\xi_1, \ldots, \xi_m$ be independent random variables distributed according to $\mu$. 

Then, for any $\delta \in (0,1]$, we have 
\[ \P \bigg( \absB{\frac{1}{m} \sum_{i=1}^m F(\xi_i) - \int_\Omega F \dd \mu} \leq \frac{10^{1/a}}{m^{1-1/a}\delta^{1/a}} \Big(\int_\Omega \absB{F - \int_\Omega F \dd \mu }^a \dd \mu\Big)^{1/a} \bigg) \geq 1 - \delta. \] 
\end{lemma} 

\begin{Proof} 
The random variables $F(\xi_1), \ldots, F(\xi_m)$ are i.i.d.\ with mean $\int_\Omega F \dd \mu$. 
Thus,
Proposition~\ref{pro:law_large_numbers} in Appendix~\ref{app:law_large_numbers} below   
implies Lemma~\ref{lem:sampling}. 
\end{Proof}

The next bound on $\bs{G} - \bs{M}$ is the last missing ingredient for the proof of Theorem~\ref{thm:global_law_for_X}. 

\begin{proposition}[Global law for $\bs{H}_\zeta$, averaged version] \label{pro:global_law_H_averaged} 
Let $X$ satisfy \ref{assum:A1} and \ref{assum:A2}. Then there is an absolute constant $K>0$ such that for all sufficiently small $\delta >0$ we have 
\begin{equation} \label{eq:global_law_averaged} 
 \abs{\avg{\bs{G} (\zeta,\eta) -\bs{M} (\zeta,\eta) }} \leq \frac{n^{K \delta} }{(1 + \eta^2) n} 
\end{equation}
with very high probability uniformly 
for all $n \in\N$, $\zeta \in \mathbb{D}_{\varphi}$ and $\eta \in [n^{-\delta},n^{100}]$.  
\end{proposition} 

Proposition~\ref{pro:global_law_H_averaged} is implied by Proposition~\ref{pro:global_law_H_zeta} below. 
The former has an important consequence, namely the following bound on the number of eigenvalues of $\bs{H}_\zeta$ close to zero.  
Note that the moduli of the eigenvalues of $\bs{H}_\zeta$ are the singular values of $X-\zeta$. 
The eigenvalues of $\bs{H}_\zeta$ are denoted by $\lambda_1(\zeta)$, \ldots, $\lambda_{2n}(\zeta)$. 
Together with Assumption~\ref{assum:A4}, the bound in the next lemma  will be used to control the integral in \eqref{eq:log_det_int_im_G} for small $\eta$. 

\begin{lemma}[Number of small singular values of $X-\zeta$] \label{lem:number_small_singular_values_global} 
Let $X$ satisfy \ref{assum:A1} and \ref{assum:A2}. Then there is $\delta >0$ such that 
\[  \# \{ i \in \llbracket 2n \rrbracket  \colon \abs{\lambda_i(\zeta)} \leq \eta \} \lesssim n \eta  \] 
with very high probability uniformly for all $\eta \in [n^{-\delta}, n^{100}]$ and $\zeta \in \mathbb{D}_{\varphi}$ with any fixed $\varphi>0$. 
Here, the constant $C_\nu$ implicit in the very high probability notion from  Definition~\ref{def:very_high_probability}  depends on $\delta$ and $\varphi$ as well as the constants in \ref{assum:A1} and \ref{assum:A2} in 
addition to $\nu$. 
\end{lemma} 

\begin{Proof}
The trace of $\bs{G}$ is bounded by $n$, $\abs{\tr \bs{G}} \lesssim n$ with very high probability, for all $\eta \in [n^{-\delta},n^{100}]$ due to \eqref{eq:global_law_averaged} and $\norm{\bs{M}} \lesssim 1$ by \eqref{eq:norm_M}. 
Therefore, setting $\Sigma_\eta := \{ i \in \llbracket 2n \rrbracket \colon \abs{\lambda_i(\zeta)} \leq \eta \}$ yields 
\[ \frac{\# \Sigma_\eta}{2 \eta} \leq \sum_{i \in \Sigma_\eta} \frac{\eta}{\eta^2 + \lambda_i(\zeta)^2} \leq \Im \tr \bs{G}(\zeta,\eta) \lesssim n. \qedhere \] 
\end{Proof} 

We will now conclude Theorem~\ref{thm:global_law_for_X} from Proposition~\ref{pro:global_law_H_averaged} and Lemma~\ref{lem:number_small_singular_values_global}.

\begin{Proof}[Proof of Theorem~\ref{thm:global_law_for_X}] 
We will show below that for all sufficiently small $\delta>0$ we have 
\begin{equation} \label{eq:quantitative_global_law_for_X}  
\absbb{\frac{1}{n} \sum_{\zeta \in \spec X} f(\zeta) - \int_{\C} f(\zeta) \sigma(\zeta) \dd^2\zeta } \leq n^{-\delta} \norm{\Delta f}_{\mathrm{L}^1} 
\end{equation} 
with very high probability uniformly for all $f \in C_0^2(\C)$ satisfying $\norm{\Delta f}_{\rm{L}^{1 + \beta}} \leq n^{D} \norm{\Delta f}_{\rm{L}^1}$ and $\supp f \subseteq \mathbb{D}_\varphi$, 
where $\beta>0$, $ D >0$ and $\varphi>0$ are some constants. 
In \eqref{eq:quantitative_global_law_for_X}, the constant $C_\nu$ in the definition \eqref{eq:def_very_high_probability} depends only on $\delta$, $ \beta$, $ D$ and $\varphi$ in addition to $\nu$ and 
the constants from \ref{assum:A1} and \ref{assum:A2}. 

Given \eqref{eq:quantitative_global_law_for_X},  we now explain how Theorem~\ref{thm:global_law_for_X} follows. 
Let $f \in C_b(\C)$ and $\eps>0$. 
Since $X$ does not have any eigenvalues outside $\mathbb{D}_{R+1}$, where $R =  \sqrt{\varrho(\mathcal{S})}$,  
with very high probability by Theorem~\ref{thm:no_eigenvalues}, 
we assume without loss of generality that $\supp f \subset \mathbb{D}_{R+2}$.  
As $R = \sqrt{\varrho(\mathcal{S})} \lesssim 1$  due to  Assumptions~\ref{assum:A1} and \ref{assum:A2}, 
we can choose a constant $\varphi$ such that $\varphi > R + 2$. 
We find $g \in C_0^2(\C)$ such that $\norm{f - g}_{\rm{L}^\infty} \leq \eps/2$, $\supp g \subset \mathbb{D}_\varphi$  and $\norm{\Delta g}_{\rm{L}^{2}} \lesssim_\eps 1$. 
Therefore, approximating $f$ by $g$ in the statement of Theorem~\ref{thm:global_law_for_X} and applying 
\eqref{eq:quantitative_global_law_for_X} to $g$ shows Theorem~\ref{thm:global_law_for_X}. 

What remains is proving \eqref{eq:quantitative_global_law_for_X}. 
We set $\Omega = \mathbb{D}_{\varphi}$. 
Combining \eqref{eq:girko_1} and \eqref{eq:sigma_integration_part_def_L} as well as using the second bound in \eqref{eq:integral_im_M} yield
\begin{equation} \label{eq:decomposition} 
 \frac{1}{n} \sum_{\xi \in \spec X} f (\xi) - \int_\C f (\zeta) \sigma(\zeta) \dd^2 \zeta 
  = \int_{\Omega} F(\zeta) \dd \mu(\zeta)   + \ord\big(T^{-1} \norm{\Delta f}_{\rm{L}^1}\big)  
\end{equation} 
for any $T>0$. 
Here, we denoted by $\mu$ the normalized Lebesgue measure on $\Omega$ and introduced 
\begin{equation} \label{eq:def_F} 
F(\zeta) := \frac{\abs{\Omega}}{2\pi}  \Delta f(\zeta) h(\zeta), \qquad h(\zeta ) : = \frac{1}{n} \sum_{\xi \in \spec X} \log \abs{\xi- \zeta} + \int_0^T \bigg( \avg{\Im \bs{M}(\zeta,\eta)} - \frac{1}{1 + \eta} \bigg) \dd\eta . 
\end{equation}  

We now apply Lemma~\ref{lem:sampling} to the first term on the right-hand side of \eqref{eq:decomposition}. Note that $\zeta \mapsto \log \abs{\xi -\zeta}$ lies in $\rm{L}^p(\Omega)$ for every $p \in [1,\infty)$. 
Hence, owing to the first bound in \eqref{eq:integral_im_M} we get that, for any $p \in [1,\infty)$, $\norm{h}_{\rm{L}^p(\Omega)} \lesssim_p 1$ uniformly for $T >0$. 
In particular, the function $F$ defined in \eqref{eq:def_F} is square-integrable on $\Omega$. 
Thus, Lemma~\ref{lem:sampling} is applicable and choosing 
 $\delta = n^{-\nu}$, $a = 1 + \beta/2$ and $m = n^{(\nu + (D + 11) a)/(a-1)}$ 
shows that 
\begin{equation} \label{eq:error_term_sampling} 
 \absbb{\int F \dd \mu - \frac{1}{m} \sum_{i=1}^m F(\xi_i)} \leq n^{-D - 10} \norm{\Delta f}_{\rm{L}^{1 + \beta}} 
\end{equation}
with very high probability, where $\xi_1$, \ldots, $\xi_m$ are independent random variables distributed according to $\mu$. 

We set $T =n^{100}$ and now show that for all sufficiently small $\delta >0$ we have 
\begin{equation} \label{eq:F_estimate_fixed_zeta} 
 \abs{F(\zeta)}\leq n^{-\delta} \abs{\Delta f(\zeta)} 
\end{equation}
with very high probability uniformly for all $\zeta \in \Omega$.  
To that end, we define $\eta_* := n^{-\delta}$ and   
\[ \begin{aligned} 
h_1(\zeta) & := \, \int_{\eta_*}^T \Im \avg{\bs{M}(\zeta, \eta) - \bs{G}(\zeta,\eta)} \dd \eta, &  \qquad & & h_2(\zeta) & :=  - \int_0^{\eta_*} \avg{\Im \bs{G}(\zeta,\eta)} \dd \eta, \\ 
h_3(\zeta) & := \, \frac{1}{4n} \sum_{i  \in \llbracket 2n \rrbracket} \log \bigg( 1 + \frac{\lambda_i(\zeta)^2}{T^2} \bigg) - \log \bigg( 1 + \frac{1}{T}\bigg) , & \qquad & & 
 h_4(\zeta) & :=  \phantom{-} \int_0^{\eta_*} \avg{\Im \bs{M}(\zeta,\eta)} \dd \eta. 
\end{aligned} \] 
Using \eqref{eq:linear_stat_log_det_H_zeta}, \eqref{eq:log_det_int_im_G} and $\int_0^T \frac{1}{1 + \eta} \dd \eta = \log( 1 + T)$, it is easy to see  
that $h(\zeta) = h_1(\zeta) + h_2(\zeta) + h_3(\zeta) + h_4(\zeta)$. 

Next, we establish individual estimates on $h_1$, \ldots, $h_4$ which hold with very high probability. 
We get $\abs{h_1(\zeta)} \leq 2 n^{-1 + K \delta}$ from \eqref{eq:global_law_averaged} as well as a union bound and a continuity argument in $\eta$. 
To estimate $h_2$, we write $\lambda_j \equiv \lambda_j(\zeta)$ and compute
\[ - h_2(\zeta) = \frac{1}{4n} \sum_{j \in \llbracket 2n \rrbracket} \log \bigg( 1 + \frac{\eta_*^{2}}{ \lambda_j^2} \bigg).  \] 
In the following, we will decompose the sum into two regimes, $\abs{\lambda_j} < \eta_*^{1/2}$ and $\abs{\lambda_j} \geq \eta_*^{1/2}$, and estimate each regime separately. 
For the first regime, Assumption~\ref{assum:A4} and Lemma~\ref{lem:number_small_singular_values_global} yield 
\begin{equation} \label{eq:bound_h_2_first_regime} 
 \frac{1}{4n} \sum_{\abs{\lambda_j} < \eta_*^{1/2}} \log \bigg( 1 + \frac{\eta_*^2}{ \lambda_j^2} \bigg) \leq \frac{ C(\log n + \abs{\log \min_{j \in \llbracket 2n \rrbracket} \abs{\lambda_j}} )}{n} 
\# \big\{ j \in \llbracket 2n \rrbracket \colon \abs{\lambda_j} \leq \eta_* \big\} \leq n^\eps \eta_* 
\end{equation}
with very high probability for all $\eps>0$ small enough. 
In the remaining regime,  $\log(1 + x) \leq x$ yields
\begin{equation} \label{eq:bound_h_2_last_regime} 
 \frac{1}{4n} \sum_{\abs{\lambda_j} \geq \eta_*^{1/2}} \log \bigg( 1 + \frac{\eta_*^2}{ \lambda_j^2} \bigg) \leq \frac{1}{4n} \sum_{\abs{\lambda_j} \geq \eta_*^{1/2}} \log ( 1 + \eta_*) \leq \eta_*.  
\end{equation}
By slightly shrinking $\delta$, these estimates imply $\abs{h_2(\zeta)} \leq n^{-\delta}$. 
For $h_3$, we obtain 
\[ \abs{h_3(\zeta)} \leq \frac{1}{4nT^2} \tr (\bs{H}_\zeta)^2  + T^{-1} = \frac{1}{2 n T^2} \sum_{i,j=1}^n (\overline{x_{ji}} - \bar \zeta \delta_{ji}) (x_{ij} - \zeta \delta_{ij}) + T^{-1} \lesssim n^\eps T^{-2} ( 1 + \abs{\zeta}^2) + T^{-1}, \] 
where we used $\log(1+ x) \leq x$ in the first step and $\abs{x_{ij}} \prec n^{-1/2}$ due to Assumption~\ref{assum:A1} in the last step. Since $\zeta \in \Omega$ we conclude $\abs{h_3(\zeta)} \lesssim n^{-10}$ with very high probability. 
Since $\norm{\bs{M}} \lesssim 1$ due to \eqref{eq:norm_M} we have $ 0 \leq h_4 (\zeta) \lesssim \eta_*$. 
Hence, the proof of \eqref{eq:F_estimate_fixed_zeta} is complete. 

Since $m$ is at most of polynomial order in $n$, a union bound over $\xi_1, \ldots, \xi_m$ and \eqref{eq:F_estimate_fixed_zeta} 
yield 
\begin{equation} \label{eq:estimate_sum_F_xi_i} 
 \frac{1}{m} \sum_{i=1}^m \abs{F(\xi_i)} \leq \frac{n^{-\delta}}{m} \sum_{i=1}^m \abs{\Delta f(\xi_i)} \leq n^{-\delta} \norm{\Delta f}_{\rm{L}^1} + n^{-D - 10} \norm{\Delta f}_{\rm{L}^{1+\beta}} 
\end{equation}
with very high probability. 
Here, we applied Lemma~\ref{lem:sampling} with $a = 1 + \beta$ in the last step and used $\supp f \subseteq \Omega$. 

Finally, we combine the relation \eqref{eq:decomposition}, the estimates \eqref{eq:error_term_sampling} and \eqref{eq:estimate_sum_F_xi_i} as well as 
$\norm{\Delta f}_{\rm{L}^{1+\beta}} \leq n^{D} \norm{\Delta f}_{\rm{L}^1}$ 
and obtain \eqref{eq:quantitative_global_law_for_X}. 
This completes the proof of Theorem~\ref{thm:global_law_for_X}. 
\end{Proof}

\section{Dyson equation and its stability}
\label{Sec:Dyson equation and its stability}

In this section, we analyse the solution $\bs{M}$ to the matrix Dyson equation \eqref{eq:MDE} and its stability against perturbations $\bs{D}$, i.e.\ we control the solution $\bs{G}(\bs{D})$ of a perturbed version of the MDE (see \eqref{Perturbed MDE} below) such that $\bs{G}(\bs{0})=\bs{M}$. 
 These results are the core of this article as they will be the basis of the proofs of Theorem~\ref{thr:Density} and Lemma~\ref{lem:sigma_as_dist_derivative} about the properties of $\sigma$ as well as the local law for $\bs{H}_\zeta$ (cf.\ Theorem~\ref{thm:local_law_H_zeta} below). 

The matrix Dyson equation and its stability have been analysed in \cite{AjankiCorrelated,Altshape}. 
However, their main regularity and stability results impose the flatness condition (see \cite[equation (2.7)]{AjankiCorrelated} and \eqref{S Flatness} below) on the self-energy operator  ${\scr{S}}$. 
This condition is not satisfied by $\scr{S}$ as defined in \eqref{definition scr S and Z}. In fact, the special structure of  $\scr{S}$, originating from the zero blocks on the diagonal of $\bs{H}_\zeta$, poses significant challenges since it leads to an instability in the Dyson equation \eqref{eq:MDE} 
 which was not present in \cite{AjankiCorrelated,Altshape}.   
Dealing with this instability is the main purpose of this section.

In \cite{Altcirc}, a similar instability was analysed, but in the  simpler setup of a random matrix with independent entries. 
This setup results in a vector-valued Dyson equation whose formulation on the commutative algebra $\C^{2n}$ with entry wise multiplication simplifies the analysis  compared to the present article. 
In particular, in the commutative setting of \cite{Altcirc} the MDE was formulated on the entire algebra $\C^{2n}$ and the contribution to the error term  in the unstable direction determined to be sufficiently small to cancel the instability in the $\eta \to 0$ limit. The corresponding algebraic manipulations are considerably harder in the non-commutative space  $\C^{2n \times 2n}$. Therefore, we develop a different strategy in the present work. Here, we identify a stable manifold $\Gamma:=\bs{G}^{-1}[\bs{E}_-^\perp]\subset \C^{2n \times 2n}$, defined as the preimage of a linear hyperspace $\bs{E}_-^\perp \subset \C^{2n \times 2n}$ under the solution map $\bs{G}$ to the perturbed MDE, such that $\Gamma \ni \bs{D} \mapsto \bs{G}(\bs{D}) \in \bs{E}_-^\perp$ is stable. Then we implicitly construct a parametrisation $\bs{E}_-^\perp \ni \wt{\bs{D}} \mapsto {\bs{D}}(\wt{\bs{D}}) \in \Gamma$ of this manifold (see \eqref{Projected Perturbed MDE}  below) and rewrite the MDE directly on the codimension one subspace $\bs{E}_-^\perp$ (see \eqref{perturbed MDE bulk regime} below). In short, we remove the unstable direction from the MDE at the beginning. In addition to removing the need to trace the unstable component of the error matrix, this strategy also implies analyticity  of $\tau \mapsto V_1(\tau),\tau \mapsto V_2(\tau)$ from \eqref{Dyson equation at eta = 0} in the bulk (see Proposition~\ref{prp:Existence and uniqueness}). The ensuing analyticity of $\sigma$ from \eqref{definition of sigma} is a new result even  for matrices with independent entries. With the strategy from \cite{Altcirc} showing only smoothness already required tracking the unstable direction to all derivative orders (cf.\ \cite[proof of Proposition~2.4]{Altcirc}).

 In Subsection~\ref{Subsec:Solution} we will establish some basic properties of the solution to \eqref{Dyson equation} and hence \eqref{eq:MDE}. Then we will prove stability of the Dyson equation in the bulk of the spectrum against small perturbations in  Theorem~\ref{thr:Stability} of  Subsection~\ref{Subsec:Stability}, using an important technical lemma that will be proven in Subsection~\ref{Subsec:Resolvent control on scr L}. 

 Since \eqref{Dyson equation} is invariant under the scaling $\cal{S} \to \lambda \2\cal{S}$, $\eta \to \lambda^{1/2}\eta$, $\tau \to \lambda\2 \tau$ and $V_i \to \lambda^{-1/2}V_i$ for any $\lambda>0$, we will assume for the rest of the paper that 
\begin{equation} \label{eq:normalization_spectral_radius} 
\spradius(\cal{S})\,=\, 1\,.
\end{equation}
Furthermore, we denote the unit disk in the complex plane centred at the origin by $\D=\D_1  = \{ \zeta \in \C \colon \abs{\zeta} < 1\}$.

\subsection{Solution}
\label{Subsec:Solution}
The first result of this subsection establishes matching upper and lower bounds on the solution of \eqref{Dyson equation} in the sense of quadratic forms. For this proposition only the following flatness property of $\cal{S}$ is needed.
Due to assumptions \ref{assum:A1}--\ref{assum:A3} the operators $\cal{S}$ and $\cal{S}^*$ are both comparable to the normalised trace in the sense of quadratic forms, i.e.\  
\bels{S Flatness}{
\cal{S} A \sim \cal{S}^* A \sim \avg{A}, \qquad A \ge 0\,,
}
for any positive semi-definite $A$. In fact, the lower bound $\cal{S} A  \gtrsim \avg{A}$ is just an equivalent formulation of assumption \ref{assum:A3}, while the upper bound $\cal{S} A  \lesssim  \avg{A}$ is a simple consequence of assumptions \ref{assum:A1} and \ref{assum:A2}. The same is true for $\cal{S}^*$. 

\begin{proposition}[Behaviour of solution] 
\label{prp:Behaviour of solution}
The solution of \eqref{Dyson equation} satisfies
\bels{avg V1 = avg V2}{
\avg{V_1(\tau,\eta)}\,=\,\avg{V_2(\tau,\eta)}\,,
}
for all $\tau \ge 0$ and $\eta >0$,
as well as the bounds
\bels{bounds on V}{
V_1(\tau,\eta) \,\sim\, V_2(\tau,\eta)\,\sim\,    
\begin{cases}
 \displaystyle (1-\tau)^{1/2} + \eta^{1/3} \,, & \tau\le 1\,,\;\eta \le 1\,, \vspace{0.2cm}
\\
\displaystyle \frac{\eta}{\tau-1+\eta^{2/3}}\,, &  \tau \ge 1 \,,\;\eta \le 1\,,\vspace{0.1cm}
 \\
 \displaystyle \frac{\eta}{\eta^2+\tau} \,, &\eta \ge 1\,.
\end{cases}
}
\end{proposition}
  
 \begin{Proof} 
 Throughout this proof, we will refer to some identities from the proof of Lemma~\ref{lem:Im_M_sim_eta}.   
To see \eqref{avg V1 = avg V2}, we multiply \eqref{U identity2} with $V_1$ from the right and with $V_2$ from the left and take the normalized trace. 
 
 Now we prove \eqref{bounds on V}. First we observe that $V_1$ and $V_2$ are both comparable to their normalized traces, which coincide as we have just shown, i.e.\
 \bels{V comparable to avg V}{
 V_1 \,\sim\, V_2 \,\sim\, \avg{V_1}\,=\, \avg{V_2}\,.
 }
 This is seen directly from the two equations \eqref{Dyson equation} since the right hand sides are both comparable to  the same  multiple of the identity  due to \eqref{S Flatness} and  \eqref{avg V1 = avg V2}. 
 
 Now let $S_2$ be the unique positive definite Perron-Frobenius eigenmatrix of $\cal{S}$ with normalisation $\avg{S_2}=1$, i.e.\ $\cal{S}S_2=S_2$. Because of \eqref{S Flatness} this eigenmatrix satisfies $S_2 \sim 1$.
We take the scalar product with $S_2$ on both sides of \eqref{eq:V_2_symmetric} and  get
\begin{equation} \label{v scaling first step}
\eta +\scalar{S_2}{V_1}\,=\, 
\eta +\scalar{S_2}{\cal{S}^*V_1}
\,=\, \scalar{S_2}{( \eta +\cal{S}^*V_1)V_2(\eta+\cal{S}^*V_1)}
+\tau \scalar{S_2}{V_1}\,.
\end{equation}

Depending on whether $\tau\le 1$ or $\tau>1$ we either subtract $\tau \scalar{S_2}{V_1}$ or $\scalar{S_2}{V_1}$ on both sides of \eqref{v scaling first step} and use \eqref{S Flatness} as well as $S_2\sim 1$ to see that 
\begin{subequations}
\label{v scaling final step}
\begin{align}
\eta + (1-\tau)\avg{V_1} \,&\sim\, \avg{V_2}(\eta + \avg{V_1})^2
\,,\quad &\tau \2\le\2 1\,,
\\
\eta 
\,&\sim\,
\avg{V_2}(\eta + \avg{V_1})^2 + (\tau-1)\avg{V_1}
\,,\quad &\tau \2>\2 1\,.
\end{align}
\end{subequations}
The claim \eqref{bounds on V} is now an immediate consequence of \eqref{v scaling final step} and \eqref{V comparable to avg V}. 
\end{Proof}

As a consequence of Proposition~\ref{prp:Behaviour of solution} we can also estimate the singular values of $U$, defined in \eqref{definition of U}, from above and below.
When multiplying \eqref{V1 equation} with $\eta + \cal{S}^*V_1$ from the left and \eqref{V2 equation} with $\eta + \cal{S}V_2$ from the right we see the identities
 \bels{U identity}{
 U \,=\, V_1\frac{1}{\eta+\cal{S}^*V_1}\,=\, \frac{1}{\eta+\cal{S}V_2}V_2\,.
  }
  Furthermore, when we multiply \eqref{V1 equation} with $V_1$ from the left and \eqref{V2 equation} with $V_2$ from the right we see that
  \bels{V equation with U}{
  1=V_1(\eta + \cal{S}V_2) +\tau\2U =(\eta + \cal{S}^*V_1)V_2 +\tau\2U\,.
  }
  Multiplying \eqref{V equation with U} by $U$ and using \eqref{U identity} also reveals
 \bels{U equation}{
{U}\,=\,V_1V_2+\tau\2{U}^2 \,.
 }
 Finally, using \eqref{bounds on V} in \eqref{U identity} shows the comparison relation 
  \bels{comparison relation for U}{
 U^*U \,\sim\, \frac{1}{1+\tau^2+\eta^4}\,
 }
 uniformly for $\eta >0$ and $\tau \geq 0$.  
 For future reference we also record the identities 
 \begin{subequations}
 \label{im of MDE with U}
\begin{align}
V_1\,&=\, \eta(V_1^2+\tau{U}{U}^*)+V_1(\cal{S}V_2)V_1+\tau{U}(\cal{S}^*V_1){U}^*,
\\
V_2\,&=\, \eta(V_2^2+\tau{U}^*{U})+V_2(\cal{S}^*V_1)V_2+\tau{U}^*(\cal{S}V_2){U}\,,
\end{align}
\end{subequations}
which result from multiplying \eqref{V1 equation} from left and right by $V_1$ and \eqref{V2 equation} by $V_2$ and then using \eqref{U identity}. 
As a consequence of \eqref{bounds on V} for $\tau \ge 1$ we can extend $V_i$ continuously  to $\eta=0$ as $V_i(\tau,0)=0$. This is summarised in the following corollary whose proof is immediate from 
 the representation of $\bs{M}$ in \eqref{definition of M} and the definition of  $U$ in \eqref{definition of U}.

 \begin{corollary}[Extension outside the spectrum] 
The functions $V_1,V_2$ admit continuous extensions to $([0,1) \times (0,\infty) )\cup([1,\infty) \times [0,\infty))$, that is, to $\eta = 0$ for $\tau \geq 1$. 
These extensions are still denoted by the same symbols.  
 Similarly, the solution $\bs{M}$ of the MDE \eqref{eq:MDE} from \eqref{definition of M} can be extended continuously to $(\D \times (0,\infty)) \cup ((\C\setminus \D) \times [0,\infty))$, i.e.\ to $\eta =0$ for $\zeta \not \in \D$. We still denote the extension by $\bs{M}$. 
The extension satisfies
 \[
  \bs{M}(\zeta,0) =  \mtwo{0 &-1/\ol{\zeta}}{-1/\zeta&0} \,, \qquad \zeta \not \in \D\,.
 \]
 \end{corollary}

\begin{Proof}[Proof of Lemma~\ref{lem:bounds_M}] 
First, we get \eqref{eq:norm_M} from \eqref{definition of M}, \eqref{bounds on V} and \eqref{comparison relation for U}. 
Second, both bounds in \eqref{eq:integral_im_M} follow directly from the estimate 
\begin{equation} \label{eq:M_minus_i_eta_inverse} 
 \norm{\bs{M}(\zeta,\eta) - \ii(1 + \eta)^{-1}} \lesssim \min\{ 1,  (1 + \abs{\zeta}) \eta^{-2} \}, 
\end{equation}
which holds uniformly for $\zeta \in \C$ and $\eta >0$ and is shown next.  
Since $\norm{\bs{M}} \lesssim {1}$ by \eqref{eq:norm_M} we trivially have $\norm{\bs{M}(\zeta,\eta) - \ii (1 + \eta)^{-1}} \lesssim {1}$. 
Multiplying \eqref{eq:MDE} by $\ii\eta^{-1} \bs{M}$ and using $\norm{\bs{M}} \lesssim \eta^{-1}$ as well as $\norm{\scr{S}} \lesssim 1$ (cf. upper bound in \eqref{S Flatness})   imply $\norm{\bs{M}(\zeta,\eta) - \ii\eta^{-1}} \lesssim (1 + \abs{\zeta})\eta^{-2}$, 
i.e.\ the missing bound in \eqref{eq:M_minus_i_eta_inverse}. This completes the proof of Lemma~\ref{lem:bounds_M}.
\end{Proof}

 \subsection{Stability}
 \label{Subsec:Stability}
 In this subsection we will establish stability of the MDE \eqref{eq:MDE} and its solution against small perturbations. As indicated at the beginning of the section, \eqref{eq:MDE} has an inherent instability due to the structure of $\scr{S}$. This instability originates from a single unstable direction and implies that stability can only be expected with respect to perturbations $\bs{D}$ that take values in a manifold of codimension $1$ in $\C^{2n \times 2n}$. Through a special choice of coordinates this manifold can be mapped to the orthogonal complement of $ \bs{E}_- \in \C^{2n \times 2n}$ defined through 
  \bels{definition bs E pm}{
 \bs{E}_\pm\,:=\, \mtwo{1 & 0}{0 & \pm1}
 }
 and thus projected out.  
 
 Before we state the stability theorem we introduce a norm that is designed to prove isotropic convergence of the resolvent $\bs{G}$ from \eqref{definition resolvent G} to $\bs{M}$, i.e.\ to prove $\abs{\scalar{\bs{x}}{(\bs{G}-\bs{M})\bs{y}}} \to 0$ for any fixed vectors $\bs{x},\bs{y} \in \C^{2n}$ in a high moment sense. A similar norm was introduced in \cite{Erdos2017Correlated} for the same purpose and to match the notation to this work we introduce coefficients $\kappa_{\scr{R}}(\alpha,\beta)$ with $\alpha, \beta \in \llbracket 2n\rrbracket^2$ associated to any linear operator $\scr{R} \colon \C^{2n \times 2n} \to\C^{2n \times 2n}$ via
 \begin{equation}\label{R and kappa R}
( \scr{R}\bs{R})_{ad} = \sum_{b,c \in \llbracket 2n\rrbracket}\msp{-7} \kappa_{\scr{R}}(ab,cd) \2r_{bc}  \,.
 \end{equation}
 Through this one to one correspondence between $\scr{R}$ and $\kappa_{\scr{R}}$ we define $\kappa_{\rm{c}}:=\kappa_{\scr{S}}$. 
 We also recall the following notation from \cite{Erdos2017Correlated}. For an expression $f_{a_1a_2\dots a_k}$ with indices $a_1,\dots, a_k$ we write $f_{\bs{x}a_2\dots a_k} = \sum_a x_a f_{a a_2\dots a_k}$ if an index is  averaged against a vector $\bs{x}$, and similarly if more than one index is averaged.  We also write $f_{\1\cdot\, a_2\dots a_k}$ for the vector $(f_{a a_2\dots a_k})_a$. In particular, $\bs{A}_{\bs{x}\bs{y}} = \sum_{i,j} x_i y_j a_{ij}$ and $\bs{A}\bs{x} = \bs{A}_{\1\cdot\,\bs{x} }$.

 Let us now fix two deterministic vectors $\bs{x},\bs{y} \in \C^{2n}$ and $K \in \N$. Then for fixed $\eta$ and $\zeta$ writing $\bs{M} =\bs{M}(\zeta,\eta)$ we recursively define the sets of vectors
 \bes{
 I_0 &:= \{\bs{x},\bs{y}\} \cup\{\bs{e}_i : i \in \llbracket 2n\rrbracket\}\,,
 \\
I_{k+1} &:= I_{k}\cup\{\bs{M}_{\bs{u}\1\cdot},\bs{Z}_{\bs{u}\1\cdot},(\scr{S} \bs{M})_{\bs{u}\1\cdot},(\ii \1 \eta\1\bs{1} +\scr{S}\bs{M})^{-1}_{\bs{u}\1\cdot}: \bs{u} \in I_k\}
\\
&\msp{45}\cup\{\kappa_{\rm{c}}(\bs{u} i,j\1\cdot\1): \bs{u} \in I_k, i,j \in \llbracket 2n\rrbracket\}\,.
 }
 Here, $\bs{e}_a$ denotes the $a$-th standard basis vector in $\C^{2n}$.  
 The $\norm{\1\cdot\1}_*$-norm is then defined as
 \bels{definition * norm}{
 \norm{\bs{A}}_*:=  \norm{\bs{A}}_*^{K, \bs{x},\bs{y}}:=\sum_{0\le k<K} n^{-k/2K} \norm{\bs{A}}_{I_k} +n^{-1/2}\max_{ \bs{u} \in I_K}\frac{\norm{\bs{A}_{\1\cdot\1 \bs{u}}}}{\norm{\bs{u}}} \,,\qquad 
 \norm{\bs{A}}_I:=\max_{\bs{u},\bs{v} \in I}\frac{\abs{\bs{A}_{\bs{u}\bs{v}}}}{\norm{\bs{u}}\norm{\bs{v}}}\,.
 }
 The definition of $\norm{\2\cdot\2}_\ast$ is chosen such that 
the arguments from \cite{Erdos2017Correlated} can be followed directly in the proof of Theorem~\ref{thm:local_law_H_zeta} below.
The norm is dominated by the standard operator norm, $\norm{\bs{A}}_*\le 2\norm{\bs{A}}$ and
by construction  and \ref{assum:A2}  it satisfies 
 \bels{scr S bounds}{
\norm{\scr S}_{\ast \to \norm{\1\cdot\1}}\lesssim 1\,, \qquad \norm{\bs{R}\bs{A}}_{\ast}\lesssim n^{1/2K}\norm{\bs{R}}\norm{\bs{A}}_\ast\,, \qquad \norm{(\scr{S}\bs{A})\bs{B}}_{\ast}\lesssim n^{1/2K}\norm{\bs{A}}_\ast\norm{\bs{B}}_\ast,
}
for all $\bs{A},\bs{B} \in \C^{2n \times 2n}$ and $\bs{R} \in \{\bs{M},\bs{Z},\scr{S}\bs{M},(\ii \1 \eta\1\bs{1} +\scr{S}\bs{M})^{-1}\}$.
The bounds from \eqref{scr S bounds} follow exactly as (73), (70b) and (70a) in \cite{Erdos2017Correlated}  from \ref{assum:A2}.  

Now we present our main stability theorem. It states that when \eqref{eq:MDE} is properly rewritten and restricted to $\bs{E}_-^\perp$ it is stable against small perturbations.  

\begin{theorem}[Stability] 
\label{thr:Stability}
For any sufficiently small  $\delta>0$ (depending on model parameters) and any $\eta \in (0,\delta^3)$, $\zeta \in \C$ with $\abs{\zeta}^2 \le 1-\delta$
 there is a unique function
\[
\bs{G}:  D\times B_1\to B_2 \,,
\]
such that $\bs{G}=\bs{G}(\zeta_1,\zeta_2,\wt{\eta},{\bs{D}})$ satisfies the equation
\bels{perturbed MDE bulk regime}{
(\ii \1\wt \eta\1\bs{1} +\scr{S} \bs{G})\pbb{
 \bs{G}+\frac{1}{\ii \1\wt \eta\1\bs{1} +\bs{Z}( \zeta_1, \zeta_2)+\scr{S}\bs{G}}
 }+{\bs{D}}\,=\,\bs{0}\,.
}
Here, $D$ is a neighbourhood of $(\zeta, \ol{\zeta}, \eta)$ in $\C^3$, $B_1$ a neighbourhood of $\bs{0}$ in $\C^{n \times n}\cap \bs{E}_-^\perp $ and $B_2$ a neighbourhood of $\bs{M}=\bs{M}(\zeta, \eta)$ in $\C^{n \times n} \cap\bs{E}_-^\perp $.
For $D,B_1$ and $B_2$ we have the following choices. Either
\bels{norm neighbourhoods of stability}{
D:= (\zeta,\ol{\zeta},\eta)+(\D_{c_1})^3\,, \qquad B_1:= B^{\norm{\1\cdot\1}}_{c_1}(\bs{0})\cap \bs{E}_-^\perp \,,
\qquad B_2:=B^{\norm{\1\cdot\1}}_{c_2}(\bs{M})\cap \bs{E}_-^\perp \,,  
}
with $c_1,c_2>0$ constants, depending only on the model parameters and on $\delta$, or 
\[
D:=(\zeta,\ol{\zeta},\eta)+(\D_{n^{-4/K}})^3 \,, \qquad B_1:=B^\ast_{n^{-4/K}}(\bs{0})\cap \bs{E}_-^\perp \,,
\qquad B_2:=B^\ast_{n^{-1/K}}(\bs{M})\cap \bs{E}_-^\perp 
\]
for sufficiently large $n$ (depending on model parameters, $\delta$ and $K$). 
Here the superscripts indicate with respect to which norm the ball $B_r(\bs{A})$ of radius $r$ around $\bs{A}$ is meant. 
The function $\bs{G}$ is analytic in all variables.
\end{theorem}

\begin{Proof} We solve the implicit equation 
\[
 \scr{J}_{ \zeta_1,\zeta_2, \wt{\eta}}[\bs{G}] + {\bs{D}} = \bs{0}
\]
for $\bs{G}=\bs{G}(\zeta_1,\zeta_2, \wt{\eta},{\bs{D}})$, where
\bels{definition scr J}{
 \scr{J}_{  \zeta_1,\zeta_2, \wt{\eta}}[\bs{G}]\,&:=\,
(\ii \1\wt{\eta}\1\bs{1} +\scr{S}\bs{G})\pbb{
 \bs{G}+\frac{1}{\ii \1\wt{\eta}\1\bs{1} +\bs{Z}( \zeta_1,\zeta_2)+\scr{S}\bs{G}}
}\,.
}
Note that $\scr{J}_{ {\zeta},{\ol{\zeta}}, \eta}[\bs{M}]=\bs{0}$ due to \eqref{eq:MDE}.
We will show that $\scr{J}$ is a well-defined bounded holomorphic function on $D \times B_2$ with values in $\bs{E}_-^\perp$. In particular, we will see that
\begin{subequations}
\label{cont bound scr J}
\begin{align}
\label{norm bound scr J}
&\norm{\scr{J}_{\zeta_1,\zeta_2, \wt{\eta}}[\bs{M}+\bs{\Delta}]}\msp{6}\lesssim \abs{ \zeta_1-\zeta}+\abs{\zeta_2-\ol{\zeta}}+\abs{\wt \eta-\eta}+\norm{\bs{\Delta}}\,,
\\
\label{ast norm bound scr J}
&\norm{\scr{J}_{\zeta_1,\zeta_2, \wt{\eta}}[\bs{M}+\bs{\Delta}]}_\ast \lesssim \abs{ \zeta_1-\zeta}+\abs{\zeta_2-\ol{\zeta}}+\abs{\wt \eta-\eta}+n^{1/2K}\norm{\bs{\Delta}}_\ast,
\end{align}
\end{subequations}
where the constants hidden in the comparison relation may depend on $\delta$ and $K$ in addition to the model parameters. We will keep this convention until the end of this proof. 
The theorem then follows from the implicit function theorem, Lemma~\ref{lmm:Implicit function theorem},  and the following bound on the inverse of the derivative $\nabla\scr{J}_{ {\zeta},\ol{\zeta}, {\eta}}:\bs{E}_-^\perp \to \bs{E}_-^\perp$ evaluated at $\bs{G}=\bs{M}$:
\bels{bounds on inverse derivative}{
\norm{(\nabla|_{\bs{G}=\bs{M}}\scr{J}_{ {\zeta},\ol{\zeta}, {\eta}})^{-1}|_{\bs{E}_-^\perp}} \lesssim 1\,, \qquad  \norm{(\nabla|_{\bs{G}=\bs{M}}\scr{J}_{ {\zeta},\ol{\zeta}, {\eta}})^{-1}|_{\bs{E}_-^\perp}}_{\ast} \lesssim n^{1/K}.
}
Note that the inverse  of the  derivative in \eqref{bounds on inverse derivative} is restricted to the hyperplane $\bs{E}_-^\perp$ and the $\norm{\1\cdot\1}_\ast$-norm on $\bs{E}_-^\perp$ is simply the restriction of the $\norm{\1\cdot\1}_\ast$-norm from   \eqref{definition * norm}  on $\C^{2n \times 2n}$. 

To see that $\scr{J}_{\zeta_1,\zeta_2, \wt{\eta}}$ leaves the hyperplane $\bs{E}_-^\perp$ invariant we compute
\bels{scr J leaves E_- perp}{
\scalar{\bs{E}_-}{\scr{J}_{\zeta_1,\zeta_2, \wt{\eta}}[\bs{G}]}= \avgB{\bs{E}_-(\ii \1\wt{\eta}\1\bs{1} +\scr{S}\bs{G})\frac{1}{\ii \1\wt{\eta}\1\bs{1} +\bs{Z}( \zeta_1,\zeta_2)+\scr{S}\bs{G}}}=0\,.
}
Here we used $\avg{\bs{E}_-(\scr{S}\bs{R})\bs{R}}=0$ for any $\bs{R} \in \C^{2n \times 2n}$  and $\bs{G} \in \bs{E}_-^\perp$  in the first identity and the general fact that by the Schur complement formula
\[
\tr\mtwo{R_{11} & R_{12}}{R_{21}& R_{22}}^{-1}\mtwo{R_{11} & 0}{ 0  & -R_{22}} =0\,,
\]
for any  invertible $2\times 2$-block matrix with square blocks in the second identity.

In the remainder of the proof we verify \eqref{cont bound scr J} and \eqref{bounds on inverse derivative} and thus the assumptions of Lemma~\ref{lmm:Implicit function theorem}. In the following we will frequently use the bounds $\bs{M}^*\bs{M} \sim 1$ and $\im \bs{M} \sim -\im \bs{M}^{-1}\sim 1$ that are a consequence of Proposition~\ref{prp:Behaviour of solution} and \eqref{definition of M}.

The inequality \eqref{norm bound scr J} is immediate when $c_1$ is chosen small enough and we apply $\norm{\scr{S}}\lesssim 1$, $\norm{\bs{M}}\lesssim 1$, as well as the fact that the singular values of $\ii \eta\bs{1} + \bs{Z} + \scr{S}\bs{M} = -\bs{M}^{-1}$ are bounded form above and below. 
For \eqref{ast norm bound scr J} we in addition employ  the bounds from \eqref{scr S bounds} with $\bs{A}=\bs{B}=\bs{\Delta}$ and $\bs{R}=\scr{S}\bs{M}$  as well as $\norm{\bs{A}}_\ast \leq 2\norm{\bs{A}}$.  
We leave the details to the reader.

The remaining part of the proof is dedicated to showing \eqref{bounds on inverse derivative}. 
Differentiating \eqref{definition scr J} with respect to $\bs{G}$ reveals that 
the derivative of $\scr{J}=\scr{J}_{\zeta, \ol{\zeta}, \eta}$  evaluated at $\bs{G}=\bs{M}$ satisfies
\bels{derivative of scr J}{
 \nabla|_{\bs{G}=\bs{M}}\scr{J}
  \,=\, -\ii \1\scr{M}\scr{L}\,,
  \qquad \scr{M} \bs{R}:= \Big(\im\frac{1}{\bs{M}} \Big) \bs{R}\,, \qquad \scr{L}\bs{R} := \bs{R}-\bs{M}(\scr{S}\bs{R}) \bs{M} \,,
}
where we used $\ii\1\eta \1\bs{1}+ \scr{S}[ \bs{M}] = \ii (\eta \bs{1}+\scr{S}[\im \bs{M}])=-\ii \im (\bs{M}^{-1})$ (cf.\ \eqref{definition scr S and Z}, \eqref{eq:MDE} and \eqref{definition of M}).
By the second bound in \eqref{scr S bounds} and
because $ -  \im (\bs{M}^{-1}) \sim  \bs{1}$ by \eqref{bounds on V} we have  
\[
\norm{\scr{M}^{-1}[\bs{R}]}_\ast
\lesssim n^{1/2K}\norm{\bs{R}}_\ast\,,
\]
 and also  $\norm{\scr{M}^{-1}}\lesssim \norm{(\im(\bs{M}^{-1}))^{-1}}\lesssim 1$. Due to \eqref{derivative of scr J}, in order to show \eqref{bounds on inverse derivative},  it therefore suffices to establish bounds on the inverse of the \emph{stability operator} $\scr{L}$, namely
\bels{bound on L inv}{
\norm{\scr{L}^{-1}|_{\scr{M}^{-1}\bs{E}_-^\perp}}\lesssim 1\,, \qquad  \norm{\scr{L}^{-1}|_{\scr{M}^{-1}\bs{E}_-^\perp}}_\ast\lesssim n^{1/2K}\,,
}
where the inverse is understood to be restricted to $\scr{M}^{-1}\bs{E}_-^\perp$ and the $\norm{\1\cdot\1}_\ast$-norm on the hyperplane $\scr{M}^{-1}\bs{E}_-^\perp$ is simply the restriction of the  $\norm{\1\cdot\1}_\ast$-norm on $\C^{2n \times 2n}$. 
The bounds \eqref{bound on L inv} are a consequence of the following three lemmas.

\begin{lemma}[Resolvent control for $\scr{L}$] \label{lmm:Resolvent control for scr L}
Let $\scr{L}$ be defined as in \eqref{derivative of scr J}. 
 For any sufficiently small  $\delta>0$ (depending on model parameters) there is a constant $\eps \sim_\delta 1$ such that uniformly in $\eta \in (0,\delta^3)$ and $\zeta \in \C$ with $\abs{\zeta}^2 \le 1-\delta$ we have the resolvent bound
\bels{resolvent bound for scr L}{
\sup
\cB{ \norm{(\scr{L}-\xi)^{-1}}_{\rm{hs}}\,:\; \xi \in \C\,,\; \xi \not \in 
 \p{2+\D_{\eps}}\cup\p{1+\D_{1-\eps}}\cup \D_\eps
}\,\lesssim_\delta\,1 \,.
}
Furthermore, the $\eps$-ball around zero contains a single isolated eigenvalue $\wh{\lambda}\ne 0$ of $\scr{L}$, i.e.\
\bels{nondegeneracy for scr L}{
\D_\eps\cap \spec(\scr{L})\,=\, \cb{\wh{\lambda}}\,,\quad \absb{\wh{\lambda}}\,\lesssim_\delta\, \eta\,,\quad
 \dim \rm{ker} \pb{\scr{L}-\wh{\lambda}}^2\,=\, 1\,.
}
Approximate right and left eigenvectors corresponding to this isolated eigenvalue of $\scr{L}$ are given by the identities
\begin{subequations}
\label{left and right eigenvector for scr L}
\begin{align}
\label{right eigenvector for scr L}
\scr{L}[\bs{E}_-\im \bs{M}]\,&=\,  \ord_{\norm{\2\cdot\2}}(\eta)\,,
\\
\label{left eigenvector for scr L}
  \scr{L}^*[\bs{E}_-\im (\bs{M}^{-1})]\,&=\,-\eta \2\bs{E}_-\,,
\end{align}
which are valid globally for $\eta >0$ and $\zeta \in \C$. 
\end{subequations}
\end{lemma}
 \begin{lemma}[Smoothing lemma] \label{lmm:Smoothing lemma} Let $\C^d$ be equipped with two norms $\norm{\2\cdot\2}_{\#}$, $\norm{\2\cdot\2}_{+}$ and $B \in \C^{d \times d}$ with 
\[
\norm{B}_{\# } + \norm{B}_{\# \to +}\norm{B}_{+ \to \#}\,\le\, C\,,
\] 
for some constant $C>0$. Then for $\xi \not \in \spec(B) \cup \{ 0 \}$ we have 
\[
\norm{(B-\xi)^{-1}}_{\# }\le\, \frac{1}{\abs{\xi}}+ \frac{C}{\abs{\xi}^2}\pb{1+ \norm{(B-\xi)^{-1}}_{+ }}\,.
\]
\end{lemma}
\begin{lemma}[Twist lemma]
\label{lmm:Twist lemma}
Let $\C^d$ be equipped with a scalar product $\scalar{\2\cdot\2}{\2\cdot\2}$ and a norm $\norm{\2\cdot\2}_{\#}$ (not necessarily induced by the scalar product), $\eps \in (0,1)$ and $A \in \C^{d \times d}$ such that $\ol{\D}_\eps \cap \spec A = \{ \alpha \}$. 
We assume that $\alpha$ is a non-degenerate eigenvalue of $A$ and $A a = \alpha a$ for 
some $a \in \C^d$ with $\norm{a}_\# = 1$. 
Let
\bels{Definition of projection P}{
P\,:=\,- \frac{1}{2\pi \ii}\oint_{\partial \D_{\eps}} \frac{\dd \zeta}{A-\zeta}\,=\, \scalar{p}{\2\cdot\2}\2a\,,
}
with some $p \in \C^d$
be the corresponding spectral projection and $b \in \C^d$ a vector such that
\bels{assumptions on orthogonality for a and b}{
\abs{\scalar{a}{b}}\,\ge\, 2\2 \eps\,,\qquad
\abs{\scalar{b}{w}}\,\le\, \norm{w}_{\#}\,,\qquad \forall\;w \in \C^d\,.
}
Suppose that $A$ has a bounded inverse on the range of $1-P$, i.e.\
\bels{assumptions on lower bound of A}{
\norm{Aw}_{\#}\,\ge\, \norm{w}_{\#}\,,\qquad \forall\; w \perp p\,.
}
Then $A$ has a bounded inverse when restricted to $b^\perp$, namely
\bels{Bound on A inverse on b perp}{
\norm{Aw}_{\#}\,\ge\, \frac{\eps}{3}\2\norm{w}_{\#}\,,\qquad \forall\; w \perp b\,.
}
\end{lemma}
Lemma~\ref{lmm:Resolvent control for scr L} is an important technical result that allows to apply analytic perturbation theory to the isolated eigenvalue  $\wh{\lambda}$  of the non-selfadjoint operator $\scr{L}$. Its proof is given in Subsection~\ref{Subsec:Resolvent control on scr L} below. 
The proof of Lemma~\ref{lmm:Smoothing lemma} is to simply take the $\norm{\1\cdot\1}_{\#}$-norm in the identity 
\[
\frac{1}{B-\xi}\,=\, -\frac{1}{\xi} -\frac{1}{\xi^2} B  +  \frac{1}{\xi^2} B \frac{1}{B-\xi}B\,.
\] 
The proof of Lemma~\ref{lmm:Twist lemma} is postponed to  Appendix~\ref{apx:Auxiliary results}. 

To show \eqref{bound on L inv} we use that by Lemma~\ref{lmm:Resolvent control for scr L} the spectral projection $\scr{P}_{\!\scr{L}}$ corresponding to the isolated eigenvalue $\wh{\lambda}$ of $\scr{L}$ close to zero has rank one and thus the form
\[
\scr{P}_{\!\scr{L}} = - \lim_{\gamma \downarrow 0} \frac{1}{2\pi \ii}\oint_{\partial \D_\gamma} \frac{\dd \xi}{\scr{L}-\wh{\lambda}-\xi}\,=\, \frac{\scalar{\bs{L}_l}{\2\cdot\2}}{\scalar{\bs{L}_l}{\bs{L}_r}}\2\bs{L}_r\,,
\]
where $(\scr{L}-\wh{\lambda})\bs{L}_r = (\scr{L}-{\wh{\lambda}})^*\bs{L}_l = \bs{0}$, i.e.\ $\bs{L}_r$ and $\bs{L}_l$ are the  unique (up to normalisation)  corresponding right and left eigenvectors of $\scr{L}$, respectively.  

Now we extend the resolvent control \eqref{resolvent bound for scr L} from the $\norm{\1\cdot\1}_{\rm{hs}}$-norm to the norms $\norm{\1\cdot\1}$ and $\norm{\1\cdot\1}_{\ast}$ with the help of Lemma~\ref{lmm:Smoothing lemma} applied to the choice $B=\rm{Id}-\scr{L}$. This is possible because
\[
\norm{B}_{\#\to \#}+\norm{B}_{\rm{hs}\to \#} \norm{B}_{\#\to \rm{hs}}  \lesssim \norm{\bs{M}}^2\norm{\scr S}_{\# \to \norm{\1\cdot\1}}
+\norm{\bs{M}}^4 \norm{\scr S}_{\rm{hs} \to \norm{\1\cdot\1}}\norm{\scr S}_{\# \to \norm{\1\cdot\1}}\,, \quad \#=\ast,\norm{\1\cdot\1}\,.
\]
and $\norm{\scr S}_{\rm{hs} \to \norm{\1\cdot\1}}+\norm{\scr S}_{\# \to \norm{\1\cdot\1}}\lesssim 1$.   In particular, we may use analytic perturbation theory in the $\norm{\1\cdot\1}$-norm and  find  
\bels{expansion of Lr and Ll}{
\bs{L}_r = \bs{E}_-\im \bs{M} +\ord_{\norm{\1\cdot\1}}(\eta)\,, \qquad \bs{L}_l = \bs{E}_-\im (\bs{M}^{-1}) +\ord_{\norm{\1\cdot\1}}(\eta)
}
according to \eqref{left and right eigenvector for scr L}. 
Applying Lemma~\ref{lmm:Twist lemma} with  the choices 
\[
A= C\2\scr{L}\,, \qquad a= \frac{\bs{L}_r}{\norm{\bs{L}_r}_{\#}}\,, \qquad p = \frac{\norm{\bs{L}_r}_{\#}}{\scalar{\bs{L}_l}{\bs{L}_r}}  \bs{L}_l \,, \qquad b = c\2 \frac{\bs{E}_-\im (\bs{M}^{-1})}{\norm{\im (\bs{M}^{-1})}} 
\]
shows  the invertibility of $\scr{L}$ on $\scr{M}^{-1}\bs{E}_-^\perp = (\bs{E}_-\im (\bs{M}^{-1}))^{\perp}$ in the $\norm{\1\cdot\1}_\#$-norm. Here, the positive constants $c$ and $C$ are chosen sufficiently small and large, respectively, in order to ensure the assumptions \eqref{assumptions on orthogonality for a and b} and \eqref{assumptions on lower bound of A} of  Lemma~\ref{lmm:Twist lemma}. 
In case of the $\# = \norm{\1\cdot\1}$ we have $c \sim C \sim 1$ and in  the  $\#=* $ case $c \sim n^{-1/2K}$ and $C\sim 1$. The expansion \eqref{expansion of Lr and Ll} is used to ensure that indeed $\abs{\scalar{a}{b}}\,\ge\, 2\2 \eps$ as required in \eqref{assumptions on orthogonality for a and b} and \eqref{assumptions on lower bound of A} follows from the resolvent control on $\scr{L}$ in $\#$-norm.
\end{Proof}
 
  \begin{corollary}[Perturbations] \label{cor:perturbations} 
Let  $\eta >0$ and $\zeta \in \C$ with $\eta + \abs{\abs{\zeta}-1} \ge \delta$ for some fixed $\delta>0$.  
For any $\bs{D} \in \C^{2n \times 2n}$ and $\bs{G} \in \bs{E}_-^\perp$ such that $\norm{\bs{G}-\bs{M}}_* + \norm{\bs{D}}_* \le n^{-7/K}$ (respectively $\norm{\bs{G}-\bs{M}} + \norm{\bs{D}} \ll 1$) that satisfy  the perturbed Dyson equation
 \bels{Perturbed MDE}{
-\bs{1}\,=\, (\ii \1\eta\1\bs{1} +\bs{Z}(\zeta,  \overline{\zeta} )+\scr{S}[ \bs{G}]) \bs{G}-\bs{D}\,,
 }
the matrix $\bs{G}$ is close to $\bs{M} = \bs{M}(\zeta,\eta)$ in the sense that for sufficiently large $n$ we have  
  \bels{Global stability}{
 \norm{\bs{G}- \bs{M}}_{*}\,\le\,  \frac{n^{6/K}}{1 + \eta}  \norm{\bs{D}}_{*} \quad \bigg(\text{respectively} \;  \norm{\bs{G}- \bs{M}}\,\lesssim\, \frac{1}{1+\eta}\norm{\bs{D}}
\bigg)\,.
 }
 \end{corollary}
 
We also introduce the Matrix Dyson equation with general spectral parameter given by 
\begin{equation} \label{eq:mde_general_spectral_parameter} 
 -\bs{M}^{-1} = z \bs{1} + \bs{Z} + \scr S \bs{M} 
\end{equation}
with $z \in \bb{H} := \{ w \in \C \colon \Im w >0 \}$ as well as   
$\scr S$ and $\bs{Z} = \bs{Z}(\zeta,\bar \zeta)$ from \eqref{definition scr S and Z} with $\zeta \in \C$. 
There is a unique solution $\bs{M} = \bs{M}(\zeta,z)$ to \eqref{eq:mde_general_spectral_parameter} 
under the constraint $\Im \bs{M} \geq 0$ \cite{Helton01012007}. 
Note that \eqref{eq:mde_general_spectral_parameter} is the counterpart of \eqref{eq:MDE}, 
where the special spectral paramter $\ii \eta \in \bb{H}$ is replaced by a general $z \in \bb{H}$. 
In particular, both solutions agree for $z = \ii \eta$.   
To \eqref{eq:mde_general_spectral_parameter}, we associate the \emph{self-consistent density of states} 
$\rho_\zeta$ of $\bs{H}_\zeta$ defined as the unique probability measure on $\R$ whose Stieltjes transform is given 
by 
\begin{equation} \label{eq:def_rho_zeta} 
 \avg{\bs{M}(\zeta,z)} = \int_{\R} \frac{\rho_\zeta(\dd \omega)}{\omega -z} 
\end{equation}
for any $z \in \bb{H}$. 

The support of $\rho_\zeta$ is called the \emph{self-consistent spectrum} of $\bs{H}_\zeta$.
By Corollary~\ref{crl:scDOS inside disk} below, $\supp \rho_\zeta$ is bounded away from 
zero for any $\zeta \notin \D_{1 + \delta}$ due to \ref{assum:A1} -- \ref{assum:A3} and our normalisation~\eqref{eq:normalization_spectral_radius}.

 \begin{Proof}
We first consider the regime $\max\{\eta , \abs{\zeta} -1\} \ge \delta$. As we will see this  corresponds to the regime away from the self-consistent spectrum and can be covered by combining existing results.  
If $\abs{\zeta} \geq 1 + \delta$ then Corollary~\ref{crl:scDOS inside disk} below implies 
that $\dist(0,\supp \rho_\zeta) \gtrsim_\delta 1$. 
Therefore, $\dist(\ii \eta, \supp \rho_\zeta) \gtrsim_\delta 1$ in the regime 
under consideration (this estimate is trivial if $\eta \geq\delta$), i.e.\ this regime 
is away from the self-consistent spectrum $\supp \rho_\zeta$. 
In particular, we may apply  Lemma~\ref{lem:stability_outside} below and \cite[eq.~(70c)]{Erdos2017Correlated}  
 to \cite[eq.~(69)]{Erdos2017Correlated} and conclude that \eqref{Global stability} holds
if $\max\{ \eta, \abs{\zeta} - 1 \} \geq \delta$. 
 
 The remaining regime $1-\abs{\zeta} \ge \delta$ and $\eta <\delta$ is treated using Theorem~\ref{thr:Stability}.
  In this case we rewrite \eqref{Perturbed MDE} in the form 
 \bels{Projected Perturbed MDE}{
  \scr{J}_{ \zeta,  \ol{\zeta},  \eta}[\bs{G}]+ \wt{\bs{D}}\,=\, \bs{0}\,, \qquad \wt{\bs{D}}:=-
  (\ii \1\eta\1\bs{1} +\scr{S}\bs{G})\frac{1}{\ii \1\eta\1\bs{1} +\bs{Z}+\scr{S}\bs{G}}\bs{D}\,,
 }
 where $\scr{J} = \scr{J}_{\zeta, \ol{\zeta}, \eta} $ is  given in \eqref{definition scr J}. We have seen in \eqref{scr J leaves E_- perp} that $\scr{J}:\bs{E}_-^\perp \to \bs{E}_-^\perp$ and thus $\scr{J}[\bs{G}] \in \bs{E}_-^\perp$ by the assumption on $\bs{G}$. In particular, \eqref{Projected Perturbed MDE} also implies $\wt{\bs{D}}\in \bs{E}_-^\perp$. By  Theorem~\ref{thr:Stability} the claim \eqref{Global stability} now follows from 
 \bels{wt D le D}{
 \norm{\wt{\bs{D}}}_* \lesssim n^{2/K}\norm{\bs{D}}_* \quad (\text{respectively}\;  \norm{\wt{\bs{D}}} \lesssim \norm{\bs{D}})\,,
 }
 because $\bs{G}$ analytically depends on $\wt{\bs{D}}$ and thus $\norm{\bs{G}(\wt{\bs{D}})-\bs{M}}_\ast\lesssim n^{3/K}\norm{\wt{\bs{D}}}_\ast$ (respectively $\norm{\bs{G}(\wt{\bs{D}})-\bs{M}}\lesssim \norm{\wt{\bs{D}}}$). 
 
 To show \eqref{wt D le D} in case of the $\norm{\1\cdot\1}_*$-norm we use the MDE \eqref{eq:MDE} and a geometric series expansion to write $\wt{\bs{D}}$ in the form
 \[
 \wt{\bs{D}} = -\bs{D} - \bs{Z}\frac{1}{\bs{M}^{-1} -  \scr{S}\bs{\Delta}} \bs{D}   =- \bs{D} -\bs{Z} \bigg(\sum_{k=0}^K\bs{M}( (\scr{S}\bs{\Delta})\bs{M})^k + 
\frac{1}{\bs{M}^{-1} -  \scr{S}\bs{\Delta}}( (\scr{S}\bs{\Delta})\bs{M})^{K+1}\bigg) \bs{D}\,,
 \]
 where $\bs{\Delta}= \bs{G}-\bs{M}$. Applying \eqref{scr S bounds} we take the $\norm{\2\cdot\2}_*$-norm on both sides and estimate 
  \[
 \norm{\wt{\bs{D}}}_* \lesssim \norm{\bs{D}}_*+\sum_{k=0}^{K} n^{(k+2)/K}\norm{\bs{\Delta}}_*^k\norm{\bs{D}}_*  + n^{1/2} \norm{\bs{\Delta}}_*^{K+1}\norm{\bs{D}}_*,
 \]
 where for the last summand we used that
$\norm{\bs{A}\bs{B}}_* \lesssim n^{1/2} \norm{\bs{A}} \norm{\bs{B}}_*$  for any pair of matrices $\bs{A},\bs{B}$. 
 Owing to the assumption $\norm{\bs{\Delta}}_*\le n^{-7/K}$, this  verifies \eqref{wt D le D}. 
  \end{Proof}

\subsection{Resolvent control on $\scr{L}$}
\label{Subsec:Resolvent control on scr L}

In this subsection we prove Lemma~\ref{lmm:Resolvent control for scr L} by considering a reduction $\cal{L}$ of $\scr{L}$ on the space of diagonal block matrices, or equivalently on $\C^{n \times n} \oplus \C^{n \times n}$. 
We introduce the short hand notation
\[
\scr{C}_{\bs{A}}\bs{B}\,:=\, \bs{A}\bs{B}\bs{A}\,,\qquad
\cal{C}_{A}B\,:=\, ABA\,,\qquad \cal{K}_AB\,:=\, A^*BA\,,
\]
 as well as the average and scalar product on $\C^{n \times n} \oplus \C^{n \times n}$ as
\bels{definition reduced scalar product}{
\avgbb{\vtwo{A}{B}}\,:=\, \frac{1}{2}(\avg{A}+\avg{B})\,, \qquad  \scalarbb{\vtwo{A_1}{B_1}}{\vtwo{A_2}{B_2}}\,:=\, \frac{1}{2}\pb{\scalar{A_1}{A_2}+\scalar{B_1}{B_2}}\,,
}
 for $A,B,A_1,B_1,A_2,B_2 \in \C^{n \times n}$. We will denote linear operators $\cal{A}$ on $\C^{n \times n} \oplus \C^{n \times n}$ by the block notation 
 \[
\mtwo
{\cal{A}_{11} & \cal{A}_{12} }
{\cal{A}_{21} & \cal{A}_{22}  }\vtwo{A}{B}\,:=\,\vtwo{\cal{A}_{11} A+ \cal{A}_{12}B}{\cal{A}_{21}A+ \cal{A}_{22}B}\,.
\]

We split the stability operator $\scr{L}$ into diagonal and off-diagonal contributions,
\bels{splitting stability operator}{
\scr{L}\,=\, \scr{P}^*\cal{L}\scr{P}+\rm{Id}-\scr{P}^*\scr{P} +\scr{Q}\,.
}
Here we introduced the projection and embedding operators 
\bels{definition of scr P}{
&\scr{P}: \C^{2n \times 2n} \to \C^{n \times n} \oplus \C^{n \times n}\,,\;\bs{A}=\mtwo{A_{11} & A_{12}}{A_{21} & A_{22}} \mapsto (A_{11},A_{22})\,,
\\
&\scr{P}^*: \C^{n \times n} \oplus \C^{n \times n} \to \C^{2n \times 2n}\,,\; (A_1,A_2) \mapsto \mtwo{A_1 & 0}{0 & A_2}\,,
}
the \emph{reduced stability operator}
 
\bels{definition cal L}{
\cal{L}\,:=\, 1-\scr{P} \2\scr{C}_{\bs{M}}\scr{S}\scr{P}^*=\mtwo{1-\tau\2\cal{K}_{U^*}\cal{S}^* & \cal{C}_{V_1}\cal{S}}
{ \cal{C}_{V_2}\cal{S}^* &1- \tau\2\cal{K}_{U}\cal{S}}, \qquad  \cal{L}[E_-^\perp] = \vtwo{\eta + \cal{S}V_2}{-\eta- \cal{S}^*V_1}^\perp\,,
}
with $E_- \in \C^{n \times n} \oplus  \C^{n \times n}$ defined in analogy to \eqref{definition bs E pm} through
\bels{definition E pm}{
E_\pm := \vtwo{1}{\pm 1}\,, 
}
and the offdiagonal contribution $\scr{Q}:\C^{2n \times 2n} \to\C^{2n \times 2n} $ to the stability operator,
\bels{definition of scr Q}{
\scr{Q}\bs{A}\,:=\, \mtwo{0& \ii\2\zeta\2(V_1(\cal{S}A_{22} )  U+U(\cal{S}^*A_{11})V_2)}{\ii\2\ol{\zeta}\2(U^*(\cal{S}A_{22})V_1+V_2(\cal{S}^*A_{11})U^*)&0}\,.
}
Similarly to \eqref{definition E pm} we also write
\[
V_{\pm}\,:=\, \vtwo{V_1}{\pm V_2}\,.
\]
On the level of the reduced stability operator the result analogous to Lemma~\ref{lmm:Resolvent control for scr L} is the following statement. 

\begin{lemma}[Resolvent control for $\cal{L}$] \label{lmm:Resolvent control for cal L} For any sufficiently small  $\delta>0$ (depending on model parameters) there is a constant $\eps \sim_\delta 1$ such that uniformly in $\eta \in (0,\delta^3)$ and $\zeta \in \C$ with $\abs{\zeta}^2 \le 1-\delta$ we have the resolvent bound 
\bels{resolvent bound for L}{
\sup
\cB{ \norm{(\cal{L}-\xi)^{-1}}_{\rm{hs}}\,:\; \xi \in \C\,,\; \xi \not \in 
\p{2+\D_\eps} \cup\p{1+\D_{1-\eps}}\cup \D_\eps
}\,\lesssim_\delta\,1 \,.
}
Furthermore, the $\eps$-ball around zero contains a single isolated eigenvalue $\lambda\ne 0$ of $\cal{L}$, i.e.\ 
\bels{nondegeneracy for L}{
\D_\eps \cap \spec(\cal{L})\,=\, \{\lambda\}\,,\quad \abs{\lambda}\,\lesssim_\delta\, \eta\,,\quad
 \dim \rm{ker} (\cal{L}-\lambda)^2\,=\, 1\,.
}
Approximate right and left eigenvectors corresponding to this isolated eigenvalue of $\cal{L}$ are given by the identities
\begin{subequations}
\label{left and right eigenvector for L}
\begin{align}
\label{right eigenvector for L}
\cal{L}V_-\,&=\,  \ord_{\norm{\2\cdot\2}}(\eta)\,,
\\
\label{left eigenvector for L}
  \cal{L}^*{\vtwo{\eta +\cal{S}V_2}{-\eta-\cal{S}^*V_1}}\,&=\, \eta \2E_-\,.
\end{align}
\end{subequations}
which are valid globally for $\eta >0$ and $\zeta \in \C$. 
\end{lemma}
The proof of Lemma~\ref{lmm:Resolvent control for cal L} requires some preparation. But first we will see how the lemma is used to establish Lemma~\ref{lmm:Resolvent control for scr L}. 
\begin{Proof}[Proof of Lemma~\ref{lmm:Resolvent control for scr L}] The identities \eqref{left and right eigenvector for scr L} follow from \eqref{left and right eigenvector for L} because  the off-diagonal component $\scr{Q}$ of $\scr{L}$ from \eqref{definition of scr Q} almost vanishes on the approximate eigenvector. More precisely,  $\scr{Q}[ 
\bs{E}_-  \im \bs{M}]= \ord(\eta)$ and  $\scr{Q}^*[\bs{E}_- \Im (\bs{M}^{-1})] = \bs{0}$ due to the definition of $U$ in \eqref{U identity}. 

For $t \in [0,1]$ we consider an interpolation $\scr{L}_t := \scr{L} -t\scr{Q} $ that removes the off-diagonal contribution.
With the help of \eqref{resolvent bound for L} we now establish the lower bound 
\bes{
\norm{(\scr{L}_t-\xi \1\rm{Id})\bs{R}}_{\rm{hs}}&= \norm{(\cal{L}-\xi )\scr{P}\bs{R}}_{\rm{hs}} + \norm{((1-t)\scr{Q}+(1-\xi) (\rm{Id}-\scr{P}^*\scr{P}))\bs{R}}_{\rm{hs}} 
\\
&\ge \frac{\norm{\scr{P}\bs{R}}_{\rm{hs}}}{\norm{(\cal{L}-\xi )^{-1}}_{\rm{hs}}}+\pb{\abs{1-\xi}\norm{(\rm{Id}-\scr{P}^*\scr{P})\bs{R}}_{\rm{hs}}-(1-t)\norm{\scr{Q}}_{\rm{hs}}\norm{ \scr{P}\bs{R}}_{\rm{hs}}}_+\gtrsim_\delta \norm{\bs{R}}_{\rm{hs}} 
}
 for any $\bs{R} \in \C^{2n \times 2n}$ and $\xi$ in the domain where the  resolvent is controlled, i.e.\ $\xi \not \in\p{2+\D_\eps}\cup \p{1+\D_{1-\eps}}\cup \D_\eps$. This finishes the proof of \eqref{resolvent bound for scr L} with the choice $t=0$. Furthermore, it shows that no eigenvalues can leave the complement of the domain where the  resolvent is controlled along the continuous interpolation. We conclude that the non-degeneracy property \eqref{nondegeneracy for scr L} holds if it can be established for $\scr{L}_1 = \scr{P}^*\cal{L}\scr{P}+\rm{Id}-\scr{P}^*\scr{P}$. But $\scr{L}_1$ leaves both, the space of diagonal and of off-diagonal block matrices, invariant and acts as $\cal{L}$ on the first and as the identity on the latter. Thus \eqref{nondegeneracy for scr L} follows from \eqref{nondegeneracy for L}. 
 
 Finally, the fact that $\wh{\lambda} \ne 0$ follows e.g.\ from the general result on the weak, i.e.\ $\eta$-dependent,  stability of the Dyson equation from Lemma~\ref{lem:stability_outside}.
\end{Proof}

To prepare the proof of Lemma~\ref{lmm:Resolvent control for cal L} we introduce some auxiliary operators. The purpose of these operators is to allow for a rewriting of the non-Hermitian reduced stability operator $\cal{L}$ in terms of Hermitian operators for which spectral information can be turned into norm bounds.

\begin{definition}
\label{def:definition of P cal T cal F and cal V} For any $\eta >0$ and $\tau \ge 0$  we define the $n \times n$-matrices
\bels{definition of P}{
P:= \frac{1}{\sqrt{V_1}}\2U\frac{1}{\sqrt{V_2}}\,,\qquad K_1\,:=\,(1+\tau P^*P)^{-1/4}
\,,\qquad K_2\,:=\, (1+\tau PP^*)^{-1/4}\,,
}
in terms of the solution $V_i=V_i(\tau,\eta)$ to \eqref{V1 equation} and $U=U(\tau,\eta)$ from \eqref{definition of U}. Furthermore, we define the linear operators $\cal{T}_{\tau,\eta }, \cal{F}_{\tau,\eta}, \cal{V}_{\tau,\eta}:  \C^{n\times n} \oplus \C^{n\times n} \to \C^{n\times n} \oplus \C^{n\times n}$ 
through
\bels{definition of cal T}{
\cal{T}_{\tau,\eta }\,:=\,
\mtwo{
-\2\cal{C}_{K_2}^{\22}
&
\tau \2\cal{C}_{K_2}\,\cal{K}_{P^*}\2\cal{C}_{K_1}
}
{
\tau\2 \cal{C}_{K_1} \,\cal{K}_{P}\2\cal{C}_{K_2}
&
-\2\cal{C}_{K_1}^{\22}
}\,,
}
as well as
\bels{definition of cal F}{
\cal{F}_{\tau,\eta}\,:=\,
\mtwo{
0
&
\cal{C}_{K_2}^{-1}\,\cal{C}_{\!\sqrt{V_1}}\2\cal{S}\2\cal{C}_{\!\sqrt{V_2}}\2\cal{C}_{K_1}^{-1}
}
{
\cal{C}_{K_1}^{-1}\,\cal{C}_{\!\sqrt{V_2}}\2\cal{S}^*\2\cal{C}_{\!\sqrt{V_1}}\2\cal{C}_{K_2}^{-1}
&
0
}\,,
}
and 
\bels{definition of cal V}{
\cal{V}_{\tau,\eta}\,:=\,
\mtwo{
\cal{C}_{K_2}\,\cal{C}_{\!\sqrt{V_1}}^{-1}
&
0
}
{
0
&
\cal{C}_{K_1}\,\cal{C}_{\!\sqrt{V_2}}^{-1}
}\,.
}
\end{definition}

The matrices from Definition~\ref{def:definition of P cal T cal F and cal V} allow to rewrite $\cal{L}$ through the formula
\bels{stability operator in terms of cal F}{
\cal{L} 
\,=\, \cal{V}^{-1}(1-\cal{T} \cal{F})\cal{V}\,.
}
The following three lemmas list important analytical properties of the operators from Definition~\ref{def:definition of P cal T cal F and cal V}.

\begin{lemma}[Properties of $P$, $K_1$, $K_2$  and $\cal{V}$] Fix $\eta >0$ and $\tau \ge 0$. 
The  matrices $K_1$, $K_2$ and  $P$ defined in \eqref{definition of P}  satisfy  the identities
\bels{eta + SV in terms of P}{
K_2^4\,=\, \sqrt{V_1} \,(\eta + \cal{S}[V_2]) \sqrt{V_1}\,,\qquad K_1^4\,=\, \sqrt{V_2}\, (\eta + \cal{S}^*[V_1]) \sqrt{V_2}\,,
}
as well as the comparison relations
\bels{comparison for 1+tau PP}{
 P P ^* \, \sim \,  P^*P\,\sim\, \frac{1}{\eta^2 +\rho^2}\,,\qquad 
 K_1^4 \, \sim \,  K_2^4\,\sim\, (1+\tau+\eta^2)\rho^2\,.
}
The operator $\cal{V}$ from \eqref{definition of cal V}  is invertible and  satisfies
\bels{bound on cal V}{
\norm{\cal{V}}_{\rm{hs}}\norm{\cal{V}^{-1}}_{\rm{hs}}\,\sim\, 1\,.
}
\end{lemma}
\begin{Proof}
The identities \eqref{eta + SV in terms of P} follow from
\bels{eta + SV in terms of P proof}{
\frac{1}{\sqrt{V_1}}\frac{1}{1+\tau PP^*}\frac{1}{\sqrt{V_1}}\,&=\, \frac{1}{V_1+\tau\2 U V_2^{-1}U^*}\,=\, 
\eta + \cal{S}[V_2]\,,
\\
\frac{1}{\sqrt{V_2}}\frac{1}{1+\tau P^*\!P}\frac{1}{\sqrt{V_2}}\,&=\, \frac{1}{V_2+\tau\2 U^* V_1^{-1}U}\,=\, 
\eta + \cal{S}^*[V_1]\,,
}
which is easily checked by inserting the definition of $P$ from \eqref{definition of P} and using \eqref{U identity} as well as the Dyson equation \eqref{Dyson equation} for $V_1$ and $V_2$. In particular, \eqref{eta + SV in terms of P proof} implies 
the  third and fourth  relation in \eqref{comparison for 1+tau PP}
by the comparison relation for $V_1$ and $V_2$ from \eqref{bounds on V}. From  these comparison relations for $1+\tau P P^*$ and  $1+\tau P^*P$ as well as  \eqref{bounds on V} the bound \eqref{bound on cal V} follows. 
The  first two relations  in \eqref{comparison for 1+tau PP}  are  immediate consequences of the definition of $P$ in \eqref{definition of P}, the identity \eqref{U identity} and \eqref{bounds on V}. 
\end{Proof}

\begin{lemma}[Properties of $\cal{F}$]
\label{lmm:Properties of cal F}
 The operator $\cal{F}$ defined in \eqref{definition of cal F} satisfies the following properties  uniformly in $\eta>0$ and $\tau \ge 0$:
\begin{enumerate}
\item It is self-adjoint with respect to the scalar product \eqref{definition reduced scalar product} and positivity preserving, i.e.\ 
\bels{cal F self-adjoint and positivity preserving}{
\cal{F}^*\,=\, \cal{F}\,,\qquad\qquad \cal{F}[\ol{\scr{C}}_+ \oplus \ol{\scr{C}}_+] \subseteq  \ol{\scr{C}}_+ \oplus \ol{\scr{C}}_+\,,
}
where $\scr{C}_+$ denotes the cone of positive definite matrices and $\ol{\scr{C}}_+$ its closure. 
\item \label{cal F Property 2} It has a positive spectral radius 
\[
\norm{\cal{F}}_{\rm{hs}}\,\sim\, \frac{1}{1+\tau+\eta^2}\,,
\]
and $\pm \norm{\cal{F}}_{\rm{hs}}$ are non-degenerate eigenvalues of $\cal{F}$ with unique corresponding eigenvectors $F_{\pm}$ of the form 
\[
F_\pm\,=\,\vtwo{F_1}{\pm F_2}\,,
\]
for some normalized ($\norm{F_1}_{\rm{hs}}=\norm{F_2}_{\rm{hs}}=1$) matrices $F_1,F_2 \in \scr{C}_+$. Both these matrices are comparable to the identity matrix
\bels{F1 and F2 comparable to 1}{
F_1 \,\sim\, 1\,,\qquad F_2\,\sim\, 1\,.
}
\item \label{cal F Property 3} The spectral gap of $\cal{F}$ is bounded away from zero, i.e.\  there exists $\eps\sim 1$ such that
\bels{spectral gap of cal F}{
\spec(\cal{F}/ \norm{\cal{F}}_{\rm{hs}})\,\subseteq\, \{-1\} \cup [-1+\eps,1-\eps]\cup\{1\}\,.
}
\item \label{cal F Property 5} The spectral radius  of $\cal{F}$ is given by the formula
\bels{formula for spectral radius of cal F}{
1-\norm{\cal{F}}_{\rm{hs}}\,=\, \frac{\scalar{F_1}{\cal{C}_{K_2}^{-1}[V_1]}+ \scalar{F_2}{\cal{C}_{K_1}^{-1}[V_2]}}{2 \scalar{F_+}{\cal{V}[V_+]} }\,\eta\,\sim\, \frac{1}{1+\tau+\eta^2}\frac{\eta}{\rho}\,.
}
\item \label{cal F Property 4} The eigenvectors $F_\pm$ satisfy
\bels{approximate eigenvectors for cal F}{
F_\pm\,=\, \frac{\cal{V}[V_\pm]}{\msp{5}\norm{\cal{V}[V_\pm]}_{\rm{hs}}\msp{-5}}+ \ord_{\rm{hs}}\pbb{\frac{1}{1+\tau+\eta^2}\frac{\eta}{\rho}} \,.
}
\end{enumerate}

\end{lemma}

\begin{Proof}
The self-adjointness of $\cal{F}$ is clear from its definition \eqref{definition of cal F} and the property of being positivity preserving is inherited from the same properties of $\cal{S}$ (cf.\ \eqref{definition cal S}). Thus \eqref{cal F self-adjoint and positivity preserving} holds true. 

Properties~\ref{cal F Property 2} and \ref{cal F Property 3} now follow from the structure
\[
\cal{F}\,=\, \mtwo{0 & \wh{\cal{F}}}{ \wh{\cal{F}}^* & 0}\,,\qquad  \wh{\cal{F}}=\cal{C}_{K_2}^{-1}\,\cal{C}_{\!\sqrt{V_1}}\2\cal{S}\2\cal{C}_{\!\sqrt{V_2}}\2\cal{C}_{K_1}^{-1}\,,
\]
given in \eqref{definition of cal F}. Thus the spectrum of $\cal{F}$ is determined by the spectrum of $\wh{\cal{F}}^*\wh{\cal{F}}$ through
\[
\spec(\cal{F})\,=\, \spec\pb{-(\wh{\cal{F}}^*\wh{\cal{F}})^{1/2}}\cup\spec\pb{(\wh{\cal{F}}^*\wh{\cal{F}})^{1/2}}. 
\]
Because of \eqref{comparison for 1+tau PP} and $V_1\sim V_2 \sim \rho$ (cf.\ \eqref{bounds on V}) the operators  $\wh{\cal{F}}^*\wh{\cal{F}}$ and $\wh{\cal{F}}\wh{\cal{F}}^*$ inherit the flatness property \eqref{S Flatness} from $\cal{S}$, i.e.,
\[
\wh{\cal{F}}^*\wh{\cal{F}}A\,\sim\, \frac{1}{1+\tau^2+\eta^4}\,\avg{A}\,,\qquad\wh{\cal{F}}\wh{\cal{F}}^*A\,\sim\, \frac{1}{1+\tau^2+\eta^4}\,\avg{A}\,,\qquad \forall \; A \in \ol{\scr{C}}_+\,.
\]
Thus we can apply  \cite[Lemma~4.8]{AjankiCorrelated}  to infer 
\bels{spectrum of abs wh cal F}{
\spec\pb{\p{\wh{\cal{F}}^*\wh{\cal{F}}}^{1/2}/ \norm{\cal{F}}_{\rm{hs}}}\,=\, \spec\pb{\p{\wh{\cal{F}}\wh{\cal{F}}^*}^{1/2}/ \norm{\cal{F}}_{\rm{hs}}} \subseteq [-1+\eps,1-\eps]\cup\{1\}\,,
}
where $\eps\sim 1$ is a bound on the spectral gap and 
\[
\norm{\wh{\cal{F}}^*\wh{\cal{F}}}_{\rm{hs}}\,=\, \norm{\cal{F}}_{\rm{hs}}^2\,\sim\,  \frac{1}{1+\tau^2+\eta^4}\,.
\]
According to the same lemma the eigenvalue $1$ in \eqref{spectrum of abs wh cal F} is non-degenerate with corresponding normalised eigenmatrices $F_1, F_2 \in \scr{C}_+$ that satisfy \eqref{F1 and F2 comparable to 1}.
 In particular,  
\[
\wh{\cal{F}}\wh{\cal{F}}^* F_1\,=\, \norm{\cal{F}}_{\rm{hs}}^2F_1\,,\qquad \wh{\cal{F}}^*\wh{\cal{F}} F_2\,=\, \norm{\cal{F}}_{\rm{hs}}^2F_2\,. 
\]
 Therefore, $F_\pm$ are eigenvectors of $\cal{F}^2$ corresponding to $\norm{\cal{F}}_{\rm{hs}}^2$ and, consequently, $\cal{F} F_\pm =  \pm \norm{\cal{F}}_{\rm{hs}} F_\pm$.

It remains to verify Properties~\ref{cal F Property 5} and \ref{cal F Property 4}.
For this purpose we will use that $\cal{V} V_{\pm}$ are approximate eigenvectors,
\bels{Vpm approximate eigenvectors of F}{
\cal{F}\cal{V}V_\pm\,=\,\pm \cal{V} V_\pm-\eta \vtwo{\pm\2\cal{C}_{K_2}^{-1}V_1}
{\cal{C}_{K_1}^{-1}V_2 }\,.
}
Indeed, \eqref{Vpm approximate eigenvectors of F} follows from 
using the definition of $\cal{V}$ in \eqref{definition of cal V} to identify the first summand on the right hand side of 
\bels{computing cal F cal VVpm}{
\cal{F}\cal{V}V_\pm\,=\,
\vtwo{
\pm\cal{C}_{K_2}^{-1}\,\cal{C}_{\!\sqrt{V_1}}\2\cal{S}V_2
}{
\cal{C}_{K_1}^{-1}\,\cal{C}_{\!\sqrt{V_2}}\2\cal{S}^*V_1
}
\,=\, 
\vtwo{\pm K_2^2}{K_1^2}
-\eta \vtwo{\pm\2\cal{C}_{K_2}^{-1}V_1}
{\cal{C}_{K_1}^{-1}V_2}\,,
}
as $\pm \cal{V} V_\pm$.
In \eqref{computing cal F cal VVpm} we used the definition of $\cal{V}$ and $\cal{F}$ for the first equality and the identities \eqref{eta + SV in terms of P} for the second equality. 

For \eqref{formula for spectral radius of cal F} we choose the $+$ in \eqref{Vpm approximate eigenvectors of F}, take the scalar product with $F_+$ and use that $\cal{F}$ is self-adjoint to obtain
\[
\norm{\cal{F}}_{\rm{hs}}\scalar{F_+}{\cal{V}V_+}\,=\, \scalar{F_+}{\cal{V}V_+} - \frac{\eta}{2} \pb{\scalar{F_1}{\cal{C}_{K_2}^{-1}V_1}+ \scalar{F_2}{\cal{C}_{K_1}^{-1}V_2}}\,.
\]

To establish \eqref{approximate eigenvectors for cal F} we apply Lemma~\ref{lmm:Resolvent control for cal S} for $\cal{S}$ replaced by $\wh{\cal{F}}\wh{\cal{F}}^*/\norm{\cal{F}}_{\rm{hs}}^2$ and $\wh{\cal{F}}^*\wh{\cal{F}}/\norm{\cal{F}}_{\rm{hs}}^2$, i.e.\  for the diagonal entries of $(\cal{F}/\norm{\cal{F}}_{\rm{hs}})^2$. Due to \eqref{Vpm approximate eigenvectors of F} the projections of $\cal{V} V_+$ to the first and second component provide approximate eigenvectors for these two operators. The resolvent control from Lemma~\ref{lmm:Resolvent control for cal S} allows us to use analytic perturbation theory and the size of the error term in \eqref{approximate eigenvectors for cal F} is a consequence of \eqref{comparison for 1+tau PP}, \eqref{bounds on V} and the definition of $\cal{V}$ in \eqref{definition of cal V}. 
This finishes the proof of the lemma.
\end{Proof}

\begin{lemma}[Spectral properties of $\cal{T}$]
 The operator $\cal{T}$ defined in \eqref{definition of cal T} satisfies the following properties uniformly for $\eta>0$ and $\tau\ge 0$:
\begin{enumerate}
\item It is self-adjoint, $\cal{T}^*= \cal{T}$.
\item Let $P = \sum_{i=1}^n \pi_i \2p_iq_i^*$ with $\pi_i \ge 0$ and orthonormal bases $(p_i)_i$ and $(q_i)_i$ of $\C^n$ be the singular value decomposition of $P$. The eigenvectors of $\cal{T}$ are 
\bels{cal T eigenvectors}{
\cal{T}\sbb{\vtwo{\msp{10}p_ip_j^*}{\pm q_i q_j^*}}\,=\, \frac{-1\pm\tau \pi_i\pi_j}{\sqrt{(1+\tau \pi_i^2)(1+\tau \pi_j^2)}}\vtwo{\msp{10}p_ip_j^*}{\pm q_i q_j^*}\,.
}
In particular, the spectrum of $\cal{T}$ is bounded away from $1$ by some $\eps_1> 0$ satisfying
\bels{spectrum of cal T}{
\spec(\cal{T}) \,\subseteq\, [-1,1-\eps_1]\,,\qquad 
\eps_1 \,\sim\, 
(1+\tau+\eta^2)\rho^2\,. 
}
\item An eigenvector of $\cal{T}$ corresponding to the eigenvalue $-1$ is given by
\bels{-1 eigenvector of cal T}{
\cal{T}\cal{V} V_{-}\,=\,-\cal{V} V_{-}\,.
}
\item On $\cal{V} V_{+}$ the operator $\cal{T}$ acts contracting, i.e.\  there is an $\eps_2>0$ such that
\bels{cal T on cal V V+}{
\norm{\cal{T}\cal{V} V_{+}}_{\rm{hs}}\,\le\, (1-\eps_2)\norm{\cal{V} V_{+}}_{\rm{hs}}\,,\qquad  \eps_2 \,\sim\, \frac{\tau}{1+\tau +\eta^2}\,.
}
\end{enumerate}
\end{lemma}

\begin{Proof} The self-adjointness of $\cal{T}$ follows immediately from its definition in \eqref{definition of cal T}. The form of the eigenvectors in \eqref{cal T eigenvectors} is a consequence of the following general fact. Let $A=\sum_{i=1}^n\alpha_i a_i x_i^*$ and $B=\sum_{i=1}^n\beta_i b_i y_i^*$ be singular value decompositions of matrices $A$ and $B$ and $\cal{A}R := ARB$ the operator that multiplies a matrix $R$ from the left by $A$ and from the right by $B$. Then $\cal{A}[x_lb_k^*] = \alpha_l \beta_k\2a_l y_k^*$.  In particular,
\begin{align*}
\cal{C}_{f(PP^*)}[p_ip_j^*]\,&=\, f(\pi_i^2)f(\pi_j^2)\,p_ip_j^*\,,\qquad \cal{C}_{f(P^*P)}[q_iq_j^*]\,=\, f(\pi_i^2)f(\pi_j^2)\,q_iq_j^*\,,
\\
\cal{K}_{ P}[p_ip_j^*]\,&=\, \pi_i\2\pi_j\,q_iq_j^*\,,\qquad \msp{75}\cal{K}_{ P^*}[q_iq_j^*]\,=\, \pi_i\2\pi_j\, p_ip_j^*\,.
\end{align*}
 for any function $f$ that is continuous on the positive reals. 
With these formulas \eqref{cal T eigenvectors} is easily verified using the definition of $\cal{T}$. The bound \eqref{spectrum of cal T} on the spectrum of $\cal{T}$ now follows from \eqref{cal T eigenvectors} and \eqref{comparison for 1+tau PP}. 

For \eqref{-1 eigenvector of cal T} and \eqref{cal T on cal V V+} we use the identities 
\[
\cal{V} V_\pm \,=\, \vtwo{ (1+\tau PP^*)^{-1/2}}{\pm (1+\tau P^*\!P)^{-1/2}}\,, \qquad \cal{T}\cal{V} V_{\pm}\,=\,
\vtwo{
\frac{-1\pm \tau PP^*}{(1+\tau PP^*)^{3/2}}
}
{
\frac{\mp 1 +\tau P^*\!P}{(1+\tau P^*\!P)^{3/2}}
}\,,
\]
that follow from  the definitions of $\cal{T}$ and $\cal{V}$ in \eqref{definition of cal T} and \eqref{definition of cal V}, respectively. 
To show \eqref{cal T on cal V V+} we also use \eqref{comparison for 1+tau PP} and \eqref{bounds on V}.
\end{Proof}

\begin{Proof}[Proof of Lemma \ref{lmm:Resolvent control for cal L}]
We start by verifying \eqref{left and right eigenvector for L}. Indeed,  owing to the representation  of $\bs{M}$ in \eqref{definition of M} we have 
\bels{right eigenvector of L proof}{
\scr{P} \2\scr{C}_{\bs{M}}\scr{S}\scr{P}^*V_-\,=\, 
\vtwo{V_1\cal{S}[V_2]V_1+\abs{\zeta}^2U\cal{S}^*[V_1]U^*}{-V_2\cal{S}^*[V_1]V_2-\abs{\zeta}^2U^*\cal{S}[V_2]U}
\,=\, V_- +\eta 
\vtwo{-V_1^2-\abs{\zeta}^2{U}{U}^*}{V_2^2+\abs{\zeta}^2{U}^*{U}}\,,
}
where we used the identities \eqref{im of MDE with U} for the second equality. By using the comparison relations \eqref{bounds on V} and \eqref{comparison relation for U} to bound the last summand on the right hand side of \eqref{right eigenvector of L proof} we conclude \eqref{right eigenvector for L}. The identity \eqref{left eigenvector for L} is verified by using the definition of $\bs{M}$ and \eqref{im of MDE with U} again.

Now we turn to the proof of the resolvent bound  \eqref{resolvent bound for L} for the reduced stability operator $\cal{L} $. We rewrite this operator  using \eqref{stability operator in terms of cal F}
and apply this representation to the resolvent of $\scr{P} \2\scr{C}_{\bs{M}}\scr{S}\scr{P}^*$ to get
\bels{L resolvent in terms of TF}{
\frac{1}{\scr{P} \2\scr{C}_{\bs{M}}\scr{S}\scr{P}^*-\xi}\,=\,\cal{V}^{-1}\frac{1}{\cal{T} \cal{F}-\xi}\cal{V}\,.
}
For $\eta \le \delta^3$ and $\tau  = \abs{\zeta}^2  \le \delta$ we use $\cal{T}= -1 + \ord_{\rm{hs}}(\tau)= -1 + \ord_{\rm{hs}}(\delta)$ which follows from the definition of $\cal{T}$ in \eqref{definition of cal T} and \eqref{comparison for 1+tau PP} as well as $V_1\sim V_2\sim 1$ in this regime (cf.\ \eqref{bounds on V}). From \eqref{bound on cal V}, the spectral properties  of $\cal{F}$,  \eqref{spectral gap of cal F}, and $\norm{\cal{F}}_{\rm{hs}} =1 +\ord(\eta)= 1 +\ord(\delta^3) $ (cf.\ \eqref{formula for spectral radius of cal F}), as well as \eqref{L resolvent in terms of TF} we infer that there is an $\eps \in (0,1/2)$ such that
\[
\sup
\cB{ \norm{(\scr{P} \2\scr{C}_{\bs{M}}\scr{S}\scr{P}^*-\xi)^{-1}}_{\rm{hs}}\,:\; \xi \in \C\,,\; \xi \not \in 
 (-1+\D_\eps) \cup \D_{1- 2\eps}\cup (1+\D_\eps)
}\,\lesssim_\delta\,1 \,, \qquad \eps \,\sim_\delta\, 1\,.
\]
In particular, \eqref{resolvent bound for L} holds true. 
The non-degeneracy \eqref{nondegeneracy for L} of the eigenvalue in $\D_\eps$
follows from the non-degeneracy of the eigenvalue $\norm{\cal{F}}_{\rm{hs}}$ of $\cal{F}$ as stated in Lemma~\ref{lmm:Properties of cal F}. The statement $\lambda=\ord_\delta(\eta)$ 
 about the non-degenerate isolated eigenvalue in \eqref{nondegeneracy for L} follows from $V_-$ being an approximate eigenvector (cf.\ \eqref{right eigenvector for L}) and the resolvent bound \eqref{resolvent bound for L}.

For $\eta \le \delta^3$ and $\tau \in [\delta,1-\delta]$ we will apply Lemma~\ref{lmm:TF Stability lemma} with the choices $F:=\cal{F}/\norm{\cal{F}}_{\rm{hs}}$ and $T:=\cal{T}$. We verify the assumptions of the lemma. The required upper bound $\norm{\cal{T}}_{\rm{hs}}\le 1$ follows from \eqref{spectrum of cal T} and \eqref{spectrum of F}  holds true because of \eqref{spectral gap of cal F}.
Furthermore according to \eqref{approximate eigenvectors for cal F} and \eqref{cal T on cal V V+} we have
\bels{upper bound on cal T F+}{
\norm{\cal{T} F_+}_{\rm{hs}}\,\le\,  1-\eps_2+\ord\pb{\eta/\rho}
\,=\, 1-\eps_2+\ord(\delta^2)
\,,\qquad \eps_2 \,\sim\, \tau \,\gtrsim\,  \delta\,,
}
 where $F_+$ is the normalized eigenvector of $\cal{F}$ corresponding to the eigenvalue $\norm{\cal{F}}_{\rm{hs}}$ and we used \eqref{bounds on V} to see the bounds in terms of $\delta$. 
We also have 
\bels{bound on 1+cal T F-}{
\norm{(1+\cal{T}) F_-}_{\rm{hs}} \,\lesssim\,  \frac{\eta}{\rho}\,\lesssim\, \delta^{5/2}\,,
}
by \eqref{approximate eigenvectors for cal F} and \eqref{-1 eigenvector of cal T}, where $F_-$ is the normalised eigenvector of $\cal{F}$ corresponding to the eigenvalue $-\norm{\cal{F}}_{\rm{hs}}$. 
Thus Lemma~\ref{lmm:TF Stability lemma} is applicable because of \eqref{upper bound on cal T F+} and \eqref{bound on 1+cal T F-} as long as $\delta\sim 1$ is chosen sufficiently small. Thus we find
\bels{TF resolvent bound}{
\sup
\cb{ \norm{(\cal{T}\cal{F}-\norm{\cal{F}}_{\rm{hs}}\zeta)^{-1}}_{\rm{hs}}\,:\; \zeta \in \C\,,\; \zeta \not \in 
 \D_{1-2\eps}\cup(1+\D_{\eps})
}\,\lesssim_\delta\, 1\,,
}
for some $\eps \sim \delta^{5/2}$. 
Since $\norm{\cal{F}}_{\rm{hs}} =1 +\ord(\eta/\rho)= 1 +\ord(\delta^{2}) $ (cf.\ \eqref{formula for spectral radius of cal F} and \eqref{bounds on V}) we infer \eqref{resolvent bound for L} from \eqref{TF resolvent bound} by using \eqref{L resolvent in terms of TF} and \eqref{bound on cal V}. The non-degeneracy of the isolated eigenvalue $\lambda$ in \eqref{nondegeneracy for L} stems from \eqref{Nondegeneracy of TF eigenvalue} and the resolvent bound \eqref{resolvent bound for L} in combination with the approximate eigenvector equation \eqref{right eigenvector for L} for $V_-$ implies $\lambda= \ord_\delta(\eta)$. 

Finally, note that $\lambda \ne 0$ because the representation \eqref{stability operator in terms of cal F} shows that with $\norm{\cal{T}}\le 1$ (cf.\ \eqref{spectrum of cal T}) the operator $\cal{L}$ is invertible as long as $\norm{\cal{F}}_{\rm{hs}}<1$, which is always true for $\eta>0$ due to the right hand side of \eqref{formula for spectral radius of cal F} not vanishing. 
\end{Proof}

\begin{Proof}[Proof of Proposition~\ref{prp:Existence and uniqueness}] For $\bs{D}=\bs{0}$ the equation \eqref{perturbed MDE bulk regime} is equivalent to \eqref{eq:MDE} and thus by Theorem~\ref{thr:Stability} for any $\abs{\zeta}^2 =\tau \in [0,\spradius)$ we can extend the solution $\bs{M}(\zeta,\eta)$ analytically to $\eta =0$. Thus also the solution $V_1,V_2$ of the Dyson equation \eqref{Dyson equation} can be analytically extended to $\eta =0$. This proves 
 the  existence of a positive definite solution to \eqref{Dyson equation at eta = 0}.

For the uniqueness, note that in the proof of  Theorem~\ref{thr:Stability} and in particular for the key input, Lemma~\ref{lmm:Resolvent control for scr L}, we never used $\eta>0$, but only that $\bs{M} \in \bs{E}_-^\perp$ solves \eqref{eq:MDE} and has positive definite imaginary part with lower and upper bounds depending on model parameters and $\delta$. Thus for any positive definite solution $V_1,V_2$ that satisfies \eqref{Dyson equation at eta = 0} and \eqref{constraint at eta = 0} we can construct a solution $\bs{M}$ of \eqref{eq:MDE} at $\eta =0$ through \eqref{definition of M} and Theorem~\ref{thr:Stability} also applies to this $\bs{M}=\bs{M}(\sqrt{\tau},0)$ with $c_1,c_2$ from \eqref{norm neighbourhoods of stability} now depending also on the lower and upper bounds on $V_1,V_2$. By analyticity of $\bs{G}$ in all variables $\bs{G}(\zeta,\ol{\zeta},\eta,\bs{0})$ has positive definite imaginary part for sufficiently small $\abs{\zeta-\sqrt{\tau}}$ and $\eta >0$. We conclude $\bs{G}(\zeta,\ol{\zeta},\eta,\bs{0}) = \bs{M}(\zeta,\eta)$ since it solves \eqref{perturbed MDE bulk regime} with $\bs{D}=\bs{0}$ and $\bs{M} = \lim_{\eta \downarrow 0}\bs{M}(\sqrt{\tau},\eta)$, establishing uniqueness of the solution to \eqref{Dyson equation at eta = 0}. 
\end{Proof}

As used in the proof of Proposition~\ref{prp:Existence and uniqueness} above, the uniformity of the statement of Theorem~\ref{thr:Stability} in $\eta>0$ allows for an extension of $\bs{M}$ as well as $V_1,V_2$ to $\eta =0$ in the following sense. 

 \begin{corollary}[Extension inside the spectrum] 
 \label{crl:Extension inside the spectrum}
 The solution $\bs{M}$ of the MDE \eqref{eq:MDE} has a unique continuous extension to $\C\times [0,\infty)$, i.e.\  to $\eta =0$. For every $\zeta \not \in \partial \D$ this extension, still denoted by $\bs{M}$, also has a continuation to a neighbourhood of $(\zeta,0)$ that is real analytic in $\re \zeta, \im \zeta, \eta$. The size of this neighbourhood only depends on the model parameters and on $\dist(\zeta,\partial\D)=\abs{\abs{\zeta} -1}$. 
 
 Similarly $V_1,V_2$ admit a continuous extension to $[0,\infty) \times [0,\infty)$ that extends to an analytic function in a neighbourhood of $(\tau,0)$ for any $\tau \in [0,\infty)\setminus\{1\}$  with the size of the neighbourhood depending only on $\abs{\tau-1}$ in addition to the model parameters. 
  \end{corollary}

\section{Self-consistent density of states}
\label{Sec:Self-consistent density of states}
In this section we use the information about the solution of the Dyson equation to control the self-consistent density of states  $\sigma$  corresponding to $X$. In Subsection~\ref{subsec:Upper and lower bounds in the bulk} we begin with establishing upper and lower bounds on the density.
These bounds rely on a novel representation of $\sigma$ in \eqref{Formula for sigma}. 
 In Subsection~\ref{Subsec:Solution close to the edge} we provide a detailed description of $V_i$ and $\sigma$ at the edge of the spectrum. We end the  subsection by summarising its results in the proof of Theorem~\ref{thr:Density}. 
 Subsections~\ref{subsec:proof_sigma_dist_derivative} and \ref{sec:proof_Brown_measure} contain the proofs of Lemma~\ref{lem:sigma_as_dist_derivative} and Proposition~\ref{pro:Brown_measure}, respectively.

\subsection{Upper and lower bounds in the bulk} \label{subsec:Upper and lower bounds in the bulk}
In this subsection we establish lower and upper bounds on the density $\sigma$ inside the spectrum, i.e.\ we show \eqref{lower and upper bounds on sigma}  away from the edge of the spectrum at $\abs{\zeta} =1$.

\begin{lemma}[Formula for density] For any $\zeta \in \D$ the density $\sigma$ admits the formula
\bels{Formula for sigma}{
\sigma(\zeta)\,=\, \frac{2}{\pi}\scalarbb{\cal{L}^{-1}\vtwo{V_1\frac{1}{\cal{S}^*V_1}V_1}{V_2\frac{1}{\cal{S}V_2}V_2}}{\mtwo{\cal{S}V_2}{\cal{S}^*V_1}}
\,=\, 
\frac{1}{\pi \tau}\scalar{Y}{(1-\cal{T} \cal{F}^2\cal{T})Y}\,,
}
where all expressions on the right hand side are evaluated at $\eta =0$ (cf.\ Corollary~\ref{crl:Extension inside the spectrum}) and $\tau=\abs{\zeta}^2$ and where
\bels{definitions of Ki}{
Y := (1- \cal{F}\cal{T})^{-1}|_{\cal{V}[V_-]^\perp} \mtwo{K_2^2}{K_1^2}.
}
and $K_i$ the matrices from \eqref{definition of P}. 
For $\tau=0$ the very right hand side of \eqref{Formula for sigma} is interpreted as its limit $\tau \downarrow 0$.
Here and in the following, the notation $(1- \cal{F}\cal{T})^{-1}|_{\cal{V}[V_-]^\perp}$ on the right-hand side of \eqref{definitions of Ki} is 
understood as $((1-\cal F \cal T)|_{\cal V[V_-]^\perp})^{-1}$. 
\end{lemma}
\begin{Proof}
By definition of $\sigma$ in \eqref{definition of sigma} and the identity \eqref{V equation with U} we have
\bels{first formula for sigma}{
\sigma(\zeta)= \frac{1}{\pi} \partial_\tau(\tau\2\avg{U(\tau,0)})|_{\tau = \abs{\zeta}^2} = -\frac{1}{2\pi}\partial_\tau(\scalar{V_1}{\cal{S}V_2}+\scalar{V_2}{\cal{S}^*V_1}) |_{\tau = \abs{\zeta}^2, \eta =0}\,,
}
for any $\zeta \in \D$. By rotational symmetry it suffices to establish \eqref{Formula for sigma} at $\zeta =\sqrt{\tau} >0$. Thus we denote $\sigma=\sigma(\sqrt{\tau})$ and $\bs{M} = \bs{M}(\sqrt{\tau},0)$. By \eqref{first formula for sigma}, 
 the representation of $\bs{M}$ from \eqref{definition of M} and 
the definition of  $\scr{S}$ in \eqref{definition scr S and Z} 
we find
\bels{second formula for sigma}{
\sigma= -\frac{1}{\pi}\partial_\tau \scalar{\bs{M}}{\scr{S}\bs{M}} = -\frac{2}{\pi}\scalar{\scr{P}\partial_\tau \bs{M}}{\scr{P}\scr{S}\bs{M}}\,,
}
where we used the structure of $\scr{S}$ and the projection $\scr{P}$ from \eqref{definition of scr P} in the second equality.
We compute the derivative of $\bs{M}$ with respect to $\tau$ by differentiating both sides of \eqref{eq:MDE} and solving for
\bels{tau derivative of M}{
\partial_\tau \bs{M} = (\rm{Id}-\scr{C}_{\bs{M}}\scr{S})^{-1}\scr{C}_{\bs{M}}\partial_\tau \bs{Z}(\sqrt{\tau},\sqrt{\tau})\,.
}
By definition of $\bs{M}$ and the identities \eqref{U identity} for $U$ we have
\bels{compute P of CM}{
\scr{P}\scr{C}_{\bs{M}}\partial_\tau \bs{Z}(\sqrt{\tau},\sqrt{\tau})= -\frac{\ii}{2} \vtwo{UV_1+V_1 U^*}{V_2 U + U^*V_2}
= -\ii\vtwo{V_1\frac{1}{\cal{S}^*V_1}V_1}{V_2\frac{1}{\cal{S}V_2}V_2}\,.
}
Thus, inserting \eqref{compute P of CM} into \eqref{tau derivative of M} and recalling the definition of $\cal{L}$ from \eqref{definition cal L}, shows
\bels{scr P partial tau M}{
\scr{P}\partial_\tau \bs{M} = -\ii\2\cal{L}^{-1}\vtwo{V_1\frac{1}{\cal{S}^*V_1}V_1}{V_2\frac{1}{\cal{S}V_2}V_2}\,.
}
We plug this into \eqref{second formula for sigma} and verify the first equality in \eqref{Formula for sigma}. 
 Note that $\cal{L}^{-1}$ is applied to the orthogonal complement of $(\cal{S}V_2,\cal{S}^*V_1)$ in \eqref{scr P partial tau M}. To check the orthogonality of the vector on the right hand side we can use \eqref{U identity} at $\eta=0$.

For the second equality in \eqref{Formula for sigma} we recall the definitions of $\cal{T}$, $\cal{F}$ and $\cal{V}$ from \eqref{definition of cal T}, \eqref{definition of cal F} and \eqref{definition of cal V}, as well as the identities \eqref{eta + SV in terms of P} that take the form
\[
1+\tau\1PP^* = K_2^{-4}\,, \qquad 1+\tau\1P^*P = K_1^{-4} \qquad \text{with}\qquad P \,=\,K_2^{-4}\sqrt{V_1}\sqrt{V_2}= \sqrt{V_1}\sqrt{V_2}\,K_1^{-4}\,,
\]
at $\eta =0$.  Then we  compute
\bels{cal V applied to}{
(\cal{V}^*)^{-1}\mtwo{\cal{S}V_2}{\cal{S}^*V_1}=
\mtwo{K_2^2}{K_1^2}
\,, \qquad
\cal{V}\vtwo{V_1\frac{1}{\cal{S}^*V_1}V_1}{V_2\frac{1}{\cal{S}V_2}V_2}
=\vtwo{\cal{C}_{\!K_2}\cal{K}_{P^*}K_1^4}{\cal{C}_{\!K_1}\cal{K}_{P}K_2^4}
=
\frac{1}{\tau}\vtwo{\cal{C}_{\!K_2}[1-K_2^4]}{\cal{C}_{\!K_1}[1-K_1^4]}\,,
}
where we used \eqref{Dyson equation at eta = 0} for the last equality. Again with \eqref{Dyson equation at eta = 0} we also have
\bels{1+cal T applied to}{
(1+\cal{T})\mtwo{K_2^2}{K_1^2}=2\vtwo{\cal{C}_{\!K_2}[1-K_2^4]}{\cal{C}_{\!K_1}[1-K_1^4]}\,.
}
Now we insert the representation \eqref{stability operator in terms of cal F} for the reduced stability operator into the middle formula of \eqref{Formula for sigma}. Afterwards we use \eqref{cal V applied to} and \eqref{1+cal T applied to} to get
\[
\sigma
=
\frac{2}{\pi \tau}\scalarbb{ (1-\cal{T} \cal{F})^{-1}\vtwo{\cal{C}_{\!K_2}[1-K_2^4]}{\cal{C}_{\!K_1}[1-K_1^4]}}{\mtwo{K_2^2}{K_1^2}}
=
\frac{1}{\pi \tau}\scalarbb{ (1-\cal{T} \cal{F})^{-1}(1+\cal{T})\mtwo{K_2^2}{K_1^2}}{\mtwo{K_2^2}{K_1^2}}\,,
\]
where the inverse of $1-\cal{T} \cal{F}$ is restricted to $\cal{V}[V_-]^\perp$. 
The vector in the second argument of the scalar product is a representation of the Perron-Frobenius eigenvector for $\cal{F}$. Indeed, by the definitions of $\cal{F}$ in \eqref{definition of cal F} and $K_i$ in \eqref{definitions of Ki} we see that
\bels{K eigenvector of cal F}{
\cal{F}\mtwo{K_2^2}{K_1^2} = \mtwo{K_2^2}{K_1^2}\,.
}
Because of \eqref{K eigenvector of cal F} we also have the identity
\[
(1-\cal{T} \cal{F})^{-1}(1+\cal{T})\mtwo{K_2^2}{K_1^2}= (1-\cal{T} \cal{F})^{-1}(1-\cal{T} \cal{F}^2\cal{T})(1-\cal{F}\cal{T} )^{-1}\mtwo{K_2^2}{K_1^2}\,,
\]
which finishes the proof of the second equality in \eqref{Formula for sigma} and, thus, the proof of the lemma. 
\end{Proof}

\begin{corollary}[Bounds on the density] \label{crl:Bounds on the density} 
For any $\delta\in (0,1)$, we have $\sigma(\zeta) \sim_\delta 1$ uniformly for $\zeta \in \D_{1-\delta}$. 
\end{corollary}
\begin{Proof}
We consider two separate regimes. First upper and lower bounds on $\sigma$ follow in a neighbourhood of $\zeta =0$ by continuity (cf.\ Corollary~\ref{crl:Extension inside the spectrum} and \eqref{definition of sigma}) of $\sigma$ and $\sigma(0) \gtrsim 1$. 
The latter is easy to see because at $\tau=\eta =0$ the Dyson equation simplifies to
\bels{}{
1=V_1 \cal{S}V_2\,, \qquad 1=V_2 \cal{S}^*V_1\,,
}
and we have
\[
U= V_1V_2\,, \qquad \cal{L} = \mtwo{1 & \cal{C}_{V_1}\cal{S}}{ \cal{C}_{V_2}\cal{S}^*&1}\,.
\]
In particular, the reduced stability has the form $\cal{L}  = 1-\cal{A}$, where $\cal{A}$ preserves the cone of positive definite matrix pairs. Thus the first identity in \eqref{Formula for sigma} implies
\[
\sigma(0)\,=\, \frac{2}{\pi}\scalarbb{\cal{L}^{-1}\vtwo{V_1\frac{1}{V_2}V_1}{V_2\frac{1}{V_1}V_2}}{\mtwo{\frac{1}{V_1}}{\frac{1}{V_2}}}\,\ge\, \frac{1}{\pi}\pb{\avgb{V_1V_2^{-1}}+\avgb{V_2V_1^{-1}}}\,.
\]
Note that we can expand $\cal{L}^{-1} = (1-\cal{A})^{-1}$ in a Neumann series because of the representation \eqref{stability operator in terms of cal F}, $\norm{\cal{T}}\le1$ and $\norm{\cal{F}}_{\rm{hs}}<1$. 

Now we consider the regime $ 1 \lesssim \tau^{1/2}=\abs{\zeta} \leq 1 - \delta$. Here,  owing to the second relation in 
\eqref{Formula for sigma}, we have the lower and upper bound
\[
\sigma(\zeta)\,\sim\, \scalarbb{K}{\frac{1}{1-\cal{A}}(1-\cal{A}\cal{A}^*)\frac{1}{1-\cal{A}^*}K}\,,
\]
where  $K = (K_2^2,K_1^2) \in (\cal{V}V_-)^\perp$ with $K_i \sim 1$ and $\cal{A}=\cal{T}\cal{F}$. Thus for $\sigma \sim 1$ it suffices to check that $\norm{\cal{A}|_{(\cal{V}V_-)^\perp}}_{\rm{hs}}\le 1-\eps$ for some $\eps\gtrsim 1$. 
We apply Lemma~\ref{lmm:TF Stability lemma}  with $T = \cal{T}$, $F = \cal{F}/\norm{\cal{F}}_{\rm{hs}}$ and $f_\pm = \cal{V} V_\pm /
\norm{\cal{V} V_\pm}_{\rm{hs}}$  and note that the non-degenerate eigenvalue $1$ of $\cal{A}$ corresponds to the eigenvector $\cal{V}V_-$ which is projected out when we take the norm. Thus we have the resolvent bound
\[
\sup_{\omega \not \in  \D_{1-\eps}}\normb{(\cal{A}-\omega)^{-1}|_{(\cal{V}V_-)^\perp}}_{\rm{hs}}  \lesssim 1\,,
\]
for some $\eps\sim 1$, 
which implies the desired norm bound. 
\end{Proof}

\subsection{Solution close to the edge}
\label{Subsec:Solution close to the edge}
In this subsection we explicitly determine the leading order of the solution $V_1,V_2$ to \eqref{Dyson equation} close to the edge $\tau=\abs{\zeta}^2=1$ of the spectrum. We use the result to determine the jump height \eqref{jump height of sigma} of the density $\sigma$ at the edge. 
Let $S_2$ and $S_1$ be the unique positive definite right and left eigenvectors of $\cal{S}$, respectively,
i.e.\ $\cal{S} S_2 = S_2$ and $\cal{S}^* S_1 = S_1$, satisfying $\avg{S_1} = \avg{S_2} = 1$.  We also write $\rhoDOS:=\rhoDOS_\zeta:=\avg{V_1}/\pi$ for the harmonic extension of the self-consistent density of states of $\bs{H}_\zeta$ to the complex upper half plane and recall that $\rhoDOS$ is comparable to the right hand side of \eqref{bounds on V}. 

\begin{proposition}[Solution at the edge] For any $\tau,\eta \in [0,2]$ we have the expansion
\bels{V expansion at edge}{
V_1 = \alpha \2S_1 + \ord(\eta +\rhoDOS^3)\,, \qquad V_2 =\alpha  \2S_2 + \ord(\eta +\rhoDOS^3)\,, \qquad \alpha:=\frac{\scalar{S_1}{V_2}}{\scalar{S_1}{S_2}}\,,
}
where $\alpha$ satisfies the cubic equation
\bels{alpha cubic equation}{
 \alpha^3\avg{(S_1S_2)^2}+\alpha\2(\tau-1) \1 \avg{S_1S_2}-\eta   \,=\,\ord(\rhoDOS^5+\eta\rhoDOS^2 )\,.
} 
\end{proposition}
\begin{Proof}
 We write $\tau=1+\eps$ for some small $\eps$. The case when $\eps \le  -c$ for some constant $c \sim 1$ is trivial since then $\rhoDOS \sim 1$ and the error term in \eqref{alpha cubic equation}  dominates. Similarly, for  $\eps \ge  c$ we have $\rhoDOS \sim \alpha \sim \eta$, i.e.\ in both regimes the proposition does not contain any information.  Solving \eqref{U equation} shows
 \bels{U expansion at edge}{
U \,=\, \frac{1}{2 (1+\eps)}\pB{1 + \sqrt{1-4 (1+\eps) V_1 V_2}}
\,=\, 
\frac{1}{1+\eps} - V_1 V_2- (1+\eps) (V_1 V_2)^2 +\ord(\rhoDOS^6)\,.
}
We use this expansion for $U$ in \eqref{V equation with U} and find
\bes{
0\,=\, V_1(\eta + \cal{S}V_2) + (1+\eps)\2 U-1
\,=\, V_1(\eta + \cal{S}V_2) -(1+\eps) V_1 V_2 - (1+\eps)^2(V_1 V_2)^2+\ord(\rhoDOS^6)\,.
}
Multiplying with $V_1^{-1}$ from the left and using the decomposition $V_i \,=\, \alpha_i S_i + \wt{V}_i$ shows 
\bels{wt V expansion}{
(1+\eps-\cal{S})\wt{V}_2\,=\, \eta -\eps\1 \alpha_2 S_2-  (1+\eps)^2V_2V_1V_2+\ord(\rhoDOS^5)\,.
}
Here  $\wt{V}_1$ and $\wt{V}_2$ are the spectral projection of $V_1$ and $V_2$ corresponding to the spectrum of $\cal{S}$ and $\cal{S}^*$ complementary to the isolated eigenvalue $1$, respectively, i.e.\    $\wt{V}_i = \cal{Q}_i V_i$, with 
\[
\cal{P}_1\,=\, \frac{\scalar{S_2}{\2\cdot\2}}{\avg{S_1S_2}} S_1\,,\qquad \cal{Q}_1\,=\, 1-\cal{P}_1\,, \qquad \cal{P}_2\,=\, \frac{\scalar{S_1}{\2\cdot\2}}{\avg{S_1S_2}} S_2\,,\qquad \cal{Q}_2\,=\, 1-\cal{P}_2\,.
\]

In particular, projecting both sides of \eqref{wt V expansion} onto the range of $\cal{Q}_2$ implies $\norm{\wt{V}_2}\,\lesssim\, \eta +\rhoDOS^3$. 
Here we used that $\norm{(1+\eps-\cal{S})^{-1}\cal{Q}_2} \lesssim 1$, which follows from Lemma~\ref{lmm:Resolvent control for cal S}. 
By exchanging the roles of $V_1$ and $V_2$ we also find $\norm{\wt{V}_1}\,\lesssim\, \eta +\rhoDOS^3$. Therefore, \eqref{wt V expansion} can be expanded further as 
\[
(1+\eps-\cal{S})\wt{V}_2\,=\, \eta -\eps \1\alpha_2 S_2-  \alpha_1\alpha_2^2S_2S_1S_2+\ord(\abs{\eps}\1\rhoDOS^3+\rhoDOS^5+\eta\rhoDOS^2 )\,.
\]
Now we apply the rank one projection $\cal{P}_2$ on both sides and get 
\[
0\,=\, \eta -\eps \1\alpha_2 \avg{S_1S_2}-  \alpha_1\alpha_2^2\avg{(S_1S_2)^2}+\ord(\rhoDOS^5+\eta\rhoDOS^2 )\,,
\]
where we used  $\avg{S_1} = 1$ and, for the error term,  $\abs{\eps}\1\rhoDOS^3\lesssim \rhoDOS^5+\eta\rhoDOS^2$  due to \eqref{bounds on V}.  
Finally \eqref{alpha cubic equation} follows from 
\begin{equation} \label{eq:alpha_1_approx_alpha_2}
\alpha_1\,=\, \alpha_2 +\ord(\eta +\rhoDOS^3)\,,
\end{equation}
 which is a consequence of \eqref{avg V1 = avg V2} 
and $\alpha_2 = \alpha$. 
 Moreover, \eqref{eq:alpha_1_approx_alpha_2} and $\norm{\wt{V}_i} \lesssim \eta  + \rhoDOS^3$ yield \eqref{V expansion at edge}.   
\end{Proof}

For the next corollary, we 
introduce $\cal{M}:\C^{n\times n}\times \C^{n\times n} \to \C^{n\times n}\times \C^{n\times n}$ defined by
\bels{definition cal M}{
\cal{M}:= \frac{1}{\avg{V_+}}\mtwo{\eta +\cal{S}V_2 & 0 }{0 & \eta + \cal{S}^*V_1}\,.
}
Moreover, we recall that $\cal{L}$ and $E_-$ were defined in \eqref{definition cal L} and \eqref{definition E pm}, respectively. 

\begin{corollary}[Resolvent control for stability operator close to the edge] 
\label{crl:Resolvent control for stability operator close to the edge}
The following holds. 
\begin{enumerate} 
\item The operator $\cal{M} \cal{L}$ has the invariant subspace $E_-^\perp \subset \C^{n\times n} \times \C^{n\times n}$, i.e.\ $\cal{M}\cal{L}E_-^\perp \subset E_-^\perp$. 
\item There is $\delta \sim 1$ such that, for any $\eta \in (0,\delta)$ and $\zeta \in \D_{1+\delta} \setminus \D_{1-\delta}$, 
the eigenvalues of $\cal{M} \cal{L}|_{E_-^\perp}$ close to zero are isolated in the sense that 
\bels{resolvent bound for scr JL}{
\sup
\cB{ \norm{(\cal{M}\cal{L}-\xi)^{-1}|_{E_-^\perp}}_{\#}\,:\; \xi  \in 
\D_{2\eps}\setminus \D_\eps
}\,\lesssim\,1 \,,
}
for some $\eps \sim 1$ and $\#=\rm{hs},\norm{\1\cdot\1}$. In fact, $\cal{M}\cal{L}|_{E_-^\perp}$ has only one eigenvalue in $\D_\eps$. This eigenvalue is simple and the spectral projection
\[
\cal{P}:=-\frac{1}{2\pi \ii}\oint_{\partial \D_\eps} \dd \xi \, (\cal{M}\cal{L}-\xi)^{-1}|_{{E}_-^\perp}
\]
has rank one. 
\end{enumerate} 
\end{corollary}
\begin{Proof}
 The invariance of $E_-^\perp$ under $\cal{M}\cal{L}$ is a direct consequence of $\cal{L}^*\cal{M}^* E_- = \eta \avg{V_+}^{-1} E_- $.  
The operator $\cal{M}\cal{L}|_{{E}_-^\perp}$ is a small perturbation of $\cal{K}$ from Corollary~\ref{crl:Resolvent control for scr J0 L0}, 
 since 
\[ \cal{M} \cal{L} = \cal{K} + \ord(\eta/\rhoDOS +\rhoDOS + \abs{\abs{\zeta} - 1} ) \] 
due to \eqref{V expansion at edge} and \eqref{U expansion at edge}. 
Thus the claim follows from  Corollary~\ref{crl:Resolvent control for scr J0 L0} 
 for $\#=\rm{hs},\norm{\1\cdot\1}$ by perturbation theory  for sufficiently small $\delta \sim 1$.  
\end{Proof}

\begin{corollary}[Density at the edge]
\label{crl:Density at the edge} 
At the edge of the spectrum the  self-consistent density of states has an expansion 
\bels{sigma expansion at the edge}{
\sigma(\zeta) = \frac{\avg{S_1S_2}^2}{\pi\avg{(S_1S_2)^2}}+\ord(1-\abs{\zeta})\,,
}
for any $\zeta \in \D$.
\end{corollary}

\begin{Proof}
We set $\eta=0$ throughout the proof.  
We use the first identity in \eqref{Formula for sigma} and insert $\cal{M}$ from \eqref{definition cal M} to find
\bels{sigma formula with ML}{
\sigma= \frac{2}{\pi\avg{V_+}}\scalarbb{(\cal{M}\cal{L})^{-1}\vtwo{(\cal{S}V_2)V_1\frac{1}{\cal{S}^*V_1}V_1}{(\cal{S}^*V_1)V_2\frac{1}{\cal{S}V_2}V_2}}{\mtwo{\cal{S}V_2}{\cal{S}^*V_1}}\,.
}
We consider $\tau = \abs{\zeta}^2 = 1-\eps$ for some $\eps>0$. Since $\eta = 0$, we get from \eqref{V expansion at edge} that 
\bels{Vi edge expansion}{
V_i = \sqrt{\kappa\2\eps} \2S_i+\ord(\eps^{3/2})\,,\qquad \kappa :=\frac{\avg{S_1S_2}}{\avg{(S_1S_2)^2}}\,,
}
where we used $\rhoDOS \sim \sqrt{\eps}$. From the expansion of $U$ in \eqref{U expansion at edge} this implies
\begin{equation} \label{eq:expansion_tau_U}
\tau \2U\,=\,1 
-\eps\2\kappa\2 S_1 S_2+\ord(\eps^2) \,.
\end{equation}
 Plugging \eqref{eq:expansion_tau_U} and \eqref{Vi edge expansion} into the definitions of $\cal{L}$ and $\cal{M}$ in \eqref{definition cal L} and \eqref{definition cal M}, respectively, yields 
\[
\cal{L}= \cal{L}_\rm{e} + \eps\2\cal{D} +\ord(\eps^{2}), \quad \cal{L}_\rm{e} :=
\mtwo{1-\cal{S}^* & 0}
{ 0 &1- \cal{S}}, 
 \qquad \cal{M} = \cal{M}_{\rm{e}} + \ord(\eps),\quad \cal{M}_\rm{e} := \begin{pmatrix} S_2 & 0 \\ 0  & S_1 \end{pmatrix}, 
\]
where the  first order perturbation  of $\cal{L}$  is given by
\[
\cal{D}\vtwo{A_1}{A_2}=
\vtwo{ -\cal{S}^*A_1+\kappa\2S_1S_2 \cal{S}^*  A_1+\kappa\2 (\cal{S}^* A_1)  S_2S_1+ \kappa\2\cal{C}_{S_1}\cal{S}A_2}
{\kappa\2\cal{C}_{S_2}\cal{S}^*A_1 - \cal{S}A_2+   \kappa\2S_2S_1 \cal{S}  A_2+\kappa\2 (\cal{S} A_2)  S_1S_2}\,.
\]

According to Corollary~\ref{crl:Resolvent control for stability operator close to the edge} the operator 
$\cal{M}\cal{L} =(\cal{M}_{\rm{e}} + \ord(\eps) )(\cal{L}_\rm{e} + \eps\2\cal{D} +\ord(\eps^{2}))$ 
has an isolated eigenvalue $\wt \lambda$ close to $0$ when restricted to $E_-^\perp$. Therefore, we can use perturbation theory to determine its value to leading order
\bels{expansion of wt lambda}{
\wt{\lambda} = \eps \frac{\scalar{E_+}{\cal{M}_{\rm{e}}\cal{D}S_+}}{\scalar{E_+}{S_+}} + \ord(\eps^2) = 2\2\eps\2\avg{S_1S_2} +\ord(\eps^2)\,,
}
where we used $\cal{L}_{\rm{e}}S_+=0$ with $S_+ =(S_1,S_2)$ the right eigenvector and $E_+ = (1,1)$ the left eigenvector of the unperturbed operator $\cal{M}_{\rm{e}}\cal{L}_{\rm{e}} = \cal{K}$ 
 (cf.\ Corollary~\ref{crl:Resolvent control for scr J0 L0}).  The spectral projection corresponding to the eigenvalue $0$ of $\cal{M}_{\rm{e}}\cal{L}_{\rm{e}}$ is 
\[
\cal{P}_{\rm{e}} = \frac{\scalar{E_+}{\2\cdot\2}}{\scalar{E_+}{S_+}}S_+\,,
\]
and thus inserting \eqref{expansion of wt lambda} and \eqref{Vi edge expansion} into \eqref{sigma formula with ML} yields
\[
\sigma= \frac{2\eps\2\kappa}{\pi\wt{\lambda}}\scalarbb{\cal{P}_\rm{e}\vtwo{S_2S_1}{S_1S_2}}{\mtwo{S_2}{S_1}}
+\ord(\eps)
=\frac{\avg{S_1S_2}^2}{\pi\avg{(S_1S_2)^2}}+\ord(\eps)
\,.
\]
This finishes the proof.
\end{Proof}

\begin{Proof}[Proof of Theorem~\ref{thr:Density}] The analyticity of $V_1,V_2$ and thus the well definedness of $\sigma$ in \eqref{definition of sigma}  was shown in Corollary~\ref{crl:Extension inside the spectrum}, the upper and lower bounds on $\sigma$ from \eqref{lower and upper bounds on sigma} in Corollary~\ref{crl:Bounds on the density} away from the edge and in Corollary~\ref{crl:Density at the edge} close to the edge. 
Integrating the definition of $\sigma(\zeta)$ over $\D_1$ and recalling $\varrho(\cal{S}) =1$ from 
\eqref{eq:normalization_spectral_radius} as well as $V_1(\tau) \to 0$ and $V_2(\tau) \to 0$ for $\tau \to 1$ due 
Corollary~\ref{crl:Extension inside the spectrum} and \eqref{bounds on V} imply that $\sigma$ is a probability density on $\C$. 
Finally the jump height \eqref{jump height of sigma} of $\sigma$ right at the spectral edge is read off from \eqref{sigma expansion at the edge}.
\end{Proof}

\subsection{Proof of Lemma~\ref{lem:sigma_as_dist_derivative}}  \label{subsec:proof_sigma_dist_derivative} 

In this subsection, we prove Lemma~\ref{lem:sigma_as_dist_derivative}, the basic property of $\sigma$ 
used in the proofs of the global and local inhomogeneous circular law, Theorem~\ref{thm:global_law_for_X} 
and Theorem~\ref{thr:Local inhomogeneous circular law}. 

\begin{Proof}[Proof of Lemma~\ref{lem:sigma_as_dist_derivative}]
Recall the normalization $\varrho(\mathcal{S}) = 1$ from \eqref{eq:normalization_spectral_radius}
and write $U(\tau) = U(\tau,0)$, $V_1(\tau) = V_1(\tau,0)$, and $V_2(\tau) = V_2(\tau,0)$.  

As a first step, we now compute in the integral in the definition of $L$ in \eqref{eq:sigma_integration_part_def_L}. This will yield 
\begin{equation} \label{eq:L_general_representation} 
 L(\zeta) =  \frac{1}{2}\bigg( \avg{V_1 \mathcal{S} V_2} - \frac{1}{2} \avgb{\log \big((\tau + \mathcal{S}^*V_1 \mathcal{S} V_2)(\tau + \mathcal{S} V_2 \mathcal{S}^* V_1) \big)}\bigg)\bigg|_{\tau = \abs{\zeta}^2,\, \eta = 0}. 
\end{equation}
To that end, let $t \mapsto A(t)$ be a differentiable map with values in the positive definite matrices. Then we have the well-known identity
\begin{equation} \label{eq:derivative_avg_log} 
 \pt_t \avg{\log A(t)} = \avg{A(t)^{-1} \pt_t A(t)} 
\end{equation}
 (see e.g.\ \cite[Lemma~1.1]{Brown1986}). 
We apply the relation \eqref{eq:derivative_avg_log} to $A = (UU^*)^{-1}$ with $t = \eta$ and obtain 
\[ \frac{1}{2} \avg{\pt_\eta \log(UU^*)^{-1}} = \Re \avg{U \pt_\eta U^{-1}} = \avg{V_1} + \avg{V_2} + \pt_\eta \avg{V_1 \mathcal{S} V_2}. \] 
Since $\avg{V_1} + \avg{V_2} = 2 \avg{\Im \bs{M}}$, this proves \eqref{eq:L_general_representation} 
due to \eqref{definition of U}, the continuity of   
$V_1(\tau,\eta)$ and $V_2(\tau,\eta)$ at $\eta = 0$, 
 $\lim_{\eta \to \infty} V_1 = \lim_{\eta \to \infty} V_2 = 0$ by \eqref{bounds on V} 
as well as $\lim_{\eta \to \infty} (UU^*)^{-1} (1 + \eta)^{-4} = 1$ by \eqref{definition of U}.

The identity \eqref{eq:L_general_representation} directly shows that $L$ is rotationally symmetric on $\C$.  
Moreover, it implies that $L$ is a continuous function of $\zeta$ on $\C$ since $V_1(\tau)$ and $V_2(\tau)$ are continuous functions of $\tau$. 

We now show that $\tau \mapsto L(\sqrt{\tau})$ is continuously differentiable on $(0,1) \cup (1,\infty)$ with 
\begin{equation} \label{eq:pt_tau_L} 
 \pt_\tau L(\sqrt{\tau}) = -  \frac{1}{2} \begin{cases} \avg{U(\tau)}, & \text{ if }\tau < 1, \\ \tau^{-1} , & \text{ if } \tau >1. \end{cases} 
\end{equation}
If $\tau <1$ then the continuous differentiability follows from the analyticity and positivity of $V_1$ and $V_2$. Moreover, from \eqref{eq:derivative_avg_log} with $A = (UU^*)^{-1}$ and $t = \tau$, we get 
\[ \frac{1}{2} \pt_\tau \avg{\log (U(\tau)U(\tau)^*)^{-1}} - \pt_\tau \avg{V_1(\tau) \mathcal{S}V_2(\tau)} = \avg{U(\tau)}, \] 
which implies the first case in \eqref{eq:pt_tau_L} due to \eqref{eq:L_general_representation}. 
If $\abs{\zeta} \geq 1$ then $\lim_{\eta \downarrow 0} V_1 = \lim_{\eta \downarrow 0} V_2 = 0$. Hence, we get from \eqref{eq:L_general_representation} that $L(\zeta) = - \log \abs{\zeta}$. 
Thus, the differentiability and the relation \eqref{eq:pt_tau_L} for $\tau > 1$ follows. 
This completes the proof of \eqref{eq:pt_tau_L}. 

Since $L$ is rotationally symmetric it suffices to show \eqref{eq:sigma_integration_part_def_L} under the same constraint on $f$. 
If $f \in C_0^2(\C)$ is rotationally symmetric then a simple change of coordinates yields 
\begin{equation} \label{eq:integral_Delta_f_L_symmetric} 
\frac{1}{2\pi} \int_\C \Delta f (\zeta) L(\zeta) \dd^2 \zeta = 2 \int_0^\infty ( \tau\pt_\tau^2 f(\sqrt{\tau}) + \pt_\tau f(\sqrt{\tau}) )  L(\sqrt{\tau}) \dd \tau,  
\end{equation}
where we employed $\Delta f(\zeta) = 4 \big( \tau \pt_\tau^2 f(\sqrt{\tau}) + \pt_\tau f(\sqrt{\tau})\big)|_{\tau = \abs{\zeta}^2}$. 

We now split up the $\tau$-integration into $(0,1)$ and $(1,\infty)$ and use the differentiability of $L$ on both domains to integrate by parts. 
More precisely, integrating by parts twice, using the continuity of $L$ and $L(1) = 0$ as well as \eqref{eq:pt_tau_L} and $\lim_{\tau \uparrow 1} \pt_\tau L(\sqrt{\tau}) = - 1/2$ yield 
\[ 2 \int_0^1 ( \tau\pt_\tau^2 f(\sqrt{\tau}) + \pt_\tau f(\sqrt{\tau}) )  L(\sqrt{\tau}) \dd \tau = f(1) - \int_0^1 f(\sqrt{\tau}) \pt_\tau \big( \tau \avg{U(\tau)} \big) \dd \tau = 
f(1) - \pi \int_0^1 f(\sqrt{\tau}) \sigma(\sqrt{\tau}) \dd \tau.  \] 
Here, we used in the last step that $\pi\sigma(\zeta) = \pt_\tau \big( \tau \avg{U(\tau)}\big)|_{\tau = \abs{\zeta}^2}$ if $\abs{\zeta} < 1$ due to the definition of $\sigma$ in \eqref{definition of sigma} 
and the definition of $U$ in \eqref{definition of U}.

Secondly, an integration by parts, \eqref{eq:pt_tau_L}, the continuity of $L$ and $L(1)=0$ imply 
\[ 2 \int_1^\infty  ( \tau \pt_\tau^2 f(\sqrt{\tau}) + \pt_\tau f(\sqrt{\tau}) ) L(\sqrt{\tau}) \dd \tau   = - f(1). \] 
By plugging these identities into \eqref{eq:integral_Delta_f_L_symmetric}, we obtain 
\[ - \frac{1}{2\pi} \int_\C \Delta f (\zeta) L(\zeta) \dd^2 \zeta =  \int_{\mathbb{D}} f(\zeta) \sigma(\zeta) \dd^2 \zeta = \int_{\C} f(\zeta) \sigma(\zeta) \dd^2 \zeta, \] 
where the last step follows from $\sigma(\zeta) = 0$ if $\abs{\zeta} \geq 1$ by definition (see \eqref{definition of sigma}). 
This proves \eqref{eq:sigma_integration_part_def_L}. 
\end{Proof} 
 
\subsection{Proof of Proposition~\ref{pro:Brown_measure}} \label{sec:proof_Brown_measure} 
 
In this section, we establish Proposition~\ref{pro:Brown_measure}.

\begin{Proof}[Proof of Proposition~\ref{pro:Brown_measure}]  
In the light of Theorem~\ref{thr:Density}, it suffices to show that \eqref{eq:def_Brown_measure} 
holds when $\mu_X(\dd \zeta)$ is replaced by $\sigma(\zeta) \dd^2 \zeta$. 
To that end, let $L(\zeta)$ be defined as in \eqref{eq:sigma_integration_part_def_L}. 
We first show that 
\begin{equation} \label{eq:L_equals_log_D} 
- L(\zeta) = \log D(X-\zeta)
\end{equation} 
 for all $\zeta \in \C$. 
Using \cite[Theorem~11 and Proposition~13 in Chapter~9]{MingoSpeicher} it is easy to see that 
$\bs{M}(\zeta,\eta) := E [(\bs{H}_\zeta - \ii \eta)^{-1}]$ satisfies \eqref{eq:MDE}, 
where $E := \mathrm{id} \otimes \tau \colon \cal M^{2n \times 2n} \to \C^{2n \times 2n}$, 
$\bs{H}_\zeta \in \cal M^{2n \times 2n}$ is defined analogously to \eqref{eq:def_H_zeta} 
with $X$ from \eqref{eq:Kronecker_circular} 
and $\scr S$ is defined as in \eqref{definition scr S and Z} with $\cal S$ and $\cal S^*$ from 
\eqref{eq:def_self_energy_operator_Brown}. 
We introduce the tracial state $\varphi := \avg{\,\cdot\,}\otimes \tau$ on $\cal M^{n\times n}$  
and the matrix $E_{22} \in \C^{2n\times 2n}$ which has the identity matrix in its lower-right $n \times n$-block and vanishes otherwise.  
Thus, the definitions of $\varphi$ and $\bs{M}$ as well as \eqref{Schur for H}  imply 
\begin{equation} \label{eq:relation_Brown_measure_and_M} 
\varphi \2\frac{\eta}{( X-\zeta)^*(X-\zeta) + \eta^2}   
= - 2 \ii \avg{E_{22} \bs{M}(\zeta,\eta) E_{22}} = 
\Im \avg{\bs{M}(\zeta,\eta)}.  
\end{equation} 
We set  $f_\eps(\zeta) := \varphi\2(\log ((X-\zeta)^* (X-\zeta) + \eps^2)^{1/2}) -\log(1+\eps)$  for $\eps >0$ and $\zeta \in \C$ and compute 
\begin{equation} \label{eq:computation_f_eps}  
\begin{aligned} 
f_\eps(\zeta) & = -\int_\eps^\infty \frac{\pt}{\pt \eta} \bigg( \frac{1}{2} \varphi \big( \log ( (X-\zeta)^* (X-\zeta) + \eta^2)\big) - \log ( 1 + \eta) \bigg) \dd \eta \\ 
& = - \int_\eps^\infty  \varphi  \frac{\eta}{( X-\zeta)^*(X-\zeta) + \eta^2}   - \frac{1}{1 + \eta} \,  \dd \eta \\ 
& = -\int_\eps^\infty \Im \avg{\bs{M}(\zeta,\eta)} - \frac{1}{1 + \eta} \, \dd \eta. 
\end{aligned} 
\end{equation} 
We remark that the integrals exist due to \eqref{eq:integral_im_M}. 
In \eqref{eq:computation_f_eps}, we used \eqref{eq:derivative_avg_log} for $\varphi$ instead 
of $\avg{\,\cdot\,}$ in the second step and \eqref{eq:relation_Brown_measure_and_M} in the third step.
Sending $\eps \downarrow 0$ this shows \eqref{eq:L_equals_log_D} by \eqref{eq:def_fuglede_kadison_det} and \eqref{eq:integral_im_M}. 

By Lemma~\ref{lem:sigma_as_dist_derivative} and standard results from potential theory (see e.g.\ \cite[Chapter~4.3]{ArmitageGardiner}), we know that 
\begin{equation} \label{eq:log_potential_of_sigma} 
 \int_\C \log \abs{\lambda-\zeta}\,  \sigma(\zeta) \dd^2 \zeta = - L(\lambda) + h(\lambda) 
\end{equation}
for all $\lambda \in \C$ and some harmonic function $h \colon \C \to \C$. 
In the proof of Lemma~\ref{lem:sigma_as_dist_derivative}, we saw that $L(\lambda) = - \log\abs{\lambda}$ if $\abs{\lambda}$ is sufficiently large. Hence, $h(\lambda) \to 0$ if $\abs{\lambda} \to \infty$, which implies $h \equiv 0$. 
Therefore, \eqref{eq:L_equals_log_D} and \eqref{eq:log_potential_of_sigma} with $h \equiv 0$ 
prove \eqref{eq:mu_X_equals_sigma} and, thus, Proposition~\ref{pro:Brown_measure}. 
\end{Proof}

 \section{Local inhomogeneous circular law}\label{sec:proof_local_law_X} 

This section is devoted to the proof of Theorem~\ref{thr:Local inhomogeneous circular law} which is based on 
the next theorem. 
Its formulation and the notation in the next arguments is simplified by the use of the following notion of high probability estimate first introduced in \cite{EJP2473}. 

\begin{definition}[Stochastic domination] 
Let $X = X^{(n)}$ and $Y = Y^{(n)}$ be two sequences of two non-negative random variables. 
We say that $X$ is \emph{stochastically dominated} by $Y$, denoted by $X \prec Y$, if,  
for any $\eps>0$ and $\nu>0$, there is $C \equiv C_{\eps,\nu}$ such that 
\begin{equation} \label{eq:def_stochastic_domination} 
 \P \big( X > n^\eps Y \big) \leq C_{\eps, \nu} n^{-\nu} 
\end{equation}
for all $n \in \N$. 
\end{definition} 

We remark that stochastic domination is compatible with basic arithmetic operation (see e.g.\ \cite[Lemma~4.4]{EJP2473}). 
The constants $C_{\eps,\nu}$ in \eqref{eq:def_stochastic_domination} will typically depend on the model parameters. 

To simplify the formulation of the next result, we fix $\tau_* \in (0,\varrho(\mathcal{S}))$ and define the spectral domains 
\[ \mathbb{D}_< = \{ \zeta \in \C \colon \abs{\zeta}^2 \leq  \varrho(\mathcal{S}) - \tau_* \}, \qquad \mathbb{D}_> = \{ \zeta \in \C \colon \tau_* \leq \abs{\zeta}^2 - \varrho(\mathcal{S}) \leq 1/\tau_*\}.  \] 

\begin{theorem}[Local law for $\bs{H}_\zeta$] \label{thm:local_law_H_zeta} 
Let $\eps \in (0,1)$, \
 $X$ satisfy \ref{assum:A1} -- \ref{assum:A3} and $\bs{M}$ be the solution of \eqref{eq:MDE}. 
Then we have the isotropic local law, 
\begin{equation} \label{eq:isotropic_local_law} 
\abs{\scalar{\bs{x}}{(\bs{G}(\zeta,\eta) - \bs{M}(\zeta,\eta))\bs{y}}} \prec \norm{\bs{x}} \norm{\bs{y}} 
\begin{cases} \frac{1}{\sqrt{n\eta}}, & \text{if }\zeta \in \mathbb{D}_<, ~\eta \in [n^{-1 + \eps}, 1], \\ \frac{1}{\sqrt{n}} , & \text{if }\zeta \in \mathbb{D}_> ,~ \eta \in [n^{-1+ \eps},1], \\ \frac{1}{\eta^2 \sqrt{n}}, & \text{if }\zeta \in \mathbb{D}_< \cup \mathbb{D}_>, ~\eta \in [1,n^{100}], \end{cases}  
\end{equation} 
uniformly for all deterministic vectors $\bs{x}$, $\bs{y} \in \C^{2n}$. Moreover, the averaged local law 
 \begin{equation}   \label{eq:averaged_local_law} 
\abs{\avg{\bs{R}(\bs{G}(\zeta, \eta)-\bs{M}(\zeta, \eta))}} \prec \norm{\bs{R}}  
\begin{cases} \frac{1}{{n\eta}}, & \text{if }\zeta \in \mathbb{D}_<, ~\eta \in [n^{-1 + \eps}, 1], \\ \frac{1}{{n}} , &\text{if } \zeta \in \mathbb{D}_> ,~ \eta \in [n^{-1+ \eps},1], \\ \frac{1}{\eta^2 n}, & \text{if }\zeta \in \mathbb{D}_< \cup \mathbb{D}_>,~ \eta \in [1,n^{100}] \end{cases}  
\end{equation} 
holds uniformly for all deterministic matrices $\bs{R} \in \C^{2n\times 2n}$.  
\end{theorem}

We will prove Theorem~\ref{thm:local_law_H_zeta} in Section~\ref{sec:proof_local_law_H} below. 
The next lemma is an application of Theorem~\ref{thm:local_law_H_zeta} and estimates the number of small, in modulus, eigenvalues of $\bs{H}_\zeta$. 
It will be used in the proof of Theorem~\ref{thr:Local inhomogeneous circular law} to control the integral in \eqref{eq:log_det_int_im_G} for small $\eta$.

\begin{lemma}[Number of small singular values of $X-\zeta$] \label{lem:number_eigenvalues} 
Let $X$ satisfy \ref{assum:A1} and \ref{assum:A2}. 
Then, for each $\eps >0$, we have 
\begin{equation} \label{eq:upper_bound_number_eigenvalues} 
 \# \big\{ i \in \llbracket 2n \rrbracket \colon \abs{\lambda_i(\zeta)} \leq \eta \big\} \prec n \eta 
\end{equation} 
uniformly for all $\eta \in [n^{-1 + \eps},n^{100}]$ and $\zeta \in \mathbb{D}_<$. 
\end{lemma}

\begin{Proof}
We follow the proof of Lemma~\ref{lem:number_small_singular_values_global} and use $\abs{\tr \bs{G}} \prec n$ for all $\eta \geq n^{-1 + \eps}$ due to \eqref{eq:averaged_local_law} instead of 
$\abs{\tr \bs{G}} \lesssim n$. This proves Lemma~\ref{lem:number_eigenvalues}. 
\end{Proof}

\begin{Proof}[Proof of Theorem~\ref{thr:Local inhomogeneous circular law}]
 We first remark that the condition $\norm{\Delta f}_{\rm{L}^{1+\beta}} \leq n^{D}\norm{\Delta f}_{\rm{L}^1}$ 
is not needed in Theorem~\ref{thr:Local inhomogeneous circular law} if the stronger Assumption~\ref{assum:bounded_density} holds. 
This can be seen by following the proof of \cite[Theorem~2.5]{Altcirc}
and using, in the proof of \cite[Lemma~5.8]{Altcirc}, Proposition~\ref{pro:smallest_singular_value_general} below instead of \cite[Proposition~5.7]{Altcirc}, \eqref{eq:upper_bound_number_eigenvalues} 
instead of \cite[Eq.~(5.22)]{Altcirc} and  
\eqref{eq:averaged_local_law} instead of \cite[Eq.~(5.4)]{Altcirc}. 

We now prove Theorem~\ref{thr:Local inhomogeneous circular law} assuming \ref{assum:A1}--\ref{assum:A4}.  
In fact, the proof is a simple refinement of the proof of Theorem~\ref{thm:global_law_for_X}
and we solely describe the necessary modifications. 
We replace $f$ by $f_{\zeta_0,\alpha}$ and choose $\Omega = \mathbb{D}_{\sqrt{\varrho(\mathcal{S})}-\tau_*/2}$. 
We remark that $\supp f_{\zeta_0,\alpha} \subseteq \Omega$ for all sufficiently large $n$ as $\alpha>0$. 
The functions $F$ and $h$ as well as the measure $\mu$ are defined analogously according to the new choices of $f$ and $\Omega$.  

In contrast to the proof of Theorem~\ref{thm:global_law_for_X}, we formulate all estimates in the proof of Theorem~\ref{thr:Local inhomogeneous circular law} 
with respect to stochastic domination $\prec$. 
In particular, analogously to \eqref{eq:error_term_sampling}, we obtain 
\begin{equation} \label{eq:error_term_sampling_local} 
 \absbb{\int F \dd \mu - \frac{1}{m} \sum_{i=1}^m F(\xi_i)} \prec n^{-A} \norm{\Delta f}_{\rm{L}^{1 + \beta}} 
\end{equation} 
for all $A>0$, 
where $m$ was chosen sufficiently large and $\xi_1$, \ldots, $\xi_m$ are independent random variables distributed according to $\mu$. 

The next step is proving that, for $T = n^{100}$ and for each $\eps >0$, we have 
\begin{equation} \label{eq:F_estimate_fixed_zeta_local} 
 \abs{F(\zeta)}\prec n^{-1 + \eps} \abs{\Delta f_{\zeta_0,\alpha} (\zeta)}
\end{equation} 
uniformly for all $\zeta \in \Omega$. 
This is the analogue of \eqref{eq:F_estimate_fixed_zeta} and shown by decomposing $h = h_1 + \ldots + h_4$, where 
$h_1$, \ldots, $h_4$ are defined as before but with the choice $\eta_* = n^{-1 + \eps}$. 
As in the proof of Theorem~\ref{thm:global_law_for_X}, we see that $\abs{h_3} \prec n^{-10}$ and $\abs{h_4} \prec n^{-1 + \eps}$ uniformly for $\zeta \in \Omega$.  
To establish $\abs{h_1(\zeta)} \prec n^{-1 + \eps}$, we distinguish the regimes $\eta \in [n^{-1 + \eps}, 1]$ and $\eta \in [1, T]$ in the integral as well as apply a union bound and a continuity argument in $\eta$ 
to \eqref{eq:averaged_local_law} with $\bs{R} = 1$. 
For the bound $\abs{h_2} \prec n^{-1+\eps}$, we decompose the sum into three regimes, $\abs{\lambda_j} < n^{-1 + \eps}$, $\abs{\lambda_j} \in [n^{-1+\eps}, n^{-1/2}]$ and $\abs{\lambda_j} > n^{-1/2}$, 
(instead of two regimes in the proof of Theorem~\ref{thm:global_law_for_X}) and estimate each regime separately.  
The first and the third regime are treated as \eqref{eq:bound_h_2_first_regime} and \eqref{eq:bound_h_2_last_regime}, respectively, using Lemma~\ref{lem:number_eigenvalues} 
instead of Lemma~\ref{lem:number_small_singular_values_global}. 
In the second regime, we restrict to the nonnegative eigenvalues of $\bs{H}_\zeta$ due to $\spec\bs{H}_\zeta = - \spec \bs{H}_\zeta$. We decompose $[n^{-1+\eps},n^{-1/2}]$ dyadically into intervals $[\eta_k,\eta_{k+1}]$
with $\eta_k := 2^k n^{-1 + \eps}$ and obtain 
\[ \frac{1}{4n} \sum_{\abs{\lambda_j} \in [n^{-1+\eps}, n^{-1/2}]} \log \bigg( 1 + \frac{n^{-2 + 2\eps}}{ \lambda_j^2} \bigg) \leq \frac{1}{2n} \sum_{k=0}^K \sum_{\lambda_i \in [\eta_k, \eta_{k+1}]} \log \bigg(1 + \frac{n^{-2 + 2\eps}}{\lambda_i^2} \bigg) \prec \frac{n^\eps}{n},   \] 
where $K = O(\log n)$. 
In the last step, we used the monotonicity of the logarithm, $\log(1 + x ) \leq x$ and $\# \{ j \colon \lambda_j \in [\eta_k, \eta_{k+1}]\} \leq \#\{ j \colon \abs{\lambda_j} \leq \eta_{k+1} \} \prec n^\eps 2^{k+1}$ due to \eqref{eq:upper_bound_number_eigenvalues}. 
This completes the proof of $\abs{h_2} \prec n^{-1  +\eps}$ and, thus, the one of \eqref{eq:F_estimate_fixed_zeta_local}. 

Therefore, following the remaining steps in the proof of Theorem~\ref{thm:global_law_for_X} yields 
\begin{equation} \label{eq:local_law_inside_final_estimate} 
\absbb{\frac{1}{n} \sum_{\xi \in \spec X} f_{\zeta_0,\alpha} (\xi) - \int_\C f_{\zeta_0,\alpha} (\zeta) \sigma(\zeta) \dd^2 \zeta} \prec n^{-1 +2 \alpha} \norm{\Delta f}_{\mathrm{L}^1}  + n^{-A} \norm{\Delta f}_{\rm{L}^{1+\beta}}
\end{equation}
for all $A>0$. 
Using the condition $\norm{\Delta f}_{\rm{L}^{1+\beta }}\leq n^D \norm{\Delta f}_{\rm{L}^1}$ in \eqref{eq:local_law_inside_final_estimate} as well as choosing $A$ appropriately complete the proof. 
\end{Proof}

\subsection{Local law for $\mathbf{H}_\zeta$ -- Proof of Theorem~\ref{thm:local_law_H_zeta}} 
\label{sec:proof_local_law_H} 
This section is devoted to the proof of Theorem~\ref{thm:local_law_H_zeta}. The local law for Hermitian random matrices with decaying correlations was established in \cite{AjankiCorrelated,Erdos2017Correlated}. In order to get the  isotropic version stated in Theorem~\ref{thm:local_law_H_zeta} we will follow the strategy from \cite{Erdos2017Correlated}. Its main result, \cite[Theorem~2.2]{Erdos2017Correlated} is not directly applicable to our current situation since Assumption~(E) from \cite{Erdos2017Correlated} is violated for $\bf{H}_\zeta$. The reason why Assumption~(E) is needed in the proof of \cite[Theorem~2.2]{Erdos2017Correlated} is to ensure the invertibility of $\scr{L}$ in the stability result \cite[Theorem~5.2]{Erdos2017Correlated} for the MDE. The purpose of this section is to show how the proof is adjusted by using our new stability  results,  Theorem~\ref{thr:Stability} and Corollary~\ref{cor:perturbations}, instead. 

The resolvent $\bs{G}=(\bs{H}_\zeta-\ii \eta)^{-1}$ satisfies the perturbed MDE
\bels{perturbed MDE with error}{
 1+ (\ii \eta + \bs{Z} + \scr{S}\bs{G})\bs{G} = \bs{D}\,, \qquad \bs{D}:= (\bs{H}_\zeta+\bs{Z} + \scr{S}\bs{G})\bs{G}\,.
}
The main input for the local law  for $\bs{H}_\zeta$,  Theorem~\ref{thm:local_law_H_zeta}, is the following estimate on the error term $\bs{D}$ in terms of the $p$-norms for random variables $Y$ and random matrices $\bs{A} \in \C^{2n \times 2n}$ defined through
\[
\norm{Y}_p := (\E\abs{Y}^p)^{1/p}\,, \qquad \norm{\bs{A}}_p := \sup_{\norm{\bs{x}},\norm{\bs{y}}\le 1}\pb{\E\abs{\scalar{\bs{x}}{\bs{A}\bs{y}}}^p}^{1/p}.
\] 
\begin{proposition}[Bound on error matrix] \label{prp:Bound on error matrix}
There is a constant $C>0$, depending only on model parameters, such that for any  $\eta \in [n^{-1},n^{100}]$, $p\in \N$, $\eps>0$, $\bs{R} \in \C^{2n \times 2n}$ and $\bs{x},\bs{y} \in \C^{2n}$ with $q:=Cp^4/\eps$ the following holds true:
\begin{align}
\label{D isotropic bound}
\norm{\scalar{\bs{x}}{\bs{D}\bs{y}}}_p &\lesssim_{\eps,p} \norm{\bs{x}}\norm{\bs{y}}n^\eps \sqrt{\frac{\norm{\im \bs{G}}_q}{n \eta}}\pb{1+\norm{\bs{G}}_q}^C\pbb{1+\frac{\norm{\bs{G}}_q}{n^{1/2-\eps}}}^{Cp}\,,
\\ \label{D averaged bound}
\norm{\scalar{\bs{R}}{\bs{D}}}_p &\lesssim_{\eps,p} \norm{\bs{R}}n^\eps  (1+\eta) {\frac{\norm{\im \bs{G}}_q}{n \eta}}\pb{1+\norm{\bs{G}}_q}^C\pbb{1+\frac{\norm{\bs{G}}_q}{n^{1/2-\eps}}}^{Cp}\,.
\end{align}
\end{proposition}

Before deriving Proposition~\ref{prp:Bound on error matrix} from \cite[Theorem~4.1]{Erdos2017Correlated}, we now explain the definition of the self-energy operator in \cite{Erdos2017Correlated} which differs from 
the self-energy operator $\scr S$ used in the present work and defined in \eqref{definition scr S and Z}. 
Instead of $\scr S$, the self-energy operator considered in \cite{Erdos2017Correlated} 
(and denoted by $\cal S$ in \cite{Erdos2017Correlated}) is 
\bels{definition wt scr S}{
\wt{\scr{S}}\bs{R}:= \E(\bs{H}_\zeta+\bs{
Z})\bs{R}(\bs{H}_\zeta+\bs{Z})= \mtwo{\cal{S}R_{22} & \cal{R} R_{21}}{ \cal{R}^*R_{12} & \cal{S}^*R_{11}}\,,
}
with $\bs{Z}=\bs{Z}(\zeta,\ol{\zeta})$ and  the operators $\cal{R},\cal{R}^*:\C^{n \times n} \to \C^{n \times n}$ defined through
\bels{definition cal T}{
\cal{R}R:=\E XRX\,, \qquad \cal{R}^*R:=\E X^*RX^*\,.
}
Moreover, \cite{Erdos2017Correlated} works with the solution $\wt{\bs{M}}$ (denoted by $M$ in \cite{Erdos2017Correlated}) of the MDE, \eqref{eq:mde_general_spectral_parameter}, with the self-energy $\wt{\scr S}$ instead of ${\scr S}$, i.e.\ 
$\wt{\bs{M}}=\wt{\bs{M}}(\zeta,z)$ satisfies 
\begin{equation}\label{eq:MDE_tilde} 
-\wt{\bs{M}}^{-1} = z\bs{1} + \bs{Z}+ \wt{\scr{S}}\;\wt{\bs{M}}
\end{equation} 
for all $z \in \bb{H}$ and all $\zeta \in \C$, where $\bs{Z} = \bs{Z}(\zeta,\ol\zeta)$.

\begin{Proof} The bounds \eqref{D isotropic bound} and  \eqref{D averaged bound} are an immediate consequence of \cite[equations (23a) and (23b)]{Erdos2017Correlated},  respectively,   with the choice $\mu=1/2-\eps$. We will use the following lemma.
\begin{lemma}\label{lmm:A2 implies C and D}
Let $X$ satisfy Assumptions~\ref{assum:A1} and \ref{assum:A2}, then $\bs{H}_\zeta$ defined in \eqref{Hzeta} satisfies \cite[Assumption~(C)]{Erdos2017Correlated}   and the following modification of \cite[Assumption~(D)]{Erdos2017Correlated}.  \\
Modification of \cite[Assumption~(D)]{Erdos2017Correlated}: With the notation from the formulation of \cite[Assumption~(D)]{Erdos2017Correlated}   the matrix $\bs{H}_\zeta$ satisfies
\bels{modified assumption D}{
\kappa(f, g_1, \dots, g_q) \le_{R,q,\mu} n^{-3q} \norm{f}_{2q}\prod_{j=1}^q\norm{g_j}_{2q}\,,
}
i.e.\ the $\norm{\2\cdot\2}_{q+1}$-norms on the right hand side of  \cite[Assumption~(D)]{Erdos2017Correlated}  are replaced by  $\norm{\2\cdot\2}_{2q}$-norms. This change does not effect any of the proof in \cite{Erdos2017Correlated}.
\end{lemma}

The proof of Lemma~\ref{lmm:A2 implies C and D} is given in  Appendix~\ref{apx:Auxiliary results} below.  
The matrix $\bs{H}_\zeta$ satisfies \cite[Assumptions~(A),(B),(C)]{Erdos2017Correlated} and the modified version \eqref{modified assumption D} of \cite[Assumption~(D)]{Erdos2017Correlated} according to Assumption~\ref{assum:A1}, \ref{assum:A2} and Lemma~\ref{lmm:A2 implies C and D}. Since the modification  \eqref{modified assumption D} does not effect any of the proofs in \cite{Erdos2017Correlated} we can apply \cite[Theorem~4.1]{Erdos2017Correlated} to $\bs{H}_\zeta$. 
 Owing to the different self-energy operator in \cite{Erdos2017Correlated} as explained above,  
the bounds from \cite[Theorem~4.1]{Erdos2017Correlated}\footnote{Note that there is a typo in the statement of \cite[equation~(23b)]{Erdos2017Correlated}. Compared  to \eqref{D averaged bound} the first $(1+\eta)$-factor on the right hand side was missing. Indeed, the right hand side of \cite[equation~(23b)]{Erdos2017Correlated} should have been multiplied with a factor $\avg{z}:=(1+\abs{z})$. In the arXiv version of \cite{Erdos2017Correlated} this typo was corrected.} are for 
\bels{definition wt D}{
\wt{\bs{D}}:=\bs{D}+(\wt{\scr{S}}\bs{G}-\scr{S}\bs{G})\bs{G}.
}
Thus to prove the proposition it suffices to show the following bounds on the additional error term
\begin{align}
\label{additional error isotropic bound}
\norm{\scalar{\bs{x}}{((\wt{\scr{S}}-\scr{S})\bs{G})\bs{G}\bs{y}}}_p  &\lesssim \norm{\bs{x}}\norm{\bs{y}}n^{\eps}\norm{\bs{G}}_q  \pbb{\frac{\norm{\im \bs{G}}_q}{n \eta}}^{1/2},
\\
\label{additional error averaged bound}
\abs{\scalar{\bs{R}}{((\wt{\scr{S}}-\scr{S})\bs{G})\bs{G}}}  &\lesssim \norm{\bs{R}}\frac{\avg{\im \bs{G}}}{n \eta}\,.
\end{align}
To see \eqref{additional error isotropic bound} we use that for any unit vectors $x,y \in \C^n$ and $R,Q \in \C^{n \times n}$ we have 
\[
\norm{\scalar{x}{(\cal{R}R){Q}y}}_p\le  \normb{{\textstyle \sum_{j,k}}(Rv_{jk})_j (Qy)_k}_p \le n^{2\eps} \norm{R}_{1/\eps}
\normb{{\textstyle \sum_{j,k}}\norm{v_{jk}}  \abs{(Qy)_k}}_{2p}
\,,
\]
where $v_{jk} :=\p{\sum_{i}\ol{x}_i\E X_{ij}X_{lk} }_{l} \in \C^n$ and we employed the general inequality for random variables $(X_i,Y_i)_{i=1}^{n^2}$ and $\eps \in (0,1/2p)$ of the form $\norm{\sum_i X_i Y_i}_p\le n^{2\eps}\sup_{i}\norm{X_i}_{1/\eps}\norm{\sum_i \abs{Y_i}}_{2p}$. Since the diagonal  
contributions of $\wt{\scr{S}}$ and $\scr{S}$ coincide, 
 we conclude that 
\[
\norm{\scalar{\bs{x}}{((\wt{\scr{S}}-\scr{S})\bs{G})\bs{G}\bs{y}}}_p \lesssim  
n^{-1/2+2\eps} \norm{\bs{G}}_{1/\eps} \norm{\bs{G}^*\bs{G}}_{2p}\,,
\]
where 
the decay $\norm{v_{jk}} \lesssim n^{-1}(1+d(j,k))^{-\nu}$ of arbitrarily high order $\nu \in \N$ was used. 
The Ward identity $\eta\2\bs{G}^*\bs{G} = \im \bs{G}$ now implies \eqref{additional error isotropic bound}. 

The remaining inequality, \eqref{additional error averaged bound}, follows from the Ward identity and $\norm{\wt{\scr{S}}-\scr{S}}_{\rm{hs}} \lesssim\norm{\cal{R}}_{\rm{hs}}\lesssim \frac{1}{n}$ (cf.\ \eqref{hs norm bound on cal R} for the bound on $\cal{R}$). This finishes the proof of Proposition~\ref{prp:Bound on error matrix}. 
\end{Proof}

\begin{Proof}[Proof of Theorem~\ref{thm:local_law_H_zeta}] 
To prove the theorem we follow the arguments from the proof of \cite[Theorems~2.1 and 2.2]{Erdos2017Correlated} in \cite[Sections~5.3 and 5.4]{Erdos2017Correlated} line by line. 
The  spectral parameter $\zeta \in \mathbb{D}_< \cup \mathbb{D}_>$ associated to $X$ (cf.\ \eqref{Hzeta}) is fixed throughout the proof. 
 The different definition of the self-energy in \cite{Erdos2017Correlated} as explained 
after Proposition~\ref{prp:Bound on error matrix} necessitates replacing a few objects in the arguments in 
\cite{Erdos2017Correlated} by their counterparts in the present setup. Indeed, 
$ \cal S$, $M$ and $D$ in \cite{Erdos2017Correlated} are 
replaced by $\scr S$, $\bs{M}$ and $\bs{D}$ from \eqref{eq:def_S_intro},  \eqref{definition of M} and \eqref{perturbed MDE with error}, respectively. 
The role of the spectral parameter $z$ in \cite{Erdos2017Correlated} is played here by $\ii \eta$ 
which is associated to $\bs{H}_\zeta$. Correspondingly the domains $\D_{\rm{out}}^\delta$  and $\D_{\gamma}^\delta$ in $\C$ from \cite{Erdos2017Correlated} are replaced by 
\[
\D_{\rm{out}}^\delta :=\{\ii\eta: \eta \in (0,n^{100}]\,, \,\eta + (\abs{\zeta}-1)_+ \ge \delta \} \; \text{ and }\;  \D_{\gamma}^\delta:=\{\ii \eta:  \eta \in [n^{-1+\gamma},n^{100}]\,, \, \eta +\abs{\abs{\zeta}-1}\ge \delta \}\,,
\] 
respectively. Here, $(\xi)_+:=\max\{0,\xi\}$ denotes the positive part. 

Furthermore, whenever \cite[Theorem 4.1]{Erdos2017Correlated} is used in \cite{Erdos2017Correlated} we will use Proposition~\ref{prp:Bound on error matrix} instead.  The now missing Assumption~(E) from \cite{Erdos2017Correlated} was used along the argument solely for the purpose of establishing stability of the MDE, i.e.\ to show that the inverse of $\scr{L}$  defined in \eqref{derivative of scr J}  is bounded  (Note that $\scr{L}$ is the analogue of $1-\cal{C}_M\cal{S}$ from \cite{Erdos2017Correlated}).  We will now point out where the boundedness of $\scr{L}^{-1}$ and the resulting stability in the form of \cite[equation (74)]{Erdos2017Correlated} has to be replaced by the use of Corollary~\ref{cor:perturbations}. 

Any direct use of \cite[equation (74)]{Erdos2017Correlated} is simply replaced by \eqref{Global stability}, 
using that $\bs{G} \in \bs{E}_-^\perp$ 
 by Lemma~\ref{lem:symmetry_inverse_block_matrix} below.  
Otherwise the boundedness of $\scr{L}^{-1}$ is only used to establish the averaged bound \cite[equation (84)]{Erdos2017Correlated}. To establish this bound in the current setting we start from the quadratic equation 
\bels{Delta difference equation}{
\scr{L}\bs{\Delta} = \wh{\bs{D}}\,, \qquad \wh{\bs{D}}:=-\bs{M}\bs{D}+\bs{M}(\scr{S}\bs{\Delta})\bs{\Delta}
}
for the difference $\bs{\Delta}:= {\bs{G}}-{\bs{M}}$ right after \cite[equation (83)]{Erdos2017Correlated}. 

Away from the self-consistent spectrum $\supp \rho_\zeta$ with $\rho_\zeta$ defined in \eqref{eq:def_rho_zeta},  i.e.\ for $\ii \eta \in \D_{\rm{out}}^\delta$ we can invert $\scr{L}$ and  follow the argument from \cite{Erdos2017Correlated}  exactly since
$\norm{\scr{L}^{-1}}_{\rm{hs}}\lesssim_\delta 1$. This bounded invertibility of $\scr{L}$ follows from \cite[Lemma~3.7]{Kronecker} because  $\dist(\ii \eta, \supp \rho_\zeta)\gtrsim_\delta 1$ by Corollary~\ref{crl:scDOS inside disk}.   
In the regime $\abs{\zeta} \le 1-\delta$ and $\eta \le \delta$ the operator $\scr{L}$ does not have a bounded inverse. Thus, we have to proceed more cautiously from \eqref{Delta difference equation} and use the operator $\scr{M}$ defined in \eqref{derivative of scr J}. 
Since $\bs{\Delta} \in \bs{E}_-^\perp$ and $\scr{M}\scr{L}$ preserve the subspace $\bs{E}_-^\perp$, we see that $\scr{M}\wh{\bs{D}} \in \bs{E}_-^\perp$ by acting with $\scr{M}$ on both sides of \eqref{Delta difference equation}. Therefore we can use \eqref{bound on L inv} to invert $\scr{L}$ on $\scr{M}^{-1}\bs{E}_-^\perp$ and after that follow \cite{Erdos2017Correlated} again until the end of \cite[Step~3 in Section~5.4]{Erdos2017Correlated}. This proves Theorem~\ref{thm:local_law_H_zeta} in the regime $\D_{\gamma}^\delta$ for any $\gamma>0$ without the $\eta^{-2}$-decay of the bound in the regime $\eta \ge 1$ on the right hand side of \eqref{eq:averaged_local_law} and \eqref{eq:isotropic_local_law}. 

For the $\eta^{-2}$-decay we replace \cite[Step~4 in Section~5.4]{Erdos2017Correlated} by Lemma~\ref{lem:no_eigenvalues_H}, the analogue of \cite[Corollary~2.3]{Erdos2017Correlated}, to see that there are no eigenvalues in a $\eps$-neighbourhood of the origin for $\zeta \in \D_{>}$ and follow \cite[Step~5 in Section~5.4]{Erdos2017Correlated}, again using Proposition~\ref{prp:Bound on error matrix} instead of \cite[Theorem 4.1]{Erdos2017Correlated} and \eqref{Global stability} instead of \cite[equation (74)]{Erdos2017Correlated}. 
This finishes the proof of Theorem~\ref{thm:local_law_H_zeta}. \end{Proof}

\subsection{Eigenvector delocalisation for $X$} \label{subsec:proof_consequences_local_law} 

In this subsection we prove Corollary~\ref{cor:delocalization}   
which is a consequence of the local law for $\bs{H}_\zeta$, Theorem~\ref{thm:local_law_H_zeta}.

\begin{Proof}[Proof of Corollary~\ref{cor:delocalization}] 
Take $v \in\C^n$ and $\eps >0$. Let $u \in \cal U_{\tau_*}$. Then there is $\zeta \in \mathbb{D}_<$ such that 
$X u = \zeta u$. With $\bs{u} := (0,u)^t \in \C^{2n}$, we obtain $\bs{H}_\zeta \bs{u} = 0$. 
Extending $\bs{u}/\norm{\bs{u}}$ to an orthonormal basis $\bs{u}/\norm{\bs{u}}$, $\bs{u}_2$, \ldots, $\bs{u}_{2n}$
of $\C^{2n}$ consisting of eigenvectors of $\bs{H}_\zeta$ associated to the eigenvalues $\lambda_1(\zeta) = 0$, 
$\lambda_2(\zeta)$, \ldots, $\lambda_{2n}(\zeta)$ and using the spectral theorem for any $\bs{v} \in\C^{2n}$ 
and $\eta >0$ yields 
\begin{equation} \label{eq:proof_delocalization_X_aux1}  
\Im \scalar{\bs{v}}{\bs{G}(\zeta,\eta)\bs{v}} = \frac{\abs{\scalar{\bs{v}}{\bs{u}}}^2}{\eta \norm{\bs{u}}^2} 
 + \sum_{i=2}^{2n} \frac{\eta \abs{\scalar{\bs{v}}{\bs{u}_i}}^2}{\lambda_i(\zeta)^2 + \eta^2} \geq 
\frac{1}{\eta} \frac{\abs{\scalar{v}{u}}^2}{\norm{u}^2},  
\end{equation}  
where, for the last step, we chose $\bs{v} := (0,v)^t$. 
Thus, for any $\eta >0$, the bound \eqref{eq:proof_delocalization_X_aux1} implies 
\begin{equation} \label{eq:proof_delocalization_X_aux2}  
 \big\{ \exists u \in \cal U_{\tau_*} \colon \abs{\scalar{v}{u}} \geq n^{-1/2 + \eps} \norm{v}\norm{u} \big\} 
\subset \{ \exists \zeta \in \mathbb{D}_< \colon \eta \abs{\scalar{\bs{v}}{\bs{G}(\zeta,\eta) \bs{v}}} \geq 
n^{-1 +  2\eps} \norm{\bs{v}}^2 \big\} 
\end{equation} 
with $\bs{v} := (0,v)^t$. 

From \eqref{eq:isotropic_local_law} in Theorem~\ref{thm:local_law_H_zeta} 
and $\norm{\bs{M}} \lesssim 1$ due to \eqref{eq:norm_M}, we conclude that, for each $\eps \in (0,1)$, 
the bound 
$\abs{\scalar{\bs{v}}{\bs{G}(\zeta,\eta)\bs{v}}} \lesssim \norm{\bs{v}}^2$ holds with very high probability 
uniformly for all $\zeta \in \mathbb{D}_<$ with $\eta = n^{-1 + \eps}$. 
Therefore, a grid- and continuity argument in $\zeta$ shows that $\sup_{\zeta \in \mathbb{D}_<} \abs{\scalar{\bs{v}}{\bs{G}(\zeta,\eta) \bs{v}}} \lesssim \norm{\bs{v}}^2$ with very high probability for $\eta =n^{-1 + \eps}$. 
We conclude that \eqref{eq:proof_delocalization_X_aux2} with $\eta = n^{-1 + \eps}$ and sufficiently small $\eps>0$ 
proves Corollary~\ref{cor:delocalization}.  
\end{Proof}

\section{Bound on the smallest singular value} \label{sec:smallest_singular_value} 

In this section we bound the smallest singular value of $X+A$ if $X$ satisfies Assumption~\ref{assum:bounded_density} and $A$ is deterministic. This is done in Proposition~\ref{pro:smallest_singular_value_general} below, which, in particular, implies Proposition~\ref{pro:smallest_singular_value}. 
Moreover, we prove Lemma~\ref{lem:A4_for_block_matrices} in the next subsection.

We recall that $\smin(R)$ denotes the smallest singular value of a matrix $R \in\C^{n \times n}$. 

\begin{proposition}[Smallest singular value] \label{pro:smallest_singular_value_general} 
Let $X = (x_{ij})_{i,j} \in \C^{n\times n}$ be a correlated random matrix satisfying $\E X = 0$ and \ref{assum:bounded_density}. 
Then, for any deterministic matrix $A \in \C^{n\times n}$, we have  
\[ \P \big( \smin(X+A) \leq u\big) \leq \pi n^{\kappa + 5-3/q} u^{1-1 /q} \] 
for all $u\in (0,1]$.   
\end{proposition} 

\begin{Proof}
 The following proof mimics the one of \cite[Lemma~4.12]{BordenaveChafai2012} that is valid for independent entries. 

Going back to \cite{RudelsonVershynin2008}, the smallest singular value is often estimated by the inequality 
\[ \smin(X+ A) \geq n^{-1/2} \min_{i \in \qq{n}} \dist(R_i, R_{-i}), \] 
where $R_1, \ldots, R_n$ are the rows of $X + A$ and $R_{-i} := \Span\{R_j \colon j \neq i \}$ 
(see also \cite[Lemma~4.16]{BordenaveChafai2012}).

Owing to this inequality and a union, we obtain  
\[ \P (\smin(X + A) \leq u) \leq n \max_{i \in \qq{n}} \P (n^{-1/2}\dist(R_i, R_{-i}) \leq u). \] 
We fix $i \in \qq{n}$. Let $y$ be a unit vector that is orthogonal to $R_{-i}$ and measurable with respect to $\{ R_j \colon j \neq i \}$. 
The Cauchy-Schwarz inequality implies 
\[ \abs{\scalar{R_i}{y}} \leq \norm{\pi_i(R_i)}\norm{y} = \dist(R_i, R_{-i} ), \] 
where $\pi_i$ is the orthogonal projection onto the orthogonal complement of $R_{-i}$. 
Therefore, we obtain 
\[ \P ( \dist(R_i, R_{-i}) \leq u n^{1/2}  )  \leq \P( \abs{\scalar{R_i}{y}} \leq u n^{1/2} ). \] 
Since $y$ is normalised, we find $j \in \qq{n}$ such that $\abs{y_j} \geq n^{-1/2}$. 
This yields 
\begin{equation} \label{eq:proof_singular_value_aux1} 
 \P( \abs{\scalar{R_i}{y}} \leq u n^{1/2} ) = \E \Big[ \E \Big[ \sum_{j\in \qq{n}} \mathbf{1}( j = \min \{ k \colon \abs{y_k} \geq n^{-1/2} \}) \P \big( \abs{\scalar{R_i}{y}} \leq un^{1/2} \bigm\vert \mathfrak{X}_{ij} 
\big) \Bigm\vert y \Big] \Big], 
\end{equation}  
 where we denote by $\frak X_{ij}$ the family of random variables 
$\frak X_{ij} := \big\{ x_{kl} \colon (k,l) \in \llbracket n\rrbracket^2 \setminus \{ (i,j)\} \big\}.$

We now estimate the conditional probability with respect to $\mathfrak X_{ij}$ for any $j \in \qq{n}$ 
such that $\abs{y_j}\geq n^{-1/2}$. 
We only consider the case that $\psi_{ij}$ is a density on $\C$. (If $\psi$ is a density on $\R$ then 
we proceed completely analogously.) 
The condition $\abs{y_j} \geq n^{-1/2}$, the identity \eqref{eq:rep_density} in \ref{assum:bounded_density} and Hölder's inequality imply 
\[  \P\big(\, \abs{\scalar{R_i}{y}} \leq u n^{1/2} \bigm| \mathfrak{X}_{ij} \, \big) = \int_\C \mathbf{1} \Big( \absb{a_{ij} + \overline{z} } \leq \frac{un}{\abs{y_j}} \Big) \psi_{ij}(z) \dd^2 z 
\leq \pi^{(q-1)/q} n^{3 (q-1) /q} u^{2(q-1)/q} \norm{\psi_{ij}}_{q} , \] 
for some $\C$-valued random variable $a_{ij}$, which is measurable with respect to $\mathfrak{X}_{ij}$. 
Thus, estimating the sum in \eqref{eq:proof_singular_value_aux1} by $n$ and using the bound on $\E\norm{\psi_{ij}}_{q}$ from \ref{assum:bounded_density} 
complete the proof of Proposition~\ref{pro:smallest_singular_value}. 
\end{Proof}

\subsection{ Proof of Lemma~\ref{lem:A4_for_block_matrices}}  
\label{sec:check_examples} 

\begin{Proof}[Proof of Lemma~\ref{lem:A4_for_block_matrices}] 
For all $i,j \in \qq{N}$ and $\alpha, \beta \in \qq{K}$, we set  
$\mathfrak{X}_{i,j,\alpha, \beta} := \{ (x_{kl})_{\gamma,\delta} \colon (k,l,\gamma, \delta) \neq (i,j,\alpha, \beta) \}$ 
and 
\[ \psi_{i,j,\alpha,\beta} (z) : = \frac{f_{ij}((z_{\gamma, \delta})_{\gamma, \delta \in \qq{K}})}{
\int_\C f_{ij} ((z_{\gamma, \delta})_{\gamma, \delta \in \qq{K}}) \dd^2 z_{\alpha \beta}} \biggm\vert_{
\begin{aligned} 
\vspace*{-0.5cm} 
\scriptstyle z_{\alpha\beta} & \scriptstyle  
= z,\\ 
\vspace{-0.8cm} 
\scriptstyle  
z_{\gamma\delta} & \scriptstyle  
= x_{ij} \text{ if } (\gamma,\delta) \neq (\alpha, \beta),\end{aligned}} \] 
with the convention $\psi_{i,j,\alpha,\beta} (z) = 0$ if the denominator vanishes. 
A simple computation shows that  
\[ \P \big( \scalar{e_\alpha}{x_{ij}e_\beta} \sqrt{NK} \in B \bigm\vert \mathfrak{X}_{i,j,\alpha, \beta}\big) = \int_B \psi_{i,j,\alpha,\beta} (z) \dd^2 z 
\] 
for all measurable $B \subset \C$, where $e_1$, \ldots, $e_K$ denote the standard basis vectors of $\C^K$. Hence, for each entry $\scalar{e_\alpha}{x_{ij}e_\beta}$ of $X$ as defined in \eqref{eq:def_Kronecker_matrix_X}, we have determined the density in \eqref{eq:rep_density}. 

From the definition of $\psi_{i,j,\alpha, \beta}$, it is easy to conclude that 
\begin{equation} \label{eq:relation_density_kronecker}
\E \norm{\psi_{i,j,\alpha,\beta}}_q = \int_{\C^{K\times K  -1}} \bigg( \int_\C f_{ij}(\mathbf{z})^q \dd^2 z_{\alpha\beta} \bigg)^{1/q} \dd^2 z_{11} \dd^2 z_{12} \ldots \wh{\dd^2 z_{\alpha\beta}} \ldots \dd^2 z_{KK}. 
\end{equation} 
Finally, applying \eqref{eq:condition_density_block_matrix} to \eqref{eq:relation_density_kronecker} implies Assumption~\ref{assum:bounded_density} for $X$.  
\end{Proof}

\appendix

\section{Exclusion of eigenvalues outside disk and global law for $\mathbf{H}_\zeta$} 
\label{app:tools_section_3} 

In this appendix we show how Lemma~\ref{lem:no_eigenvalues_H} and Proposition~\ref{pro:global_law_H_averaged} can be derived from existing results. 
We recall that the self-consistent density of states $\rho_\zeta$ was defined in  
\eqref{eq:def_rho_zeta} and the self-consistent spectrum is $\supp \rho_\zeta$. 

The following corollary to Lemma~\ref{lem:Im_M_sim_eta} states that 
the self-consistent spectrum $\supp \rho_\zeta$ is bounded away from zero 
for any spectral parameter $\zeta$ outside the  disk of radius $\sqrt{\varrho(\cal{S})}$. 

 \begin{corollary}\label{crl:scDOS inside disk}
 Let $\zeta \in \C$ with $\abs{\zeta}^2 \ge \varrho(\cal{S})+\delta$ for some $\delta>0$. 
Assuming \ref{assum:A1}, \ref{assum:A2} and \ref{assum:S_inverse_bound}, the self-consistent spectrum $\supp \rhoDOS_\zeta$ 
is bounded away from zero, i.e.\ $\dist(0, \supp {\rhoDOS}_\zeta)\gtrsim_\delta 1$. 
\end{corollary}

\begin{Proof} 
The corollary follows  from Lemma~\ref{lem:Im_M_sim_eta} and the implication \cite[(i) implies (v) in Lemma~D.1]{Altshape}.
\end{Proof} 

Lemma~\ref{lem:no_eigenvalues_H} and Proposition~\ref{pro:global_law_H_averaged}  will follow from 
\cite[Corollary~2.3]{Erdos2017Correlated} and \cite[Theorem~2.1]{Erdos2017Correlated}, respectively. 
As explained after Proposition~\ref{prp:Bound on error matrix}, the self-energy operator $\wt{\scr S}$ used in \cite{Erdos2017Correlated} (cf.\ \eqref{definition wt scr S}) 
differs slightly from $\scr S$ defined in \eqref{definition scr S and Z} and used in the present work. 
Therefore, applying results from \cite{Erdos2017Correlated} requires controlling the 
difference between $\bs{M}(\zeta,z)$ and $\wt{\bs{M}}(\zeta,z)$, the solutions of the MDE's 
\eqref{eq:mde_general_spectral_parameter} and \eqref{eq:MDE_tilde}, respectively. This is done in the next lemma. 
In analogy to $\rho_\zeta$, we define $\wt{\rho_\zeta}$ as the unique probability measure on $\R$
with Stieltjes transform $z \mapsto \avg{\wt{\bs{M}}(\zeta,z)}$. 

\begin{lemma} [Properties of $\wt{\bs{M}}$] \label{lmm:Properties of whM} 
Assume \ref{assum:A1} and \ref{assum:A2}. Let 
$\bs{M}=\bs{M}(\zeta,\ii\eta)$ and $\wt{\bs{M}}=\wt{\bs{M}}(\zeta,\ii\eta)$ for some $\zeta \in \C$ and $\eta>0$. 
If $\eta \ge n^{-\eps}$ for some small enough $\eps>0$ or $\abs{\zeta}^2\ge \varrho(\cal{S})+\delta$ for some $\delta>0$ then the following holds.
\begin{enumerate}[label=\emph{(\roman*)}] 
\item \label{item:tilde_M_M_operator_norm} The solutions  are close in operator norm: $\norm{\wt{\bs{M}}-\bs{M}}\lesssim  \frac{n^{C\eps}}{(1+ \eta^2)\sqrt{n}}$ 
for some universal constant $C>0$.
\item \label{item:tilde_M_M_hs_norm} The solutions  are close in $\rm{hs}$-norm: $\norm{\wt{\bs{M}}-\bs{M}}_{\rm{hs}}\lesssim  \frac{n^{C\eps}}{(1+ \eta^2){n}}$ for some universal constant $C>0$.
\item \label{item:supp_tilde_rho} If $\abs{\zeta}^2\ge \varrho(\cal{S})+\delta$ and we also assume \ref{assum:S_inverse_bound}, then $0$ is outside the self-consistent spectrum associated to $\wt{\bs{M}}$, i.e.\ $\dist(0, \supp  \wt{\rhoDOS}_\zeta )\gtrsim 1$. 
\end{enumerate}
\end{lemma}

Before establishing Lemma~\ref{lmm:Properties of whM} we use it to show Lemma~\ref{lem:no_eigenvalues_H}. 

\begin{Proof}[Proof of Lemma~\ref{lem:no_eigenvalues_H}]
Given Lemma~\ref{lmm:Properties of whM} \emph{\ref{item:supp_tilde_rho}}, 
Lemma~\ref{lem:no_eigenvalues_H} is a direct consequence of \cite[Corollary~2.3]{Erdos2017Correlated}. 
\end{Proof} 

For the reader's convenience we record the auxiliary result proven in \cite[Lemma~3.4(i)]{Kronecker} and \cite[Lemma~3.7(ii), (iii)]{Kronecker}.

\begin{lemma} \label{lem:stability_outside} 
Assume \ref{assum:A1} and \ref{assum:A2}. 
Let $\scr L=\scr L(\zeta,\eta)$ be defined as in \eqref{derivative of scr J}. 
The following holds. 
\begin{enumerate}[label=\emph{(\roman*)}] 
\item For all $\zeta \in \C$ and $\eta >0$, we have 
\[ \norm{\bs{M}(\zeta,\eta)} \leq \frac{1}{\dist(\ii \eta, \supp \rho_\zeta)}. \] 
\item 
 \label{item:Bound on inverse cal L}
There is a universal constant $K>0$ such that, for all $\eta>0$ and $\zeta \in \C$ with $\abs{\zeta} \lesssim 1$, we have  
\[ \norm{\scr L^{-1}}_{\rm{hs}} + \norm{\scr L^{-1}} + \norm{(\scr L^{-1})^*} \lesssim 1 + \frac{1}{
\dist(\ii \eta, \supp \rho_\zeta)^K}. \] 
\end{enumerate} 
\end{lemma} 

\begin{Proof}[Proof of Lemma~\ref{lmm:Properties of whM}]
We start by showing that the operators $\cal{R}$ and $\cal{R}^*$ from \eqref{definition cal T} that constitute the off-diagonal entries of $\wt{\scr{S}}-\scr{S}$ can be considered small perturbations. 
Indeed, we will prove that 
\bels{op to hs bound for cal T}{
\norm{\cal{R}R} +\norm{\cal{R}^*R} \lesssim \frac{1}{\sqrt{n}}\norm{R}_{\rm{hs}}\le \frac{1}{\sqrt{n}}\norm{R}\,
}
for every $R \in \C^{n\times n}$.  
To check \eqref{op to hs bound for cal T} we simply use that
\[
\abs{\scalar{x}{(\cal{R}R)y}}\le \frac{1}{\sqrt{n}}\bigg(\sum_{i,j,k,l}\abs{x_i y_jx_ly_k}K_{ ij,lk}\bigg)^{1/2} \norm{R}_{\rm{hs}}. 
\]
Here, we introduce the coefficients $K_{ij,lk}$ given by 
\[
K_{ij,lk}:=\sum_{u,v}\absb{\rm{Cov} (\overline{X_{iu}}\sqrt{n} ,X_{vj}\sqrt{n} )\rm{Cov} (\overline{X_{lu}}\sqrt{n} ,X_{vk}\sqrt{n} )} \lesssim \frac{1}{(1+d(i,l) +d(j,k))^\nu}\,,
\]
where we used Assumptions~\ref{assum:A1} and  \ref{assum:A2} as well as Young's inequality to see 
that they still have a  polynomial decay of arbitrarily high order $\nu \in \N$. 
Thus, the volume growth condition \eqref{eq:growth_condition} for the metric $d$ implies \eqref{op to hs bound for cal T}.

Since the difference $\bs{\Delta} = \wt{\bs{M}}-\bs{M}$ satisfies the quadratic equation \eqref{Delta difference equation} with the error matrix $\bs{D}:=(\wt{\scr{S}}\wt{\bs{M}}-\scr{S}\wt{\bs{M}})\wt{\bs{M}}$ that satisfies the bound 
$\norm{\bs{D}} \lesssim \norm{\cal{R}}\norm{\wt{\bs{M}}}^2  \lesssim n^{-1/2 + 2\eps}/(1 + \eta^2)$ 
due to \eqref{op to hs bound for cal T} and the trivial bound $\norm{\wt{\bs{M}}} \le \frac{1}{\eta}$, 
   we use the invertibility of the stability operator $\scr{L}$ from\emph{ \ref{item:Bound on inverse cal L} }of  Lemma~\ref{lem:stability_outside}
to conclude \emph{\ref{item:tilde_M_M_operator_norm}} of the lemma in case $\eta \ge n^{-\eps}$. In case $\abs{\zeta}^2\ge \varrho(\cal{S})+\delta$ and $\eta \le 1$ the invertibility of $\cal{L}$ is still guaranteed by \emph{ \ref{item:Bound on inverse cal L} }of  Lemma~\ref{lem:stability_outside} and we have 
\[
\norm{\bs{\Delta}}\lesssim \norm{\cal{R}}\norm{\wt{\bs{M}}}^2+\norm{\bs{\Delta}}^2\lesssim n^{-1/2}+\norm{\bs{\Delta}}^2\,,
\]
where we used $\norm{\cal{R}}\lesssim n^{-1/2}$ by \eqref{op to hs bound for cal T} and $\norm{\wt{\bs{M}}}\le \norm{{\bs{M}}}+\norm{{\bs{\Delta}}}$.
Thus we can bootstrap the bound $\norm{\bs{\Delta}} \lesssim n^{-1/2}$ from the regime $\eta \ge 1$. 
 
 For the proof of \emph{\ref{item:tilde_M_M_hs_norm}}, we show the improved norm bound on $\cal{R}$ in the $\rm{hs}$-sense
 \bels{hs norm bound on cal R}{
 \norm{\cal{R}}_{\rm{hs}} \lesssim \frac{1}{n}\,.
 }
 To show \eqref{hs norm bound on cal R}, for each $R \in \C^{n \times n}$,  we estimate the $\rm{hs}$-norm through
\[
\norm{\cal{R}R}_{\rm{hs}}^2\le \frac{1}{n^3}\sum_{u,v,u',v'} \wh{K}_{uv,u'v'} \abs{R_{uv}R_{u'v'}}\lesssim 
\frac{1}{n^2} \norm{R}_{\rm{hs}}^2\,,
\]
where the second bound holds because for any $\nu \in \N$ the coefficients $\wh{K}_{uv,u'v'}$ satisfy
\[
\wh{K}_{uv,u'v'}:= \sum_{i,j}\absb{\rm{Cov} (\overline{X_{iu}}\sqrt{n} , X_{vj}\sqrt{n})\rm{Cov} (\overline{X_{iu'}}\sqrt{n} , X_{v'j}\sqrt{n})}\lesssim \frac{1}{(1+d(v,v') +d(u,u'))^\nu}\,.
\]

As above, from \eqref{hs norm bound on cal R}, we get $\norm{\bs{D}}_{\rm{hs}}\lesssim \norm{\cal{R}}_{\rm{hs}} \norm{\wt{\bs{M}}}^2 \lesssim n^{-1 + 2\eps}/(1 + \eta^2)$
and infer \emph{\ref{item:tilde_M_M_hs_norm}} of the lemma from \emph{ \ref{item:Bound on inverse cal L} }of  Lemma~\ref{lem:stability_outside}. 

Now we verify \emph{\ref{item:supp_tilde_rho}}. 
First we have $\dist(\ii \eta, \supp  {\rhoDOS}_\zeta )\gtrsim 1$ by Corollary~\ref{crl:scDOS inside disk}. 
We use the implication \cite[(v) implies (ii) in Lemma~D.1]{Altshape}. By \emph{\ref{item:tilde_M_M_operator_norm}} of Lemma~\ref{lmm:Properties of whM} the property \cite[(ii) in Lemma~D.1]{Altshape} is satisfied for $\wt{\bs{M}}$ whenever it is satisfied for $\bs{M}$ 
 due to their closeness. Finally, by the implication \cite[(ii) implies (v) in Lemma~D.1]{Altshape} we see that  property \cite[(v) in Lemma~D.1]{Altshape} holds for $\wt{\bs{M}}$, i.e.\ $\dist(\ii \eta, \supp  \wt{\rhoDOS}_\zeta )\gtrsim 1$. 
\end{Proof}

 The next proposition is a generalization of Proposition~\ref{pro:global_law_H_averaged}.  

\begin{proposition}[Global law for $\bs{H}_\zeta$, general version]  \label{pro:global_law_H_zeta} 
Let $X$ satisfy \ref{assum:A1} and \ref{assum:A2}. Then there is $C>0$ such that for  all $\varphi >0$ and   all sufficiently small $\delta >0$ we have 
\begin{subequations} \label{eq:global_law} 
\begin{align} 
\absb{\scalar{\bs{x}}{(\bs{G} - \bs{M})\bs{y}}} \leq \norm{\bs{x}} \norm{\bs{y}} \frac{n^{C \delta} }{(1 + \eta^2) \sqrt{n}} 
, \label{eq:global_law_isotropic} \\ 
\absb{\avg{ \bs{R}(\bs{G} - \bs{M})}} \leq \norm{\bs{R}} \frac{n^{C \delta}}{(1 + \eta^2) n} 
\label{eq:global_law_averaged_general}  
\end{align} 
\end{subequations} 
with very high probability uniformly 
for all $n \in\N$, $\zeta \in \mathbb{D}_{\varphi}$ and $\eta \in [n^{-\delta},n^{100}]$ as well as deterministic vectors $\bs{x}$, $\bs{y} \in \C^{2n}$ and deterministic matrices $\bs{R} \in \C^{2n\times 2n}$. 
Here $K$ is some absolute constant and the constant $C_\nu$  implicit in Definition~\ref{def:very_high_probability} of 'very high probability'  depends only on $\delta$ and $\varphi$ as well as the constants from \ref{assum:A1} and \ref{assum:A2}, in addition to~$\nu$. 
\end{proposition}

\begin{Proof}
The proposition is an immediate consequence of \cite[Theorem~2.1]{Erdos2017Correlated} since $\eta \ge n^{-\delta}$ means that the spectral parameter in the MDE is sufficiently far away from the self-consistent spectrum associated to $\bs{H}_\zeta$. 
As alluded to  after Proposition~\ref{prp:Bound on error matrix}, the self-energy in \cite{Erdos2017Correlated} is $\wt{\scr S}$ instead of $\scr S$. Consequently, the resolvent $\bs{G}$ is 
compared to $\wt{\bs{M}}$, the solution of \eqref{eq:MDE_tilde}, instead of $\bs{M}$.  
Thus \eqref{eq:global_law_isotropic} and \eqref{eq:global_law_averaged_general} follows from the closeness of $\wt{\bs{M}}$ to $\bs{M}$ from \emph{\ref{item:tilde_M_M_operator_norm}} and  \emph{\ref{item:tilde_M_M_hs_norm}}   in Lemma~\ref{lmm:Properties of whM}, respectively. 
\end{Proof}

\section{Quantitative law of large numbers} \label{app:law_large_numbers} 

In this section, we state a law of large numbers with an explicit rate of convergence for random variables with only $a$--moments for some $a >1$. 

\begin{proposition}[Quantitative law of large numbers] \label{pro:law_large_numbers} 
Let $m \in \N$. 
Let $(X_i)_{i=1}^m$ be centred i.i.d.\ random variables with $\E\1\abs{X_1}^{a} <\infty$ for some $a >1$. Then, for any $\delta \in (0,1]$, we have  
\[
\P\pbb{\absbb{\frac{1}{m}\sum_{i=1}^m X_i} \leq \bigg(\frac{10 \E \abs{X_1}^a}{ m^{a-1} \delta} \bigg)^{1/a}}  \geq 1 - \delta.  
\]
\end{proposition}

For the convenience of the reader, we provide a short proof of Proposition~\ref{pro:law_large_numbers}, 
 which is a quantitative variant of the standard proof of the law of large numbers. 

\begin{proof}
We set $\mu_a := \E \abs{X_1}^a$ and $\eps^a := \frac{10 \mu_a}{m^{a - 1} \delta}$. 
 We split into different terms and estimate
\begin{equation} \label{eq:decomposition_LLN} 
\P\pbb{\absbb{\frac{1}{m}\sum_{i=1}^m X_i} >\eps}  \le \P\pb{\abs{E_1}>\eps /3}+\P\pb{\abs{E_2}>\eps /3}+\bbm{1}\pb{\abs{E_3}>\eps /3}\,,
\end{equation}
where we introduced the random variables 
\[
E_1 :=\frac{1}{m}\sum_{i=1}^m (X_i -Y_{i}) \,, \qquad E_2:= \frac{1}{m}\sum_{i=1}^m  (Y_{i}-\E\1 Y_{i}) \,, \qquad E_3 :=\frac{1}{m}\sum_{i=1}^m  \E\1 Y_{i}\,, \qquad Y_{i}:=X_i\bbm{1}_{\abs{X_i}\le \eps\1m}\,.
\]
We now estimate the different terms in \eqref{eq:decomposition_LLN} separately. 
As a preparation, we conclude from Markov's inequality that 
\begin{equation} \label{eq:bound_X_1_markov} 
\P(\abs{X_1} \ge t )\le \frac{\mu_a}{t^a}
\end{equation} 
for any $t >0$. 
Hence, a simple union bound for the first term in \eqref{eq:decomposition_LLN} and $X_i = Y_i$ if $\abs{X_i} \leq \eps n$ by definition of $Y_i$ yield  
\[
\P(\abs{E_1} > \eps/3) \le \sum_{i=1}^m \P (\abs{X_i - Y_i} > \eps /3 ) \leq m \P(\abs{X_1}> \eps\1m)\le \frac{\mu_a}{\eps^a m^{a-1}}\,.
\]
The second term in \eqref{eq:decomposition_LLN} is bounded by Chebyshev's inequality using independence, i.e.\ by
\[
\P(\abs{E_2} > \eps/3) \le\frac{9\var(Y_1)}{\eps^2\1m} \le \frac{9\1\mu_a}{\eps^a\1m^{a-1}}\,,
\]
where in the last step we used that $\var(Y_1) \leq \E[\abs{X_1}^{2-a} \abs{X_1}^a \bbm{1}_{\abs{X_1}\le \eps m}] \leq (\eps m)^{2-a} \mu_a$.

Finally, since $\E X_1 = 0$, Hölder's inequality and \eqref{eq:bound_X_1_markov} imply 
\[
 \abs{E_3} \le  \abs{\E Y_1} = \abs{\E (Y_1 - X_1)} = \abs{\E X_1 \bbm{1}_{\abs{X_1} > \eps m}} \leq \frac{\mu_a}{\eps^{a-1} m^{a - 1}}. 
\]
Altogether we conclude 
\[
\P\pbb{\absbb{\frac{1}{m}\sum_{i=1}^m X_i} >\eps}  \le\frac{10\1\mu_a}{\eps^a\1m^{a-1}} + \bbm{1}\pbb{\frac{3\1\mu_a}{m^{a-1}} > \eps^a}\,,
\]
which completes the proof as the indicator function vanishes due to the definition of $\eps^a$ and $\delta \leq 1$. 
\end{proof}

 \section{Auxiliary results}  

\label{apx:Auxiliary results}

\begin{Proof}[Proof of Lemma~\ref{lmm:A2 implies C and D}]
We start by verifying  \cite[Assumption~(C)]{Erdos2017Correlated}. For the definition of the norms used inside this proof we refer to \cite{Erdos2017Correlated}.
To show $\tnorms[0]{\kappa}^{\rm{iso}}_2\lesssim 1$  we split the covariances $\kappa(\alpha,\beta):= \E w_\alpha w_\beta$ with double indices $\alpha=(a_1,a_2)$, $\beta=(b_1,b_2) \in \qq{2n}^2$ and  
\[
\bs{W}  := (w_\alpha)_{\alpha \in \qq{2n}^2}  := \bs{H}_\zeta+\bs{Z} = \mtwo{0 & X}{X^* & 0}\,,
\]
 into two summands $\kappa = \kappa_{\rm{d}}+\kappa_{\rm{c}}$ through 
 \[
 \kappa_{\rm{c}}(a_1a_2,b_1b_2):= \kappa(a_1a_2,b_1b_2)\bbm{1}((a_1, b_2) \in \qq{n}^2 \cup (n +\qq{n})^2 )\,.
 \]
We remark that with the definition of $\kappa_{\scr{R}}$ from \eqref{R and kappa R} and $\wt{\scr{S}}$ from \eqref{definition wt scr S} we get $\kappa=\kappa_{\wt{\scr{S}}}$, $\kappa_{\rm{d}}= \kappa_{\wt{\scr{S}}-\scr{S}}$ and  $\kappa_{\rm{c}}= \kappa_{\scr{S}}$.  Now we verify that $\tnorms[0]{\kappa_{{\#}}}_{{\#}}\lesssim 1$  for $\# = \rm{c}$, $\rm{d}$. If $\# = \rm{d}$ then we estimate 
\bes{
\tnorms[0]{\kappa_\rm{d}}_{\rm{d}} &= \sup_{\norm{\bs{x}}\le1}\textstyle \normbb{\pbb{\pB{\sum_{b_1} \absb{\sum_{a_1}x_{a_1}\kappa_\rm{d}(a_1a_2,b_1b_2)}^2}^{1/2}}_{a_2,b_2}}
\\
&  \lesssim \sup_{\norm{\bs{x}}\le1}\textstyle\normB{\pB{\frac{1}{(1+d(a_2,b_2))^\nu}}_{a_2,b_2}}\sum_{a_1,a_1'}\abs{x_{a_1}x_{a'_1}}\pB{\sum_{b_1}\frac{1}{(1+d(a_1,b_1))^\nu(1+d(a'_1,b_1))^\nu}}^{1/2}\lesssim 1\,,
}
where the norm on the right side of the equality refers to the standard operator norm of the matrix indexed by $a_2,b_2$ and  where  we used the decay of correlation from Assumption~\ref{assum:A2} via
\bels{kappa d bound}{
\abs{\kappa_\rm{d}(a_1a_2,b_1b_2)} \lesssim \frac{1}{(1+d(a_1,b_1))^\nu(1+d(a_2,b_2))^\nu}\,.
}
The case $\# = \rm{c}$ is seen by interchanging the roles of $b_1$ and $b_2$ and using
\bels{kappa c bound}{
\abs{\kappa_\rm{c}(a_1a_2,b_1b_2)} \lesssim \frac{1}{(1+d(a_1,b_2))^\nu(1+d(a_2,b_1))^\nu}\,.}
The bound
$\tnorms[0]{\kappa}^{\rm{av}}_2\le\tnorms[0]{\kappa_{\rm{d}}}^{\rm{av}}_2+\tnorms[0]{\kappa_\rm{c}}^{\rm{av}}_2 \lesssim 1$ also follows from \eqref{kappa d bound} and   \eqref{kappa c bound}.
 This implies $\tnorms[0]{\kappa}_2\lesssim 1$. 

The proof of $\tnorms[0]{\kappa}_k\lesssim n^\eps$ for $k\ge 3$ relies on  \cite[Lemma~A.1]{Erdos2017Correlated}. We demonstrate the strategy for these bounds for  $\tnorms[0]{\kappa}_3^{\rm{av}}$ and leave the other simpler cases to the reader. Writing the third order cumulant of three centred matrices $\bs{R}_1,\bs{R}_2,\bs{R}_3$  with $\bs{R}_i =(r^{(i)}_\alpha)_{\alpha \in \qq{2n}^{2}}$ as $\kappa_{\bs{R}_1\bs{R}_2\bs{R}_3}(\alpha,\beta,\gamma):=\E r^{(1)}_\alpha r^{(2)}_\beta r^{(3)}_\gamma$ we split $\kappa := \kappa_{\bs{W}\bs{W}\bs{W}}$
 into four summands $\kappa= \kappa_{\rm{dd}}+\kappa_{\rm{cc}}+\kappa_{\rm{dc}}+\kappa_{\rm{cd}}$. This split is performed by plugging in
\[
\bs{W}=\bs{X} + \bs{X}^*\,, \qquad \bs{X} = \mtwo{0 & X}{0 & 0}
\]
for each of the three $\bs{W}$-factors in the definition of $\kappa$, multiplying out and then grouping the summands according to
\bes{
& \kappa_{\rm{dd}}:= \kappa_{\bs{X}\bs{X}\bs{X}}+\kappa_{\bs{X}^*\bs{X}^*\bs{X}^*}\,, \qquad 
  \kappa_{\rm{cc}}:= \kappa_{\bs{X}\bs{X}^*\bs{X}}+\kappa_{\bs{X}^*\bs{X}\bs{X}^*}\,, 
  \\ 
 &\kappa_{\rm{dc}}:= \kappa_{\bs{X}\bs{X}\bs{X}^*}+\kappa_{\bs{X}^*\bs{X}^*\bs{X}}\,, \qquad 
    \kappa_{\rm{cd}}:= \kappa_{\bs{X}\bs{X}^*\bs{X}^*}+\kappa_{\bs{X}^*\bs{X}\bs{X}}\,.
}
Since all cases $\tnorms[0]{\kappa_{\rm{\#_1\#_2}}}_{\#_1\#_2}\lesssim n^\eps$  with $\#_i = \rm{d},\rm{c}$ are proven similarly by simply interchanging the role of  certain indices, we only show the case $\#_1=\#_2 = \rm{d}$. Due to \cite[Lemma~A.1]{Erdos2017Correlated} and Assumption~\ref{assum:A2} we have for any fixed $\nu \in \N$ and $\eps>0$  that
\[
\kappa_{\rm{dd}}(\alpha,\beta,\gamma)\lesssim n^{-\nu}\quad \text{whenever }\quad d\times d(\alpha,\beta) + d\times d(\alpha,\gamma)+ d\times d(\beta,\gamma) \ge n^\eps\,.
\] 
Thus using \eqref{eq:growth_condition} we conclude
\bes{
 \tnorms[0]{\kappa_{\rm{dd}}}_{\rm{dd}}^2 &= \textstyle\frac{1}{n^2}\sum_{b_2,c_1}\pb{\sum_{b_1,c_2}\sum_{a_1,a_2}\kappa_{\rm{dd}}(a_1a_2,b_1b_2,c_1c_2)}^2
\\
&\lesssim 
\max_{b_2,c_1}\abs{\{(a_1,a_2,b_1,c_2): d(c_1,a_1)+d(c_1,b_1)+d(b_2,a_2)+d(b_2,c_2) \le n^\eps\}}^2+\,n^{-\nu}\lesssim n^{C\2\eps}.
}

We proceed by verifying the modification of \cite[Assumption~(D)]{Erdos2017Correlated} described in the lemma, where the constant $\mu>0$ from the formulation of the assumption can be chosen arbitrarily.  With the choice of  nested neighbourhoods  $\cal{N}_k(\alpha):=\{\beta : d\times d (\alpha,\beta )\le k n^{(1-3\mu)/4p}\}$, where $p\in \N$ is from  \eqref{eq:growth_condition},  \cite[Assumption~(D)]{Erdos2017Correlated} is satisfied. Indeed, with the functions $f,g_1, \dots, g_q$ from the formulation of the assumptions we have
\bels{estimate on kappa f gi}{
\kappa(f,g_1, \dots, g_q)\lesssim n^{-\nu}\norm{f}_2\norm{g_1}_{2q} \dots \norm{g_q}_{2q}\,,
}
for any $\nu \in \N$. To see \eqref{estimate on kappa f gi} we follow the proof of  \cite[Lemma~A.1]{Erdos2017Correlated} with the choice $\ul{w}_A=(f)$ and $\ul{w}_B = (g_1, \dots, g_q)$. The covariance term in the last equation of the proof we estimate using \eqref{Decay of correlation} with $f_1= \Pi \ul{w}_{\cal{P}_i \cap A}=f$ and $f_2= \Pi \ul{w}_{\cal{P}_i \cap B}$. Since $d \times d(\supp f_1,\supp f_2) \ge n^{(1-3\mu)/4p}$  we get \eqref{estimate on kappa f gi} after applying H\"older inequality to $\norm{f_2}_2$ on the right hand side of \eqref{Decay of correlation}. 
\end{Proof}

\begin{Proof}[Proof of Lemma~\ref{lmm:Twist lemma}]
Let $w \in \C^d$ with $w \perp b$. To prove \eqref{Bound on A inverse on b perp} we use the spectral projection $P$ from \eqref{Definition of projection P} and its complementary projection $Q:=1-P$ as well as $\norm{a}_{\#}=1$ to estimate
\bels{starting estimate on Aw}{
\norm{Aw}_{\#}\,\ge\, \norm{AQw}_{\#}-\norm{APw}_{\#}\,\ge\, \norm{Qw}_{\#}-\abs{\alpha}\abs{\scalar{p}{w}}
\,.
}
Since $w$ is orthogonal to $b$ we have the identity
\bels{b orthogonal to w}{
0\,=\, \scalar{b}{w}\,=\, \scalar{b}{a}\scalar{p}{w}+ \scalar{b}{Qw}\,.
}
In particular, we find an upper bound on $\abs{\scalar{p}{w}}$ in terms of $\norm{Qw}_{\#}$, namely
\bels{first lower bound on Qw}{
\abs{\scalar{p}{w}}\,\le\, \frac{\abs{\scalar{b}{Qw}} }{\abs{\scalar{b}{a}}}\,\le\, \frac{1}{2\eps}\,\norm{Qw}_{\#}\,,
}
where we  used  the assumption from \eqref{assumptions on orthogonality for a and b}.
Continuing from \eqref{starting estimate on Aw} we see that 
\bels{lower bound on Aw in terms of Qw}{
\norm{Aw}_{\#}\,\ge\, \pbb{1-\frac{\abs{\alpha}}{2\eps}}\norm{Qw}_{\#}\,\ge\, \frac{1}{2}\,\norm{Qw}_{\#}\,,
}
because $\abs{\alpha}\le \eps$ by assumption. 

To finish the proof of \eqref{Bound on A inverse on b perp} we use
\bels{second lower bound on Qw}{
\norm{Qw}_{\#}\,\ge\,\norm{w}_{\#}-\norm{Pw}_{\#}\,\ge\, \norm{w}_{\#}- \abs{\scalar{p}{w}}\,.
}
Combining the two lower bounds \eqref{first lower bound on Qw} and \eqref{second lower bound on Qw} on $\norm{Qw}_{\#}$
and optimizing over the values of $\abs{\scalar{p}{w}}$ while using that $\eps\le 1$ yields 
\[
\norm{Qw}_{\#}\,\ge\,\frac{2 \eps}{3}\,\norm{w}_{\#}\,.
\]
Together with \eqref{lower bound on Aw in terms of Qw} this finishes the proof of Lemma~\ref{lmm:Twist lemma}.
\end{Proof}

\begin{lemma}[Quantitative implicit function theorem]
\label{lmm:Implicit function theorem}
Let ${T}: \C^{A}\times \C^{D} \to\C^{A}$ be a continuously differentiable function with invertible derivative $\nabla^{(1)}T(0,0)$ at the origin with respect to the first argument and ${T}(0,0)=0$. Suppose $\C^{A}$ and $\C^{D}$ are equipped with norms that we both denote by $\norm{\2\cdot\2}$ and let the linear operators on these spaces be equipped with the corresponding induced operator norms.
Let $\delta>0$ such that 
\bels{IFT 1-DD bound}{
\sup_{\quad(a,d\1)\1\in\1 B^A_\delta \times B^D_\delta}
\normb{\2\rm{Id}_{\C^{A}}-(\nabla^{(1)}T(0,0))^{-1}\nabla^{(1)}T(a,d\1)}\,\le\, \frac{1}{2}\,,
}
where $B_\delta^\#$ is the $\delta$-ball around $0$ with respect to $\norm{\2\cdot\2}$ in $\C^\#$. 
Suppose 
\[
\norm{(\nabla^{(1)}T(0,0))^{-1}}
\,\le\, C_1
\,,\qquad 
\sup_{\quad(a,d\1)\1\in\1 B^A_\delta \times B^D_\delta}
\norm{\nabla^{(2)}T(a,d\1)}
\,\le\, C_2
\,,
\]
for some positive constants $C_1,C_2$,  where $\nabla^{(2)}$ is the derivative with respect to the second variable.   Then there is a constant $\eps>0$, depending only on $\delta$, $C_1$ and $C_2$, and a unique function $f : B^D_\eps \to B^A_\delta$ such that $T(f(d),d\1)=0$ for all $d \in B^D_\eps$. The function $f$ is continuously differentiable. If $T$ is analytic, then so is $f$. 
\end{lemma}

\begin{lemma} 
\label{lmm:TF Stability lemma}
Let $\norm{a}$ denote the Euclidean norm of a vector $a \in \C^d$ and $\norm{A}$ the induced operator norm for a matrix $A \in \C^{d \times d}$.
Fix $\eps_1,\eps_2,\eps_3\in(0,1)$ with $100\eps_3 \le \eps_1\eps_2^2$. Let $F,T \in \C^{d \times d}$ be self-adjoint matrices such that $\norm{T}\le 1$ and  
\bels{spectrum of F}{
\spec(F) \subseteq \{-1\}\cup[-1+\eps_1,1-\eps_1]\cup\{1\}\,,
}
where $\pm 1$ are non-degenerate eigenvalues of $F$ with corresponding normalized eigenvectors $f_{\pm}$, i.e.\ 
$Ff_{\pm}=\pm f_\pm$. Suppose that $\norm{Tf_+}\le 1-\eps_2$ and $\norm{(1+T)f_-}\le\eps_3$. Then the resolvent of $TF$ satisfies 
\bels{resolvent bound for TF}{
\sup
\cbb{ \norm{(TF-\zeta)^{-1}}\,:\; \zeta \in \C\,,\; \zeta \not \in 
 \D_{1-6\1\eps_3}\cup (1+\D_{3\eps_3})
}\,\le\, \frac{4}{\eps_3}\,.
}
Furthermore, there is a single eigenvalue $\zeta_0$  close to $1$ and this eigenvalue is non-degenerate, more precisely, 
\bels{Nondegeneracy of TF eigenvalue}{
\spec(TF)\cap (1+\D_{3\1\eps_3})\,=\, \{\zeta_0\}\,,\qquad \dim \rm{ker} (TF-\zeta_0)^2\,=\, 1\,.
}
\end{lemma}

\begin{Proof}
First we realize that $f_-$ satisfies approximate eigenvalue equations for both $TF$ and $FT$, namely 
\bels{almost eigenequation for f-}{
\norm{(1-TF)f_-}\,\le\, \eps_3\,,\qquad \norm{(1-FT)f_-}\,=\,\norm{F(1+T)f_-}\,\le\,\eps_3\,.
}
 We now prove that when  restricted to the orthogonal complement of $f_-$, the matrix $TF$ is strictly smaller than $1$. More precisely,  we will establish that   
\bels{TF le 1 on f_- perp}{
\norm{TFa}\,\le\,\pbb{1-\frac{\eps_1\eps_2^2}{8}}\norm{a}\,,\qquad a \perp f_-\,.
}
To show \eqref{TF le 1 on f_- perp} we fix a  unit  vector $a \in f_-^\perp$,  $\norm{a} = 1$,  and decompose it according to  $f_+$ and its orthogonal complement,
\[
a\,=\, \sqrt{1-\alpha^2} \2f_+ +\alpha\2\wt{a}\,,\qquad \wt{a}\perp f_+\,,\; \norm{\wt{a}} \,=\, 1\,,
\]
for some $\alpha\in [0,1]$. Because $\norm{T}\le 1$ and $F$ has a spectral gap (cf.\ \eqref{spectrum of F}) we see that $\norm{TFa}$ is bounded from above by
\bels{first bound on TFa}{
\norm{TFa}\,\le\, \norm{Fa}\,\le\, \sqrt{1-\alpha^2+\alpha^2(1-\eps_1)^2}\,\le\, 1- \frac{\eps_1\alpha^2}{2}\,.
}
On the other hand, by using the assumption $\norm{Tf_+}\le 1-\eps_2$ we also get a second bound,
\bels{second bound on TFa}{
\norm{TFa}\,\le\,  \sqrt{1-\alpha^2}\,\norm{Tf_+}+\alpha \2\norm{F\wt{a}}\,\le\,
(1-\eps_2)\sqrt{1-\alpha^2}
+(1-\eps_1)\alpha\,\le \,1-\eps_2+(1-\eps_1)\alpha\,.
}
For $\alpha \le \eps_2/2$ we use \eqref{second bound on TFa} while for $\alpha \ge \eps_2 /2$ we use  \eqref{first bound on TFa} to infer \eqref{TF le 1 on f_- perp}.

With the help of \eqref{almost eigenequation for f-} and \eqref{TF le 1 on f_- perp} we represent $TF$ with respect to $f_-$ and an orthonormal basis of $f_-^\perp$. Thus we see that there is a unitary matrix $U \in \C^{d \times d}$ as well as $\alpha \in \C$, $b,a \in \C^{d-1}$ and $B \in \C^{(d-1)\times (d-1)}$ such that $A:=U^*TFU$ has the structure
\bels{unitary equivalent to TF}{
A\,=\, \mtwo{1+\eps_3\1\alpha & \eps_3\1 b^*}{\eps_3 \1a &B}\,,\qquad \abs{\alpha}\,\le\, 1\,,\; \norm{b}\,\le\, 1 \,,\; \norm{a}\,\le\, 1\,,\;\norm{B}\le 1-2\eps_4\,, \;\eps_4 :=\frac{\eps_1\eps_2^2}{16}\,.
}

Therefore it suffices to prove the resolvent bound \eqref{resolvent bound for TF} for any matrix $A$ with the structure \eqref{unitary equivalent to TF} in place of $TF$. For this purpose we fix a spectral parameter $\zeta$ with
\bels{bound on spectral parameter for A}{
\abs{\zeta}\2\ge\2 1-\eps_4\,,\qquad \abs{1-\zeta}\,\ge\, 2\eps_5\,, \qquad \eps_5 := \eps_3 +\frac{\eps_3^2}{\eps_4},
}
and use the Schur complement formula for $A-\zeta$ with respect to the block structure \eqref{unitary equivalent to TF}, i.e.\  we write
\[
(A-\zeta)^{-1}\,=\, \mtwo{((A-\zeta)^{-1})_{11}& ((A-\zeta)^{-1})_{1\perp}}{((A-\zeta)^{-1})_{\perp 1}& ((A-\zeta)^{-1})_{\perp\perp}}\,,
\] 
where $\perp$ refers to the component in the orthogonal complement of the first canonical basis vector $e_1$ of $\C^d$. 
The Schur complement itself is
\bels{Schur complement for A}{
 s_A(\zeta)\, :=\, 1+\eps_3\1\alpha-\zeta -\eps_3^2\2 b^* \1(B-\zeta)^{-1}a\,,
}
and because of \eqref{bound on spectral parameter for A} and the bound on $B$ from \eqref{unitary equivalent to TF}, we find
\[
 \abs{1-\zeta-s_A(\zeta)} \,\le\, \eps_5 \,.
\]
We conclude that $\abs{((A-\zeta)^{-1})_{11}} \le \eps_5^{-1}\le \eps_3^{-1}$ and also 
\[
\norm{((A-\zeta)^{-1})_{1\perp}}\,\le\, \frac{\eps_3}{\eps_4\eps_5}\le \frac{1}{\eps_3}\,,\quad \norm{((A-\zeta)^{-1})_{\perp 1}}\,\le\,  \frac{\eps_3}{\eps_4\eps_5}\le \frac{1}{\eps_3}\,,\quad
\norm{((A-\zeta)^{-1})_{\perp\perp}}\,\le\,\frac{1}{\eps_4}+\frac{\eps_3^2}{\eps_4^2\eps_5}\le \frac{2}{\eps_4}\,,
\]
which implies \eqref{resolvent bound for TF}, since $2\eps_4^{-1} \le \eps_3^{-1}$ and  $\zeta \not \in \D_{1-\eps_4}\cup\D_{2 \eps_5} \subseteq \D_{1-6\eps_3}\cup\D_{3 \eps_3} $ due to $\eps_4 \ge 6 \eps_3 \ge 4 \eps_5$.

To show \eqref{Nondegeneracy of TF eigenvalue} we use a simple interpolation argument. Consider the family of matrices
\[
A_\omega\,:=\, e_1\1e_1^* + \omega(A-e_1\1e_1^*)\,,\qquad  \omega \in [0,1]\,,
\]
interpolating between  $A_0=e_1\1e_1^*$ and $A_1=A$.  Since every element of this family has the same block structure \eqref{unitary equivalent to TF} as $A$, we conclude that \eqref{resolvent bound for TF} holds with $TF$ replaced by $A_\omega$. Since the eigenvalues of $A_\omega$ (as the $d$ zeros of the characteristic polynomial counted with multiplicity) depend continuously on $\omega$ and they cannot enter the regime in which the resolvent of $A_\omega$ is bounded, we conclude that the number of eigenvalues for $A_1=A$ within $1+\D_{3\1\eps_3}$ is that same as for $A_0$. Thus \eqref{Nondegeneracy of TF eigenvalue} is proven.
\end{Proof}

\begin{lemma}[Resolvent control for $\cal{S}$] 
\label{lmm:Resolvent control for cal S}
Let $\cal{S}: \C^{n \times n}\to \C^{n \times n}$ be a positivity preserving operator such that $\spradius(\cal{S})=1$ and $c \avg{A} \le \cal{S}A \le C \avg{A}$  for any $A \in \scr{C}_+$. Then $\cal{S}$ satisfies the resolvent control
\bels{resolvent bound for cal S}{
\sup
\cB{ \norm{(\cal{S}-\xi)^{-1}}_{\#}\,:\; \xi \in \C\,,\; \xi \not \in 
 \D_{1-2\eps}\cup (1+\D_\eps)
}\,\lesssim_\eps\,1 \,,
}
for any sufficiently small $\eps> 0$ (depending on the constants $c$ and $C$) and $\# = \rm{hs}, \norm{\1\cdot\1}$.  The algebraic multiplicity of the eigenvalue $\spradius(\cal{S})=1$ is one and the corresponding left and right Perron Frobenius eigenvectors,  $S_1$ and $S_2$,   satisfy
\bels{Si comparable with 1}{
S_1 \sim 1\,, \qquad S_2 \sim 1\,,
}
where $\cal{S}^* S_1=S_1$ and $\cal{S} S_2=S_2$ as well as  $\avg{S_1}=\avg{S_2}=1$.   
\end{lemma}
\begin{Proof}
We start the proof for $\#=\rm{hs}$.
 We denote by $S_1 $ and $S_2$ the positive definite left and right Perron Frobenius eigenvectors of $\cal{S}$ with normalisation $\avg{S_1}=\avg{S_2}=1$. 
The assumption $\cal{S}A\sim \avg{A}$  immediately implies \eqref{Si comparable with 1} and also the statement  about the  multiplicity. 
 Instead of studying $\cal{S}$ we study $\Sigma:\C^{n \times n}\times \C^{n \times n}\to \C^{n \times n}\times \C^{n \times n}$ defined as
\bels{eq:definition_Sigma}{
{\Sigma}\,:=\,  
\mtwo{\cal{S}^* & 0}
{ 0 &\cal{S}}
\,=\, {\cal{V}}^{-1}\cal{T} \cal{F}{\cal{V}}\,,
} 
where the representation on the very right  is in terms of  the invertible operator  
\[
\cal{V}\,:=\,
\mtwo{
\cal{C}_{\!\wt{K}_2}\,\cal{C}_{\!\sqrt{S_1}}^{-1}
&
0
}
{
0
&
\cal{C}_{\!\wt{K}_1}\,\cal{C}_{\!\sqrt{S_2}}^{-1}
}\,,
\] 
and the two self-adjoint operators
\bels{}{
\cal{T}\,:=\,
\mtwo{
0
&
 \cal{C}_{\!\wt{K}_2}\,\cal{K}_{\wt{P}^*}\2\cal{C}_{{\wt{K}_1}}
}
{
 \cal{C}_{\! \wt{K}_1} \,\cal{K}_{\wt{P}}\2\cal{C}_{\wt{K}_2}
&
0
},
\quad 
\cal{F}\,:=\,
\mtwo{
0
&
\cal{C}_{\!\wt{K}_2}^{-1}\,\cal{C}_{\!\sqrt{S_1}}\2\cal{S}\2\cal{C}_{\!\sqrt{S_2}}\2\cal{C}_{\!\wt{K}_1}^{-1}
}
{
\cal{C}_{\!\wt{K}_1}^{-1}\,\cal{C}_{\!\sqrt{S_2}}\2\cal{S}^*\2\cal{C}_{\!\sqrt{S_1}}\2\cal{C}_{\!\wt{K}_2}^{-1}
&
0
}.
}
Here we introduced a short hand notation for the matrices 
\[
\wt{K}_1:=(\sqrt{S_2}S_1\sqrt{S_2})^{1/4}\,, \quad \wt{K}_2=(\sqrt{S_1}S_2\sqrt{S_1})^{1/4},\quad 
\wt{P} :=\frac{1}{\sqrt{S_2}\sqrt{S_1}}\,.
 \]
 Note that the definitions of $\cal{V},\cal{F}$ and $\cal{T}$ above are compatible with \eqref{definition of cal V}, \eqref{definition of cal F} and \eqref{definition of cal T} in the limit $\tau \to 1, \eta \downarrow 0$, while with the same limit we have $\wt{K}_i:= \lim K_i/\sqrt{\avg{V_1}}$ and $\wt{P}:= \lim P \avg{V_1}$.

 Since $\Sigma$ from \eqref{eq:definition_Sigma} is a direct sum of $\cal{S}$ and $\cal{S}^*$, the claim \eqref{resolvent bound for cal S} is equivalent to the same statement with $\cal{S}$ replaced by $\Sigma$. 
Owing to \eqref{Si comparable with 1} we have $\norm{\cal{V}}_{\rm{hs}}\norm{\cal{V}^{-1}}_{\rm{hs}}\sim \norm{\cal{V}}\norm{\cal{V}^{-1}}\sim 1$.  
Therefore, \eqref{resolvent bound for cal S} for $\Sigma$ now follows  from 
 the following facts about $\cal{T}$ and $\cal{F}$:  
\[
 \cal{T}\cal{V}S_\pm=\pm \cal{V}S_\pm\,, \qquad \cal{F}\cal{V}S_\pm=\pm \cal{V}S_\pm\,, \qquad 
 \norm{\cal{T}}_{\rm{hs}} \leq 1,  \qquad \norm{\cal{F}|_{(\cal{V}S_+)^\perp \cap (\cal{V}S_-)^\perp}}_{\rm{hs}} \le 1-2\eps\,,
\]
for some $\eps \sim 1$, where $S_\pm = (S_1, \pm S_2)$. 
Here, the last bound is obtained from \cite[Lemma~4.8]{AjankiCorrelated}
in the same way as \eqref{spectral gap of cal F} in Lemma~\ref{lmm:Properties of cal F} was obtained.
Since via Lemma~\ref{lmm:Smoothing lemma} we can lift the resolvent control to the other norm $\# = \norm{\1\cdot\1}$, 
this finishes the proof of the lemma. 
\end{Proof}

\begin{corollary}[Resolvent control for edge stability operator] \label{crl:Resolvent control for scr J0 L0} 
Let $\cal{S}$ and $S_1,S_2$ be as in Lemma~\ref{lmm:Resolvent control for cal S} and define the operator $\cal{K}: \C^{n \times n} \times \C^{n \times n} \to \C^{n \times n} \times \C^{n \times n}$ via
\[
\cal{K}:= \mtwo{S_2(1-\cal{S}^*)& 0}{0&S_1(1-\cal{S})}\,.
\]
Then this operator satisfies 
 the resolvent control
\bels{resolvent bound for scr J0L0}{
\sup
\cB{ \norm{(\cal{K}-\xi)^{-1}}_{\#}\,:\; \xi  \in 
\D_{2\eps}\setminus \D_\eps
}\,\lesssim_\eps\,1 \,,
}
for any sufficiently small $\eps>0$ and $\# = \rm{hs}, \norm{\1\cdot\1}$. 
Furthermore, the eigenvalue $0$ has algebraic and geometric multiplicity equal to $2$ and corresponding right and left eigenvectors
\bels{eigenvector for scr L 0}{
\cal{K}S_\pm\,=\, {0}\,, \qquad   \cal{K}^*E_\pm\,=\, {0} \,, \qquad 
S_\pm \,:=\,\vtwo{S_1 }{\pm S_2}
\,,
}
where $E_\pm$ are defined in \eqref{definition E pm}.
\end{corollary}
\begin{Proof} 
Since the operator $\cal{K}$ separately acts on the first and second component of a pair of matrices,  the assertions about the multiplicity of 0 and \eqref{eigenvector for scr L 0} 
follow from Lemma~\ref{lmm:Resolvent control for cal S} and a simple computation.  
Similarly,  it suffices to prove the resolvent control \eqref{resolvent bound for scr J0L0} for each component, i.e.\  to show it for $S_1(1-\cal{S})$ and $S_2 (1-\cal{S}^*)$. We will only consider the first since the latter is treated similarly with the roles on $\cal{S}$ and $\cal{S}^*$ interchanged.  We define the projections
\[
\cal{P}A\,:=\, \frac{\scalar{S_1}{A}}{\scalar{S_1}{S_2}}S_2\,, \quad \wt{\cal{P}}A\,:=\,\frac{\scalar{1}{A}}{\scalar{1}{S_2}}S_2
\,, \quad \cal{P}^\perp A\,:=\,\frac{\scalar{S_2}{A}}{\scalar{S_2}{S_2}}S_2
\,.
\]
and their complements
\[
\cal{Q}:=1-\cal{P}\,, \qquad \wt{\cal{Q}}:=1-\wt{\cal{P}}\,, \qquad \cal{Q}^\perp:=1-\cal{P}^\perp.
\]
Due to Lemma~\ref{lmm:Resolvent control for cal S} the rank one projections $\cal{P}$ and $\wt{\cal{P}}$ are the spectral projections associated to the non-degenerate eigenvalue $0$ of $1- \cal{S}$ and $S_1(1- \cal{S})$, respectively. The claim follows now because the operator $S_1(1- \cal{S})$ has a bounded inverse on the image of $\wt{\cal{Q}}$, i.e.\  
\bels{invertibility on range wt Q}{
\norm{S_1(1- \cal{S})A}_{\rm{hs}}=\norm{S_1(1- \cal{S})\cal{Q}A}_{\rm{hs}}\gtrsim  \norm{\cal{Q}A}_{\rm{hs}}\sim \norm{\wt{\cal{Q}}A}_{\rm{hs}}\,,
}
where for the inequality we used $S_1 \sim 1$ and Lemma~\ref{lmm:Resolvent control for cal S} and for the last relation
\[
 \norm{\cal{Q}A}_{\rm{hs}}^2=\norm{\cal{Q}\wh{A}}_{\rm{hs}}^2= \norm{\wh{A}}^2_{\rm{hs}} +\frac{\abs{\scalar{S_1}{\wh{A}}}^2}{\scalar{S_1}{S_2}^2}\sim \norm{\wh{A}}_{\rm{hs}}^2\sim  \norm{\wt{\cal{Q}}\wh{A}}_{\rm{hs}}^2 =  \norm{\wt{\cal{Q}}A}_{\rm{hs}}^2.
\]
Here we used the short hand $\wh{A} =  \cal{Q}^\perp A$,  $\cal{Q}A = \cal{Q} \wh{A}$, $\wt{\cal{Q}}\wh{A} = \wt{\cal{Q}}A$  and the second comparison relation holds for the same reason  as  the first. This finishes the proof of \eqref{resolvent bound for scr J0L0} for $\#=\rm{hs}$. For $\#=\norm{\1\cdot\1}$ we use Lemma~\ref{lmm:Smoothing lemma}.
\end{Proof}

\begin{lemma} \label{lem:symmetry_inverse_block_matrix} 
Let $X \in \C^{n\times n}$ be an arbitrary matrix. Then, for all $z \in \C\setminus \R$, we have
\[ \begin{pmatrix} z & X \\ X^* & z \end{pmatrix}^{-1} \in \bs{E}_-^\perp. \] 
\end{lemma} 

\begin{Proof} Schur's complement formula directly implies that  
\begin{equation} 
\label{Schur for H}
\begin{pmatrix} z & X \\ X^* & z \end{pmatrix}^{-1} = \begin{pmatrix} z(z^2 - XX^*)^{-1}&  - (z^2 - XX^*)^{-1} X \\ - X^* (z^2 - XX^*)^{-1} & 
z(z^2 - X^* X)^{-1} \end{pmatrix}.  \end{equation}
As $XX^*$ and $X^*X$ have the same eigenvalues and their multiplicities coincide, this proves the lemma.
\end{Proof}

{\small 
\providecommand{\bysame}{\leavevmode\hbox to3em{\hrulefill}\thinspace}
\providecommand{\MR}{\relax\ifhmode\unskip\space\fi MR }
\providecommand{\MRhref}[2]{%
  \href{http://www.ams.org/mathscinet-getitem?mr=#1}{#2}
}
\providecommand{\href}[2]{#2}

} 


\begin{thebibliography}{10}

\bibitem{AjankiSingularities}
O.~H. Ajanki, L.~Erd\H{o}s, and T.~Kr\"{u}ger, \emph{Singularities of solutions
  to quadratic vector equations on the complex upper half-plane}, Comm. Pure
  Appl. Math. \textbf{70} (2017), no.~9, 1672--1705.

\bibitem{AjankiQVE}
O.~H. Ajanki, L.~Erd\H{o}s, and T.~Kr\"{u}ger, \emph{Quadratic vector equations
  on complex upper half-plane}, Mem. Am. Math. Soc. \textbf{261} (2019),
  no.~1261.

\bibitem{AjankiCorrelated}
O.~H. Ajanki, L.~Erd\H{o}s, and T.~Kr\"{u}ger, \emph{Stability of the matrix
  {D}yson equation and random matrices with correlations}, Probab. Theory
  Related Fields \textbf{173} (2019), no.~1-2, 293--373.

\bibitem{Allesina:2015ux}
S.~Allesina and S.~Tang, \emph{The stability--complexity relationship at age
  40: a random matrix perspective}, Population Ecology \textbf{57} (2015),
  no.~1, 63--75.

\bibitem{Altcirc}
J.~Alt, L.~Erd\H{o}s, and T.~Kr\"{u}ger, \emph{Local inhomogeneous circular
  law}, Ann. Appl. Probab. \textbf{28} (2018), no.~1, 148--203.

\bibitem{Altshape}
J.~Alt, L.~Erd\H{o}s, and T.~Kr\"{u}ger, \emph{The {D}yson equation with linear
  self-energy: spectral bands, edges and cusps}, Doc. Math. \textbf{25} (2020),
  1421--1539.

\bibitem{AltSpecRad}
J.~Alt, L.~Erd\H{o}s, and T.~Kr\"{u}ger, \emph{Spectral radius of random
  matrices with independent entries}, to appear in Probab. Math. Phys. (2021),
  arXiv:1907.13631.

\bibitem{Kronecker}
J.~Alt, L.~Erd\H{o}s, T.~Kr\"{u}ger, and Yu. Nemish, \emph{Location of the
  spectrum of {K}ronecker random matrices}, Ann. Inst. Henri Poincar\'{e}
  Probab. Stat. \textbf{55} (2019), no.~2, 661--696.

\bibitem{AZind}
G.~W. Anderson and O.~Zeitouni, \emph{{A CLT for a band matrix model}}, Probab.
  Theory Related Fields \textbf{134} (2005), no.~2, 283--338.

\bibitem{MR2417889}
G.~W. Anderson and O.~Zeitouni, \emph{A law of large numbers for finite-range
  dependent random matrices}, Comm. Pure Appl. Math. \textbf{61} (2008), no.~8,
  1118--1154.

\bibitem{ArmitageGardiner}
D.~H. Armitage and S.~J. Gardiner, \emph{Classical potential theory}, Springer
  Monographs in Mathematics, Springer-Verlag London, Ltd., London, 2001.

\bibitem{bai1997}
Z.~D. Bai, \emph{Circular law}, Ann. Probab. \textbf{25} (1997), no.~1,
  494--529.

\bibitem{MR3332852}
M.~Banna, F.~Merlev\`ede, and M.~Peligrad, \emph{On the limiting spectral
  distribution for a large class of symmetric random matrices with correlated
  entries}, Stochastic Process. Appl. \textbf{125} (2015), no.~7, 2700--2726.

\bibitem{BelinschiSniadySpeicher2018}
S.~T. Belinschi, P.~\'{S}niady, and R.~Speicher, \emph{Eigenvalues of
  non-{H}ermitian random matrices and {B}rown measure of non-normal operators:
  {H}ermitian reduction and linearization method}, Linear Algebra Appl.
  \textbf{537} (2018), 48--83.

\bibitem{BianeLehner2001}
P.~Biane and F.~Lehner, \emph{Computation of some examples of {B}rown's
  spectral measure in free probability}, Colloq. Math. \textbf{90} (2001),
  no.~2, 181--211.

\bibitem{BordenaveChafai2012}
C.~Bordenave and D.~Chafa\"{\i}, \emph{Around the circular law}, Probab. Surv.
  \textbf{9} (2012), 1--89.

\bibitem{Bourgade2014}
P.~Bourgade, H.-T. Yau, and J.~Yin, \emph{Local circular law for random
  matrices}, Probab. Theory Related Fields \textbf{159} (2014), no.~3-4,
  545--595.

\bibitem{BYY_circular2}
P.~Bourgade, H.-T. Yau, and J.~Yin, \emph{The local circular law {II}: the edge
  case}, Probab. Theory Related Fields \textbf{159} (2014), no.~3-4, 619--660.

\bibitem{MR1431189}
A.~Boutet~de Monvel, A.~Khorunzhy, and V.~Vasilchuk, \emph{Limiting eigenvalue
  distribution of random matrices with correlated entries}, Markov Process.
  Related Fields \textbf{2} (1996), no.~4, 607--636.

\bibitem{Brown1986}
L.~G. Brown, \emph{Lidski\u{\i}'s theorem in the type {II} case}, Geometric
  methods in operator algebras ({K}yoto, 1983), Pitman Res. Notes Math. Ser.,
  vol. 123, Longman Sci. Tech., Harlow, 1986, pp.~1--35.

\bibitem{CipolloniCLT}
G.~Cipolloni, L.~Erd\H{o}s, and D.~Schr\"{o}der, \emph{Central limit theorem
  for linear eigenvalue statistics of non-{H}ermitian random matrices},
  preprint (2019), arXiv:1912.04100.

\bibitem{CipolloniEdge}
G.~Cipolloni, L.~Erd\H{o}s, and D.~Schr\"{o}der, \emph{Edge universality for
  non-{H}ermitian random matrices}, preprint (2019), arXiv:1908.00969.

\bibitem{CipolloniCLTreal}
G.~Cipolloni, L.~Erd\H{o}s, and D.~Schr\"{o}der, \emph{Fluctuation around the
  circular law for random matrices with real entries}, preprint (2020),
  arXiv:2002.02438.

\bibitem{Cook2018}
N.~Cook, W.~Hachem, J.~Najim, and D.~Renfrew, \emph{Non-{H}ermitian random
  matrices with a variance profile ({I}): deterministic equivalents and
  limiting {ESD}s}, Electron. J. Probab. \textbf{23} (2018), Paper No. 110, 61.

\bibitem{DHKBrownMeasure}
B.~K. Driver, B.~C. Hall, and T.~Kemp, \emph{The {B}rown measure of the free
  multiplicative {B}rownian motion}, preprint (2019), arXiv:1903.11015.

\bibitem{Erdos2017Correlated}
L.~Erd\H{o}s, T.~Kr\"{u}ger, and D.~Schr\"{o}der, \emph{Random matrices with
  slow correlation decay}, Forum Math. Sigma \textbf{7} (2019), e8, 89.

\bibitem{EJP2473}
L.~Erd{\H o}s, A.~Knowles, H.-T. Yau, and J.~Yin, \emph{The local semicircle
  law for a general class of random matrices}, Elect. J. Probab. \textbf{18}
  (2013), no.~59, 1--58.

\bibitem{FugledeKadison1952}
B.~Fuglede and R.~V. Kadison, \emph{Determinant theory in finite factors}, Ann.
  of Math. (2) \textbf{55} (1952), 520--530.

\bibitem{Girko1985}
V.~L. Girko, \emph{Circular law}, Theory Probab. Appl. \textbf{29} (1985),
  no.~4, 694--706.

\bibitem{MR1887675}
V.~L. Girko, \emph{Theory of stochastic canonical equations. {V}ol. {I}},
  Mathematics and its Applications, vol. 535, Kluwer Academic Publishers,
  Dordrecht, 2001.

\bibitem{Gotze2017}
F.~G{\"o}tze, A.~A. Naumov, and A.~N. Tikhomirov, \emph{Local laws for
  non-hermitian random matrices}, Doklady Mathematics \textbf{96} (2017),
  no.~3, 558--560.

\bibitem{Guionnet-GaussBand}
A.~Guionnet, \emph{{Large deviations upper bounds and central limit theorems
  for non-commutative functionals of Gaussian large random matrices}}, Annales
  de l'IHP Probabilit{\'e}s et statistiques \textbf{38} (2002), 341--384.

\bibitem{GWZBrownMeasure}
A.~Guionnet, P.~Wood, and O.~Zeitouni, \emph{Convergence of the spectral
  measure of non normal matrices}, Proceedings of the American Mathematical
  Society \textbf{142} (2014), no.~2, 667--679.

\bibitem{HaagerupLarsen2000}
U.~Haagerup and F.~Larsen, \emph{Brown's spectral distribution measure for
  {$R$}-diagonal elements in finite von {N}eumann algebras}, J. Funct. Anal.
  \textbf{176} (2000), no.~2, 331--367.

\bibitem{MR2191967}
W.~Hachem, P.~Loubaton, and J.~Najim, \emph{The empirical eigenvalue
  distribution of a {G}ram matrix: from independence to stationarity}, Markov
  Process. Related Fields \textbf{11} (2005), no.~4, 629--648.

\bibitem{Helton01012007}
J.~W. Helton, R.~Rashidi~Far, and R.~Speicher, \emph{Operator-valued
  semicircular elements: solving a quadratic matrix equation with positivity
  constraints}, Int. Math. Res. Not. IMRN (2007), no.~22, Art. ID rnm086, 15.

\bibitem{Larsenthesis}
F.~Larsen, \emph{Brown measures and {R}-diagonal elements in finite von
  {N}eumann algebras}, Ph.D. thesis, University of Southern Denmark, 1999.

\bibitem{may1972will}
R.~M. May, \emph{Will a large complex system be stable?}, Nature \textbf{238}
  (1972), 413--414.

\bibitem{MingoSpeicher}
J.~A. Mingo and R.~Speicher, \emph{Free probability and random matrices},
  Fields Institute Monographs, vol.~35, Springer, New York; Fields Institute
  for Research in Mathematical Sciences, Toronto, ON, 2017.

\bibitem{NicaSpeicher1997}
A.~Nica and R.~Speicher, \emph{{$R$}-diagonal pairs---a common approach to
  {H}aar unitaries and circular elements}, Free probability theory ({W}aterloo,
  {ON}, 1995), Fields Inst. Commun., vol.~12, Amer. Math. Soc., Providence, RI,
  1997, pp.~149--188.

\bibitem{PasturShcerbinaAMSbook}
L.~A. Pastur and M.~Shcherbina, \emph{{Eigenvalue Distribution of Large Random
  Matrices}}, Mathematical Surveys and Monographs, vol. 171, Amer. Math. Soc.,
  2011.

\bibitem{PhysRevLett.97.188104}
K.~Rajan and L.~F. Abbott, \emph{Eigenvalue spectra of random matrices for
  neural networks}, Phys. Rev. Lett. \textbf{97} (2006), 188104.

\bibitem{MR2444540}
R.~Rashidi~Far, T.~Oraby, W.~Bryc, and R.~Speicher, \emph{On slow-fading {MIMO}
  systems with nonseparable correlation}, IEEE Trans. Inform. Theory
  \textbf{54} (2008), no.~2, 544--553.

\bibitem{RudelsonVershynin2008}
M.~Rudelson and R.~Vershynin, \emph{The {L}ittlewood-{O}fford problem and
  invertibility of random matrices}, Adv. Math. \textbf{218} (2008), no.~2,
  600--633.

\bibitem{rudelson2015}
M.~Rudelson and R.~Vershynin, \emph{Delocalization of eigenvectors of random
  matrices with independent entries}, Duke Math. J. \textbf{164} (2015),
  no.~13, 2507--2538.

\bibitem{MR2155229}
J.~H. Schenker and H.~Schulz-Baldes, \emph{Semicircle law and freeness for
  random matrices with symmetries or correlations}, Math. Res. Lett.
  \textbf{12} (2005), no.~4, 531--542.

\bibitem{ShlyakhtenkoGBM}
D.~Shlyakhtenko, \emph{{Random Gaussian band matrices and freeness with
  amalgamation}}, International Mathematics Research Notices (1996), no.~20,
  1013--1015.

\bibitem{Sompolinsky1988}
H.~Sompolinsky, A.~Crisanti, and H.-J. Sommers, \emph{Chaos in random neural
  networks}, Phys. Rev. Lett. \textbf{61} (1988), no.~3, 259--262.

\bibitem{Speicher1998}
R.~Speicher, \emph{Combinatorial theory of the free product with amalgamation
  and operator-valued free probability theory}, Mem. Amer. Math. Soc.
  \textbf{132} (1998), no.~627, x+88.

\bibitem{tao2015}
T.~Tao and V.~Vu, \emph{Random matrices: Universality of local spectral
  statistics of non-hermitian matrices}, Ann. Probab. \textbf{43} (2015),
  no.~2, 782--874.

\bibitem{tao2010}
T.~Tao, V.~Vu, and M.~Krishnapur, \emph{Random matrices: Universality of {ESD}s
  and the circular law}, Ann. Probab. \textbf{38} (2010), no.~5, 2023--2065.

\bibitem{Voiculescu1995}
D.~Voiculescu, \emph{Operations on certain non-commutative operator-valued
  random variables}, no. 232, 1995, Recent advances in operator algebras
  (Orl\'{e}ans, 1992), pp.~243--275. \MR{1372537}

\bibitem{Wigner1955}
E.~P. Wigner, \emph{Characteristic vectors of bordered matrices with infinite
  dimensions}, Ann. of Math. \textbf{62} (1955), no.~3, 548--564.

\bibitem{Y_circularlaw}
J.~Yin, \emph{The local circular law {III}: general case}, Probab. Theory
  Related Fields \textbf{160} (2014), no.~3-4, 679--732.

\end{thebibliography}
\end{document}